\documentclass[11pt,final]{amsart}
\usepackage{amsmath, amsthm, amsfonts, ifpdf}
\usepackage{color}
\usepackage{graphicx}
\usepackage{multicol}
\usepackage{epstopdf}
\usepackage{verbatim}

\usepackage[utf8x]{inputenc}

\usepackage{tikz}
\usetikzlibrary{arrows,shapes,trees,backgrounds}
\usetikzlibrary{intersections}

\theoremstyle{plain}
\newtheorem{theorem}{Theorem}[section]
\newtheorem*{theorem*}{Theorem}

\newtheorem{pro}[theorem]{Proposition}
\newtheorem{Def}[theorem]{Definition}
\newtheorem{lem}[theorem]{Lemma}

\newtheorem{cor}[theorem]{Corollary}

\theoremstyle{definition}
\newtheorem*{Def*}{Definition}

\newtheorem{Rem}[theorem]{Remark}

\numberwithin{equation}{section}

\newcommand{\bpo}{\begin{pro}}
\newcommand{\epo}{\end{pro}}
\newcommand{\be}{\begin{equation}}
\newcommand{\ene}{\end{equation}}
\newcommand{\br}{\begin{Rem}}
\newcommand{\er}{\end{Rem}}
\newcommand{\bl}{\begin{lem}}
\newcommand{\el}{\end{lem}}
\newcommand{\bd}{\begin{Def}}
\newcommand{\ed}{\end{Def}}
\newcommand{\ben}{\begin{enumerate}}
\newcommand{\een}{\end{enumerate}}
\newcommand{\bp}{\begin{proof}}
\newcommand{\ep}{\end{proof}}
\newcommand{\beq}{\begin{equation*}}
\newcommand{\eeq}{\end{equation*}}
\newcommand{\bear}{\begin{eqnarray*}}
\newcommand{\eear}{\end{eqnarray*}}
\newcommand{\bt}{\begin{theorem}}
\newcommand{\et}{\end{theorem}}
\newcommand{\bst}{\begin{split}}
\newcommand{\est}{\end{split}}

\newcommand{\bal}{\begin{aligned}}
\newcommand{\eal}{\end{aligned}}

\renewcommand{\P}{\partial}
\newcommand{\F}[2]{\frac{#1}{#2}}
\newcommand{\la}{\langle}
\newcommand{\ra}{\rangle}
\newcommand{\mT}{\mathcal{T}}

\newcommand{\R}{\mathbb{R}}

\newcommand{\bnb}{\bar{\nabla}}
\newcommand{\nb}{\nabla}

\newcommand{\RM}{Riemannian manifold}

\newcommand{\wrt}{with respect to}
\newcommand{\Sc}{\varepsilon}

\newcommand{\Ta}{\Theta}

\newcommand{\vp}{\varphi}

\newcommand{\PLH}{{\mkern-1mu\times\mkern-1mu}}

\newcommand{\DP}{Dirichlet problem}
\newcommand{\Wg}{\omega}
\newcommand{\cf}{\lambda}
\newcommand{\AF}{\mathfrak{F}}
\newcommand{\Bc}{Besicovitch}
\newcommand{\Bh}{\tilde{h}}
\newcommand{\mt}{\mu_{\alpha}}
\newcommand{\Ca}{Caccioppoli}

\newcommand{\Lp}{Lipschitz}

\newcommand{\vf}{vol}
\newcommand{\mC}{\mathcal{C}}
\newcommand{\mW}{\mathcal{W}}
\newcommand{\ap}{a priori estimate}
\newcommand{\tg}{\tilde{g}}
\newcommand{\mg}{\mathcal{GT}}
\newcommand{\fc}{characteristic function}
\newcommand{\ms}{\mathcal{S}}

\def\XXint#1#2#3{{\setbox0=\hbox{$#1{#2#3}{\int}$}
    \vcenter{\hbox{$#2#3$}}\kern-.5\wd0}}

\makeatletter
\def\@citestyle{\m@th\upshape\mdseries}
\def\citeform#1{{\bfseries#1}}
\def\@cite#1#2{{%
  \@citestyle[\citeform{#1}\if@tempswa, #2\fi]}}
\@ifundefined{cite }{%
  \expandafter\let\csname cite \endcsname\cite
  \edef\cite{\@nx\protect\@xp\@nx\csname cite \endcsname}%
}{}
\makeatother
	
\begin{document}
		\title[Dirichlet problems of translating equations]{Generalized solutions to the Dirichlet problem of translating mean curvature equations}
		\author{Hengyu Zhou}
		\date{\today}
		\address{College of Mathematics and Statistics, Chongqing University, No. 55, DaxueCheng South Rd., Shapingba, Chongqing, 401331, P. R. China}
		\email{zhouhyu@cqu.edu.cn}
		\date{\today}
		\begin{abstract}
			 In this paper we study the Dirichlet problem of translating mean curvature equations over domains in {\RM}s with dimension $n$. Imitating the generalized solution theory of Miranda-Giusti, we define a new conformal area functional and a generalized solution to this Dirichlet problem.  The existence of generalized solutions to this problem on bounded Lipschitz domains is established. If the domain is mean convex and bounded with $C^2$ boundary, its closure does not contain any closed minimal hypersurface except a singular set with its Hausdorff dimension at most $n-7$ and the boundary data is continuous,  the generalized solution is the desirable classical smooth solution. The non-minimal condition could not be removed in general. 
			\end{abstract}
		\maketitle
		\section{Introduction}
	 In this paper we consider the {\DP} of the following quasilinear equation
			\be\label{main:equation}
			\left\{\begin{split}
		div(\F{Du}{\omega}) &=\F{\alpha}{\omega},\quad x\in \Omega,\quad \omega =\sqrt{1+|Du|^2}\\
		u(x)&=\psi(x)\quad x\in \P\Omega
		\end{split}\right. 
		\ene
	 Here $\Omega$ is a bounded open set with Lipschitz boundary in a Riemannian manifold $N$,  $Du$ denotes the gradient of $u$,  $\alpha$ is a fixed constant, $div$ is the  divergence of $N$ and  $\psi(x)$ is a continuous or $L^1$ function on $\P\Omega$. The  equation in \eqref{main:equation} is called as a translating mean curvature equation.
	 \subsection{Background} The motivation of our paper are two-fold. First it naturally arises in the study of type II mean convex singularities of mean curvature flows or translating mean curvature flow in Euclidean spaces (i.e see White \cite{Whi15}, Wang \cite{Wang11}, Spruck-Xiao \cite{SX17}, Hoffman-Ilmanen-Martin-White \cite{HIMW18} etc). Thus we name \eqref{main:equation} in terms of  ``translating mean curvature equation".  Second we are interested in the Dirichlet problem of mean curvature equations in Riemannian manifolds. The generalization of various boundary problems of mean curvature equations is not trivial ( see Appendix \ref{Example:Section}). \\
		\indent The Dirichlet problem of  \eqref{main:equation} on mean convex domains in $\R^n$ with continuous boundary data is sloved by \cite{Whi15}, Ma \cite{Ma18} and \cite{Wang11} via geometric measure theory, mean curvature type flows and the 
		Monge-Ampere equation theory respectively.  A  key fact in $\R^n$ is that  there is a global solution to the translating mean curvature equation $div(\F{Du}{\omega})=\F{\alpha}{\omega}$ in $\R^n$. This provides a $C^0$ {bound for any $C^2$ solution to \eqref{main:equation}.  \\
    \indent  However in Riemannian manifolds no similar fact exists in general except manifolds with a similar warped metric structure as $\R^n$. As a result no similar $C^0$ bound estimate is known for the solution to \eqref{main:equation} in general Riemannian manifolds. In Appendix \ref{Example:Section} we give an example of  a domain such that no $C^2$ solution exists for \eqref{main:equation} when $\alpha=n$.  We need an approach to attack the Dirichlet problem  \eqref{main:equation} without requiring $C^0$ {\ap} different with the previous methods in \cite{GT01, Whi15,Wang11,Ma18}. \\
    \indent One candidate for our purpose is the generalized solution theory of the Dirichlet problem to the following mean curvature equation
    	\be\label{before:equation}
    div(\F{Du}{\omega})=f(x),\quad  x\in \Omega \quad  u(x)=\psi(x) \quad x\in \P\Omega
    \ene   developed by Miranda and Giusti in \cite{Mir64, Mir77, Giu78, Giu80}. Its key idea is described as follows. One defines an area functional on bounded variation (BV) functions with the property that if $u$ is  smooth and minizes this functional then $u$ satisfies \eqref{before:equation}. One defines a generalized solution $u(x)$ to \eqref{before:equation} taking possible infinity values which locally minimizes the corresponding functional. Without  requiring $C^0$ estimate they showed the existence of the generalized solution to the corresponding Dirichlet problem.  Suppose the domain is mean convex in $\R^n$ and the boundary data is continuous they showed that such generalized solutions are classical. More details are given in \cite{Giu84}. For other applications of the Miranda-Giusti theory see Scholz \cite{Sch04A} and Liang \cite{Lia03}. 
    \subsection{Main Results}Now we follow the idea of Miranda-Giusti generalized solution theory. Before stating the definition, we need some notation. Suppose $N$ is  a Riemannian manifold with metric $g$ and $dim N=n$. Let $\Omega$ be an open set in $N$.
Fix $\alpha > 0$. Suppose $BV(\Omega)$ is the set of bounded variation (BV) functions in $\Omega$ (Definition \ref{Def:BV}). Let $C_0(\Omega)$ and $T_0\Omega$ be the set of all smooth functions and smooth vector fields with compact support in $\Omega$ respectively.\\
\indent  For $u\in BV(\Omega)$, define a conformal area functional 
$$
\AF_{\alpha}(u,\Omega)=\sup_{h \in C_0(\Omega), X\in T_0\Omega}\{\int_{\Omega}e^{\alpha u}(h+\F{1}{\alpha}div(X))d\vf: h^2+\la X,X\ra \leq 1 \}
$$
Observe that 
  when $u\in C^2(\Omega)$ and minimizes $\AF_\alpha(u,\Omega)$ locally, it solves $div(\F{Du}{\omega})=\F{\alpha}{\omega}$. Moreover for $u\in C^1(\Omega)$,  $\AF_\alpha(u,\Omega)$ is the area of the graph of $u(x)$ in $Q_\alpha$ which is the manifold
  $
 N\PLH\R$ with the metric $ e^{2\F{\alpha}{n}}(g+dr^2)$. \\
 \indent  Let $U$ be the subgraph of $u(x)$ defined by $\{(x,t):x\in \Omega, t<u(x)\}$ in $Q_\alpha$. Inspired by Miranda-Giusti \cite{Giu84} we proposed the following definition.
\begin{Def} [Definition \ref{Def:gs}]	Let $\Omega\subset \subset\mathcal{B}$ be two open bounded set in $N$. Let $u(x), \psi(x)$ be two measurable functions taking values in $[-\infty, \infty]$ such that $\Psi$, the subgraph of $\psi(x)$, is a {\Ca} set (Definition \ref{Def:Cs}) in $Q_\alpha$. 
	We say that $u(x)$ is a generalized solution to the {\DP} \eqref{main:equation} on $\Omega$ with boundary data $\psi(x)$ on $\P\Omega$ if
	\begin{itemize} 
		\item [(1)] the subgraph of $u(x)$, $U$,  coincides with  $\Psi$ outside $\bar{\Omega}\PLH\R$ ;
		\item [(2)]for any {\Ca} set $F$ satisfying $F\Delta U\subset K$ where $K$ is a compact set in $\bar{\Omega}\PLH \R$, it holds that 
		\be \label{eq:minimize}
		\int_{K}d|D\lambda_U|_{Q_\alpha} \leq  \int_{K}d|D\lambda_F|_{Q_\alpha}
		\ene 
	\end{itemize}
	where $|D\cf_U|_{Q_\alpha}$ and $|D\cf_F|_{Q_\alpha}$ are Radon measures of generated by $\lambda_U$ and $\lambda_F$ in $Q_\alpha$ respectively (See Theorem \ref{thm:st}). 
\end{Def}
A key point of this definition is Condition (2) above. It says  that $u$ locally minizes  $\AF_\alpha(.,\Omega)$(see Theorem \ref{perimeter}). In our setting  a generalized solution always exists. 
\bt [Theorem \ref{mt:A}]\label{mt:AA}Let $\Omega\subset \subset\mathcal{B}$ be two open bounded set in a  Riemannian manifold $N$. Suppose  $\P\Omega$ is Lipschitz and  $\psi(x)$ is a  measurable function such that its subgraph $\Psi$ is a {\Ca} set in $Q_\alpha$. Then there is a generalized solution to the Dirichlet problem \eqref{main:equation} with boundary data $\psi(x)$. 
\et 
\br The assumption on $\psi(x)$ is not restrictive. For any $\psi(x)\in L^1(\P\Omega)$ we can extend $\psi(x)$ as the trace of a BV function on $\mathcal{B}$ such that its subgraph is a {\Ca} set in $Q_\alpha$. For its proof see Remark \ref{re:mk}.  Remark \ref{Rk:not:gs} gives an example of a nontrival generalized solution to the Dirichlet  problem \eqref{main:equation} with finite boundary data. 
\er 

Whether generalized solutions are classical depends heavily on the geometry of the closure of  $\Omega$, the condition (2) below. 
\bt[Theorem \ref{mt:B}]\label{mt:BB}  Let $N$ be a Riemannian manifold with dimension $n$. Suppose $\Omega$ is a bounded open domain with $\mC^2$ boundary satisfying 
\begin{enumerate} 
	\item $ H_{\P\Omega}=div(\vec{v})\geq 0$ on $\P\Omega$ where $\vec{v}$ is the outward normal vector of $\P\Omega$ ; 
	\item \begin{enumerate}
		\item if $n\leq 7$, no closed embedded minimal hypersurface exists in $\bar{\Omega}$;
		\item if $n>7$, no closed embedded minimal hypersurface with a closed singular set $S$ with $H^{k}(S)=0$ for $k>n-7$ exists in $\bar{\Omega}$ where $H^k$ denotes the $k$-dimensional hausdorff measure on $N$;
		\end{enumerate}
\end{enumerate}
 Then the Dirichlet problem \eqref{main:equation} admits an unique solution $u\in \mC^2(\Omega)\cap \mC(\bar{\Omega})$ for any continous function $\psi(x)$ on $\P\Omega$. Here $\mC^k$ denotes the property of $k$-th differential. 
\et 
\br
In the Miranda-Giusti generalized solution theory \cite{Giu84}, the condition (1) is sufficient to transfer a  generalized solution into a classical solution. However in our setting the condition (2) can not be removed in general. In Appendix \ref{Example:Section} we show that in the upper half sphere $S^+_n$ in the sphere $S_n$ no classical solution to the Dirichlet problem of \eqref{main:equation} exists
  for continuous boundary data for $\alpha=n$. Its boundary $\P S_n^+$ is minimal. 
  \er 
  \br In Euclidean spaces and Hyperbolic spaces all bounded domains satisfy the condition (2) by the maximum principle of singular stationary hypersurface in Ilmanen \cite{Ilm96}. Thus Theorem \ref{mt:BB} generalizes the corresponding results in Euclidean spaces by \cite{Wang11, Whi15, Ma18}. 
  \er
  Now we present the ideas to show Theorem \ref{mt:AA} (Theorem \ref{mt:A}) and Theorem \ref{mt:BB} (Theorem \ref{mt:B}). \\
  \indent One observes that  $\AF_\alpha(u,\Omega)$ is equal to the perimeter of the subgraph of $u$ in $Q_\alpha$ (Theorem \ref{perimeter}). Thus the local minimizer of $\AF_\alpha(u,\Omega)$ corresponds to a local minimizer of the perimeter in $Q_\alpha$ (Theorem \ref{key:lemma}). We called those results as the Miranda's observation.  Then the strategy in Theorem 16.11 of \cite{Giu84} shall work in our setting and yields Theorem \ref{mt:A}. \\
  \indent With the assumption of Theorem \ref{mt:B}  it admits a generalization $u(x)$ with continous boundary data $\psi(x)$ by Theorem \ref{mt:A}. Consider the set $P_+$ given as $\{x\in \Omega: u(x)=\infty\}$.  The finiteness of $\psi(x)$ and the smooth boundary of $\Omega$ implies that $ \P P_+\PLH\R$ is an almost minimal boundary in $Q_\alpha$. This is very different with the case of \cite{Giu84} (see Remark \ref{Rk:P:Ex}). By the mean convex condition and the property of generalized solutions, the regularity of  the almost minimal boundary implies that $\P P_+$ satisfies the opposite of Condition (2) in Theorem \ref{mt:B} (see Theorem \ref{thm:infty}).  This means $P_+$ ($P_-$)  is empty and $u(x)$ is locally bounded. By Theorem \ref{thm:ctwo} $u(x)\in \mC^2(\Omega)$. The boundary continuity of $u(x)$ is established in Lemma \ref{last:step}. This finally yields Theorem \ref{mt:B}.\\
   \subsection{Outline} 
\indent Our paper is organized as follows. In Section 2 we discuss the properties  of BV functions.  In Section 3 we discuss various properties of product area functional $\AF(u,\Omega$) and conformal area fuctional $\AF_\alpha(u,\Omega)$ (see Definition \ref{Def:smooth}).  It is techniquely involved because the method of Theorem 1.17 in \cite{Giu84} can be directly applied. \\
   \indent In Section 4 we record results of the trace of BV functions on Lipschitz domains in Riemannian manifolds.  In Section 5 we show the Miranda's observation upon the perimeter in $Q_\alpha$ and $\AF_\alpha(u,\Omega)$. In Section 6 we show Theorem \ref{mt:A}.  In Section 7 we study various properties of  generalized solutions only taking infinity values from a given generalized solution (Lemma \ref{lm:minimize}). All results are summarized in Theorem \ref{thm:infty}.  In Section 8 we prove Theorem \ref{mt:BB} ( Theorem \ref{mt:B}).  \\
   \indent In Appendix A we record a decomposition result for Radon measures in Riemannian manifolds.  It is useful to prove the $\mC^\infty$ approximation of $\AF_\alpha(u,\Omega)$ in Section 3. In Appendix B we show that the Dirichlet problem \eqref{main:equation} is solvable in sufficiently small normal balls in Riemannian manifolds. It is a key result to show the regularity result in Theorem \ref{thm:ctwo}. In Appendix C we give an example of a domain such that no classical solution to the Dirichlet problem \eqref{main:equation} exists.  In Appendix D we collect some facts on mean curvatures of hypersurfaces.\\
   \indent At last we should remind the reader that our paper depends deeply on the book of Giusti \cite{Giu84}. 
   \subsection{Acknowledgement} This project is supported by the National Natural Science Foundation of China, Grant No.
   11801046.  The author is very grateful to inspired conversations with Xiantao Huang in Sun Yat-sen University, China. 
   \section{BV functions in Riemannian manifolds}
   In this section we discuss the definition of  BV function on Riemannian manifold. Our viewpoint is that a BV function  corresponds  a Radon measure. We also define the convolution of functions  and vector fields in a normal ball for later use. We mainly follow from Chapter 1 of Simon \cite{Simon83}, the books of \cite{Giu84} and \cite{EG15}. 
      \subsection{BV functions}\label{section:BV}
   Let $N$ be a Riemannian manifold with a metric $g$. Let $\la ,\ra$ be the corresponding inner product.  Write $div$ and $d\vf$ for the divergence and volume of $N$ respectively. Suppose $\Omega\subset N$ is an open set. Let $T_0\Omega$ be the set of smooth vector fields  with compact support in $\Omega$ and $d\vf$ be the volume form of $\Omega$. 
      \begin{Def}\label{Def:BV} Let $u\in L^1(\Omega)$, we define 
   	\be
   	||Du||_{N}(\Omega)=\sup\{\int_{\Omega}udiv(X)d\vf, X\in T_0\Omega \quad \la X,X\ra\leq 1\}
   	\ene
   	If $||Du||_{N}(\Omega)<\infty$, we say that $u$ has bounded variation in $\Omega\subset N$ and $u\in BV(\Omega)$ or $u\in BV_N(\Omega)$ to emphasize the manifold $N$. \\
   	\indent If $u\in BV_N(\Omega')$ for any bounded open set $\Omega'\subset \Omega$, we say $u\in BV _{loc, N}(\Omega)$.
    \end{Def}
   \br If $u\in C^1(\Omega)$, the divergence theorem implies that 
\be 
\int_{\Omega} udiv(X)d\vf=-\int_{\Omega}\la X, Du\ra d\vf
\ene 
for any $X\in T_0\Omega$ where $Du$ is the gradient of $u$ on $\Omega$. 
\er 
  A straightforward verification shows that BV functions have a lower semicontinuity as follows. 
  \bt\label{BV:semicontinuity} Suppose $\{u_j\}_{j=1}^\infty\in BV(\Omega)$ and converges to $u$ in $L^1(\Omega)$ as $j\rightarrow +\infty$.  Then 
\be 
||Du||_N(\Omega)\leq \lim_{j\rightarrow +\infty}||Du_j||_N(\Omega)
\ene 
\et
\subsection{Radon measures} Now we shall see that BV functions naturally induce Radon measures. For more details we refer to Chapter 1 in \cite{Simon83}. 
\begin{Def} Let $X$ be a locally compact Hausdorf measure.  A Radon measure on $X$ is an outer measure $\nu$ on $X$ having three properties:
\begin{enumerate}
	\item $\nu$ is Borel regular and $\nu(K)<\infty$ for any compact set $K\subset X$;
	\item $\nu(A)=\inf\{ \nu(U): A\subset U,\quad\text{$U$ open} \}$ for each subset $A\subset X$;
	\item $\nu(U)=\sup\{\nu(K):\text{$K$ compact}\subset U\}$ for each open $U$ in $X$.
\end{enumerate}
\end{Def}
For any  set $X$ let $K_+(X)$ denote the set of non-negative continuous functions $f:X\rightarrow [0,\infty)$ with compact support. . 
  \bt [Remark 4.3 in \cite{Simon83}]
  \label{Rieze:rep}Suppose $X$ is a locally compact Hausdorff space and $\lambda: K_+(X)\rightarrow [0,\infty)$ satisfies  $\lambda(cf)=c\lambda(f)$, $\lambda(f+h)=\lambda(f)+\lambda(h)$ for any constant $c\geq 0$ and $f,g\in K_+(\Omega)$. Then there is a Radon measure $\nu$ on $X$ given by 
\be
\begin{split}
	\nu(U)=\sup\{\lambda(f),f\in K_+(X), supp(f)\subset U,f\leq 1\}\, \forall \,\, \text{open}\,\, U\in X
\end{split}
\ene
such that 
\be \label{eq:increa}
\lambda (f)=\int_{X}fd\nu \quad \forall f\in K_+(\Omega)
\ene 
\et 
Suppose $u(x)\in BV(\Omega)$. Set a nonnegative functional $\lambda_{u}: K^+(\Omega)\rightarrow[0,+\infty)$ as 
\be 
\lambda_{u}(h)=\sup\{\int udiv(X)d\vf, X\in T_0\Omega, \la X,X\ra\leq h^2\}
\ene 
for every $h\in K^+(\Omega)$. It is clear that 
\be 
\lambda_{u}(ch)=c\lambda_{u}(h),\quad  \lambda_{u}(h+h_1)=\lambda_{u}(h)+\lambda_{u}(h_1)
\ene 
where $c$ is any positive constant, $h$ and $h_1\in K^+(\Omega)$. 
Thus Theorem \ref{Rieze:rep} gives a Radon measure from a BV function as follows. 
\bt\label{thm:st}  Let $\Omega$ be an open set in a Riemannian manifold $N$. Suppose $u\in BV_{loc, N}(\Omega)$. \begin{enumerate} 
	\item There is a Radon measure $|Du|_N$ on $\Omega$ such that 
\be 
\int_{\Omega'}fd|Du|_N=\sup\{\int_{\Omega'} udiv(X)d\vf, X\in T_0\Omega', \la X,X\ra\leq f^2\}
\ene 
for any bounded open set $\Omega'\subset \Omega$ and any nonnegative function $f\in L^1(|Du|_N, \Omega)$. 
\item Moreover there is a vector field $\nu$ on $\Omega$ satisfies 
\be 
\int_{\Omega}u div(X)d\vf=-\int_{\Omega}\la X,\nu\ra d|Du|_N 
\ene 
where $\la \nu,\nu\ra=1$ a.e. $|Du|_N$ for any $X\in T_0\Omega$. 
\end{enumerate} 
\et 

 \bp The existence and definition of $|Du|_N$ are from Theorem \ref{Rieze:rep}. Then $|Du|_{N}$ is a Radon measure. Similar as the proof of Theorem 5.10 in Chapter 1 of \cite{Simon83}, there is a monotone nonegative increasing sequence $\{f_j\}_{j=1}^\infty$ such that each $f_j\in K^+(\Omega)$, $f_j\leq f$ and $f_j$ converges to $f$ in $L^1(|Du|_N,\Omega)$. Let $\Omega'$ be any open set in $\Omega$. By \eqref{eq:increa}, 
 $$
 \int_{\Omega'}f_jd|Du|_N =\sup\{\int_{\Omega'} udiv(X)dvol, X\in T_0\Omega', \la X,X\ra\leq f_j^2\}
  $$
  Letting $j\rightarrow +\infty$ on both sides we obtain the conclusion.  The conclusion (2) is from a version of the Riesz representatioin theorem (see Theorem 4.1 in \cite{Simon83}). 
 \ep  
   A relationship between BV functions on the manifold and its conformal manifold is given as follows. 
   \begin{Def} \label{Def:con}Let $M$ be a Riemannian manifold with a metric $g$. Let $\phi(x)>0$ be a smooth positive function on $M$. A conformal manifold $M_{\phi}$ is the smooth manifold $M$ with the metric $\phi^2(x)g$. 
   	\end{Def}
   \bt \label{thm:BV} Take the notation in Definition \ref{Def:con}.  Suppose $\Omega$ is an open set in $M$ and $u\in BV_{loc,M}(\Omega)$. Then 
   \be 
   ||Du||_{M_{\phi}}(\Omega)=\int_{\Omega}\phi^{m-1}(x) d|Du|_{M}
   \ene 
   where $m$ is the dimension of $M$, $|Du|_M$ is the Radon measure given in Theorem \ref{thm:st}.
   \et 
   \br\label{Remark:locbv} Notice that the metric $g$ of $M$ can be written as $\phi^{-2}\phi^2g$. A consequence of Theorem \ref{thm:BV} is that $u\in BV_{loc,M}(\Omega)$ if and only if $u\in BV_{loc, M_\phi}(\Omega)$. 
   \er
   \bp Let $div_\phi$ and $dvol_\phi$ be the divergence and the volume of $M_{\phi}$ respectively. Then $dvol_{\phi}=\phi^{m}dvol$ where $dvol$ is the volume form of $M$. By the definition of the divergence (Page 423 in \cite{Lee2013}), we have 
   \begin{align}
   div_{\phi}(X)dvol_{\phi}&=d(X\llcorner dvol_{\phi})\notag\\ 
   &=(\phi^mdiv(X)+m\phi^{m-1}\la X,\nb \phi\ra )dvol\notag\\
   &=div(\phi^m X)dvol\label{con:vec}
   \end{align}
   where $\nb\phi$ is the gradient of $\phi$ in $M$. 
By Theorem \ref{thm:st} one has that 
   \begin{align}
   ||Du||_{M_\phi}(\Omega)&=\sup\{\int_{\Omega}udiv_{\phi}(X)d\vf_{\phi}: \phi^2\la X,X\ra\leq 1, X\in T_0\Omega\}\notag\\
   &=\sup\{\int_{\Omega} udiv(X')d\vf: \la X',X'\ra\leq \phi^{2m-2}, X\in T_0\Omega\ra \notag\\
   &=\int_{\Omega}\phi^{m-1}(x)d|Du|_{M}\label{eq:representation}
   \end{align}
   The proof is complete.   
   \ep  
    
 \subsection{The convolution of functions and vector fields}  Now we consider how to approximate a function and a smooth vector field in a sufficiently small normal embedded ball in a Riemannian manifold. 
 \begin{Def}\label{defopen} Fix any point $p$ in a Riemannian manifold $M$. Let $exp_p$ be the exponential map near $p$. There is a Euclidean ball $B_r(0)$ centered at $0$ in $\R^{n}$ such that $exp_p: B_r(0)\rightarrow B_r(p)\subset M$ is a diffeomorphism. Via the exponential map we can identify $B_r(p)$ with $B_r(0)$. Moreover the metric of $M$ is represented with 
 $$
 g=g_{ij}dx^i dx^j
 $$
 with the coordinate in $\R^n$. Such ball $B_{r}(p)$ is called as a normal (open) ball. 
 \end{Def} 
 \indent Let $\vp(x)$ be a symmetric smooth mollifier in $\R^n$. i.e. $\vp(x)=\vp(-x)$ and $\vp(x)$ has a compact support in the Euclidean unit ball $B_n(1)$ and 
 \be 
 \int_{\R^n}\vp(x)dx=1
 \ene 
 where $dx$ is the standard Euclidean volume in $\R^n$.  \\
 \indent Suppose $W$ is an open set in $\R^n$. Let $h(x)$ denote a measurable function on $W$ and $X$ denote a tangent vector field on $W$ written as 
 \be \label{def:X}
 X=X^i\F{\P}{\P x_i}
 \ene 
 where $\{\F{\P}{\P x_i}\}^n_{i=1}$ is the standard orthonormal coordinate vector fields in $\R^n$.
 \begin{Def}\label{def:apr} Let $\sigma>0$ be a sufficiently small positive constant. Then $ \vp_{\sigma}*h(x)$, the convolution of $h(x)$, is given as follows.
 \be\label{con:def:fun}
\vp_{\sigma}*h(x)=\int_{\R^n}\F{1}{\sigma^n}\vp(\F{x-y}{\sigma})h(y)dy, \quad x\in W
 \ene 
 where we extend $h(x)$ outside $W$ as $h(x)=0$ for $x\notin W$. For $X$ in \eqref{def:X}, $\vp_{\sigma}*X(x)$ is defined as 
 \be\label{con:def:vec}
 \vp_{\sigma}*X(x)=\vp_{\sigma}*X^i\F{\P}{\P x_i},\quad x\in W
 \ene 
 \end{Def}
  A useful property about the  convolution is 
 \be 
  \int_{\R^n}u(x)\vp_{\sigma}*h(x)dx=\int_{\R^n}h(x)\vp_{\sigma}*u(x)dx
 \ene 
  \bt \label{estimate:thm}Let $B$ be a normal open ball in a {\RM}  with a metric $g=g_{ij}dx^i dx^j$. Let $f$ be a nonnegative continuous function on $B$. Let $h\in \mC(B)$ and let $X$ be a smooth vector field satisfying 
 \be 
 h^2(x)+\la X, X\ra(x)\leq f^2(x)\quad \forall\quad x\in B
 \ene 
 where $\la,\ra$ is the inner product determined by $g$. 
 Then for any $\Sc>0$ and any compact set $K\subset B$ there exists a $\sigma_0=\sigma_0(f, K,g,\Sc)$ s.t. for all $\sigma<\sigma_0$, 
 \be 
 h'^2(x)+\la Y,Y\ra(x)\leq (f(x)+\Sc)^2\quad x\in K
 \ene 
 where $det(g):=det(g_{ij})$, $h'$ and $Y$ are given by 
 \begin{gather}
 h'=\F{1}{\sqrt{\det(g)}}\vp_{\sigma}*(\sqrt{\det(g)}h)\label{set:h}\\
  Y=\F{1}{\sqrt{\det(g)}}\vp_{\sigma}*(\sqrt{\det(g)}X)\label{set:Y}
 \end{gather} 
 \et
 \bp  Let $\sigma_1$ be a positive constant less than the Euclidean distance between $\P\Omega$ and $K$. Since $K$ is compact, for all $\sigma<\F{\sigma_1}{2}$ the function $h'$ in \eqref{set:h} and the tangent vector $Y$ in \eqref{set:Y} are well-defined for $x\in K\subset B$. \\
 \indent Let $\Sc'$ be a small constant determined later. For any $x_0\in K$, there is a positive constant $\sigma_2=\sigma_2(f,g,K,\Sc')<\sigma_1$, such that for all $\sigma<\F{\sigma_2}{2}$ and $y,y'\in B_{x_0}(2\sigma)$ 
 \begin{gather}
 \F{1}{1+\Sc'}g_{ij}(y')\leq g_{ij}(y)\leq (1+\Sc')g_{ij}(y')\label{e:1}\\
 \max_{y,y'\in B_{2\sigma}(x_0)}\F{\sqrt{det(g)}(y')}{\sqrt{det(g)}(y)}\leq 1+\Sc'\label{e:2}\\
 f(y)\leq f(y')+\Sc'\quad\text{ for } y,y'\in B_{2\sigma}(x_0)\label{e:3}
 \end{gather}
 Here $B_{2\sigma}(x_0)$ is the Euclidean ball of $x_0$ with radius $2\sigma$ in $B$. 
 By Definition \ref{def:apr} and \eqref{def:X}, we have 
 \be \label{eq:Yform}
  Y^i=\F{1}{\sqrt{\det(g)}}\vp_{\sigma}*(\sqrt{\det(g)}X^i) \quad \text{and}\quad Y=Y^i\F{\P}{\P x_i}
 \ene 
Fix any point $y\in B_{x_0}(\sigma)$. With a rotation we can assume that $g_{ij}(y)=\sigma_{ik}\sigma_{kj}$ where $(\sigma_{ik})$ is a positive definite matrix.
By \eqref{set:Y}, \eqref{e:1},  \eqref{e:2} and \eqref{eq:Yform} for any $\sigma<\F{\sigma_1}{2}$ we have 
 \begin{align}
 g_{ij}(y)Y^iY^j(y)&= \F{1}{\det(g)(y)}(\vp_{\sigma}*(\sqrt{\det(g)}\sigma_{ik}X^i))^2(y)\notag\\
 &\leq (1+\Sc')^2(\vp_{\sigma}*(\sigma_{ik}X^i)^2)\notag\\
 &= (1+\Sc')^2(\vp_{\sigma}*(g_{ij}(y)X^iX^j)\notag\\
 &\leq (1+\Sc')^3\vp_{\sigma}*(g_{ij}X^iX^j)\label{eq:step}
 \end{align} 
 By \eqref{set:h} a similar derivation implies that 
 \be\label{eq:step:two}
 (h')^2(y)\leq (1+\Sc')^3\vp_\sigma*h^2 
 \ene 
 Combining \eqref{eq:step} with   \eqref{eq:step:two} we obtain 
 \begin{align*}
 (h')^2(y)+ g_{ij}Y^iY^j(y)&\leq (1+\Sc')^2\vp_\sigma* (h^2+g_{ij}X^iX^j)\\
 &\leq (1+\Sc')^3\vp_\sigma*f^2\\
 &\leq (1+\Sc')^3(f(y)+\Sc')^2\quad \text{ by \eqref{e:3}}
 \end{align*}
Because $K$ is compact, we can choose $\Sc'$ small enough such that $(1+\Sc')^3(f(y)+\Sc')^2\leq (f(y)+\Sc)^2$ for all $y\in B_\sigma(x_0)$ and $x_0\in K$. For such fixed $\Sc'$, define $\sigma_0=\sigma_2(f,g,K,\Sc')$. Thus for any  $x_0\in K$, $y\in B_\sigma (x_0)$ and $\sigma<\F{1}{2}\sigma_0$ we have
 \be
 (h')^2(y)+ g_{ij}Y^iY^j(y)\leq f(y)+\Sc
 \ene 
 where $\sigma_1=\sigma_1(f,g,K,\Sc)$. We complete the proof. 
 \ep 
  The following technique result will be very useful in the next section. 
    \bl \label{lm:apr} Let $B$ be a normal open ball with a metric  $g=g_{ij}dx^idx^j$. Suppose $u\in BV(B)$ and $q(x)$ is a smooth function with compact support in $B$. Let $X$ be a smooth vector field on $B$ satisfying $\la X, X\ra\leq 1$ . Then for any $\Sc>0$, there is a $\sigma_0=\sigma_0(u,g,q)>0$ such that for all $\sigma\in (0,\sigma_0)$ 
   \be\label{eq:eve}
   \int_{B}\vp_\sigma*(qu)div(X)d\vf\leq 
   \int_{B}udiv(q Y_\sigma)d\vf-\int_{B}u\la X,\nb q\ra d\vf+ \Sc
   \ene 
   where $Y_\sigma$ is $\F{1}{\sqrt{\det(g)}}\vp_{\sigma}*(\sqrt{\det(g)}X)$ and we assume $X=0$ outside $B$.  
    \el 
    \bp Notice that $d\vf=\sqrt{det(g)}dx$ where $d\vf$ and $dx$ are the volume form of $B$ with respect to $g$ and the Euclidean metric respectively. Moreover 
    $$
    div(X)d\vf=div_{\R^n}(\sqrt{det(g)}X)dx
    $$
    where $div$ and $div_{\R^n}$ are the divergence of $B$ and $\R^n$ respectively. 
    We also view $B$ as an open set in Euclidean space $\R^{n}$. Thus $\vp_\sigma*(qu)$ is well-defined if we choose sufficiently small $\sigma$.
 Then 
    \begin{align}
    \int_{B}\vp_\sigma*(q u) div(X)d\vf &=\int_{\R^n} \vp_\sigma*(qu) div_{\R^n}(\sqrt{det(g)}X)dx\notag\\
    &=\int_{\R^n}q(x)u(x)div_{\R^n}(\vp_\sigma*(\sqrt{det(g)}X))(x)dx\notag\\
    &=\int_{B}u(x)div(qY_\sigma)d\vf-\int_B u\la Y_\sigma,\nb q \ra d\vf \label{eq:comparison}
    \end{align}
    where $Y_\sigma=\F{1}{\sqrt{det(g)}}\vp_{\sigma}*(\sqrt{det(g)}X)$. On the other hand, we have 
    \be 
    \begin{split}
    \int_{B} u \la Y_\sigma, \nb q\ra d\vf &=\int_{\R^n} ug_{ij}\vp_\sigma*(\sqrt{det(g)}X^i)\nb^j q dx\\
    &=\int_{\R^n} X^i \vp_\sigma*(ug_{ij}\nb^j q)\sqrt{det(g)}dx
    \end{split}
    \ene 
    Since $\la X, X\ra\leq 1$,there is a $\sigma_0=\sigma_0(u,g,q)>0$ independent of $X$ such that 
    \begin{align*}
    -\int_{\R^n} X^i \vp_\sigma*(ug_{ij}\nb^j q)\sqrt{det(g)}dx &\leq 
    -\int_{\R^n} X^i ug_{ij}\nb^j q\sqrt{det(g)}dx+\Sc\\
    &=-\int_B u\la X,\nb q\ra d\vf+\Sc
    \end{align*}
Combining the above two inequalities together we obtain \eqref{eq:eve}. The proof is complete. 
  \ep
   \section{Conformal area functionals and their $\mC^\infty$ approximation}\label{sec:conformal-area-functionals:I}
   In this section we define product area functionals and conformal area functionals. Then we obtain their $\mC^\infty$ approximation properties in Theorem \ref{thm:app}. Our proof deeply depends on a decomposition result of Radon measures from the {\Bc} Covering Theorem in Riemannian manifolds (see Theorem \ref{thm:good:decomposition}). 
   \subsection{The conformal area functional} 
    Throughout this section we adopt the following notation. Let $N$ be a complete Riemannian manifold with metric $g$. Write $div$ and $d\vf$ for the divergence and volume form of $N$ respectively. Fix $\Omega$ as an open bounded set in $N$. Let $C_0(\Omega)$ and $T_0(\Omega)$ denote the sets of all smooth functions and smooth vector fields with compact supports in $\Omega$ respectively.  
    \begin{Def}\label{Def:smooth} Let $u(x)$ be a measurable function on $\Omega$. The product area functional of $u$, $\AF(u,\Omega)$, is defined by 
    \be 
    	\begin{split}
    	\AF(u,\Omega):&\equiv \sup\{\int_{\Omega}(h+u div(X))d\vf: h\in C_0(\Omega), X\in T_0\Omega,\\
    	& h^2+\la X, X\ra\leq 1\}
    	\end{split}\label{eq:AF1}
    \ene 
    Let $\alpha>0$ be a fixed a constant. The conformal area product functional $\AF_\alpha(u,\Omega)$, is defined by 
    \be
    	\begin{split}
    	\AF_\alpha(u,\Omega):&\equiv \sup\{\int_{\Omega}e^{\alpha  u}(h+\F{1}{\alpha}div(X))d\vf: h\in C_0(\Omega), X\in T_0\Omega,\\
    	& h^2+\la X, X\ra\leq 1\}
    	\end{split}\label{eq:AF2}
    	\ene
    	\end{Def} 
    The terms ``product area functional" and ``conformal area functional" come from the following fact. \bl\label{lm:smooth} Suppose $u\in \mC^1(\Omega)$. Let $Du$ be the gradient of $u$ in $\Omega$ and $|Du|^2=\la Du,Du\ra$. 
    \begin{enumerate}
    	\item Then 
   $\AF(u,\Omega)$ is the area of the graph of $u(x)$ in the product manifold $N\PLH\R$ with the metric $g+dr^2$ which is $\int_{\Omega}\sqrt{1+|Du|^2}d\vf$.  
   \item And $
\AF_\alpha(u,\Omega)
$
is the area of the graph of $u(x)$ in the product manifold $N\PLH\R$ with the metric $e^{2\alpha\F{ r}{n}}(g+dr^2)$ which is $\int_{\Omega}e^{\alpha u}\sqrt{1+|Du|^2}d\vf$. 
\end{enumerate}
\el
\br The proof is straightforward so we skip it here. If $u\in C^2(\Omega)$ locally minimizes of $\AF_\alpha(u,\Omega)$, then $u$ satisfies the translating mean curvature equation 
$$ 	div(\F{Du}{\omega}) =\F{\alpha}{\omega}$$
on $\Omega$ where $\omega=\sqrt{1+|Du|^2}$. This is the reason we study the conformal area functional $\AF_\alpha (u,\Omega)$ in this paper. 
\er  For any $f\in K^+(\Omega)$, we define two nonnegative functionals:
    \begin{gather}
    \begin{split}
    \lambda_{u,0}(f)&:\equiv \sup\{\int_{\Omega}(h+u div(X))d\vf: h\in C_0(\Omega), X\in T_0\Omega,\\
    & h^2+\la X, X\ra\leq f^2(x)\}
    \end{split}\label{nf:one}\\
    \begin{split}
     \lambda_{u,\alpha}(f)&:\equiv \sup\{\int_{\Omega}e^{\alpha u}(h+\F{1}{\alpha}div(X))d\vf: h\in C_0(\Omega), X\in T_0\Omega,\\
    & h^2+\la X, X\ra\leq f^2(x)\}\label{nf:two}
    \end{split}
    \end{gather}
    It is clear that both $\lambda_{u,0}(.)$ and $\lambda_{u,\alpha}(.)$ are linear on $K^+(\Omega)$. Rewriting the definitions of $\AF(u,\Omega)$ and $\AF_\alpha(u,\Omega)$ we obtain the following two representations:
    \be\label{eq:spt:one}
     \AF(u,\Omega):\equiv\sup\{\lambda_{u,0}(f):f\in K^+(\Omega), f\leq 1\}
 \ene 
     \be\label{eq:spt:two}
\AF_\alpha(u,\Omega):\equiv\sup\{\lambda_{u,\alpha}(f):f\in K^+(\Omega), f\leq 1\}
    \ene 
    In fact the above two formulas are true for any open set in $\Omega$. 
    By Theorem \ref{Rieze:rep}, as in the case of BV functions, we obtain two Radon measures from $\AF(u,\Omega)$ and $\AF_\alpha(u,\Omega)$ as follows. 
    \bt\label{Randon:measure} Let $\Omega$ be an open set in $N$. 
    \begin{enumerate} 
    	\item Suppose $u\in BV_{loc,N}(\Omega)$ such that $\AF(u,\Omega')$ is finite for any bounded open set $\Omega'\subset \Omega$. Then there is a unique Radon measure $\mu_0$ on $\Omega$ satisfying 
    	 $$
    	 \mu_0(\Omega')=\AF(u,\Omega')\quad \forall \Omega'\subset\Omega
    	 $$
    	\item Suppose $u\in BV(\Omega)$ such that $\AF_\alpha(u,\Omega')$ is finite for for any bounded open set $\Omega'\subset \Omega$. Then there is a unique Radon measure $\mu_\alpha$ on $\Omega$ satisfying 
    	$$
    	\mu_\alpha(\Omega')=\AF_\alpha(u,\Omega')\quad \forall \Omega'\subset\Omega
    	$$
    	\end{enumerate}
    \et 
    \bp By Theorem \ref{Rieze:rep} on $\lambda_{u,0}(.)$ in \eqref{nf:one}  and $\lambda_{u,\alpha}(.)$ in \eqref{nf:two} we obtain the existence of $\mu_0$ and $\mu_\alpha$. The two conclusions follow from the definitions of $\AF(u,.)$ in \eqref{eq:spt:one} and $\AF_\alpha(u,.)$ in \eqref{eq:spt:two} respectively. 
    \ep 
  The semicontinous properties are also valid for $\AF(u,\Omega)$ and $\AF_\alpha(u,\Omega)$. 
    \bt\label{thm:semi} Let $\Omega$ be a bounded open set in $N$.
    \begin{enumerate}
    	\item  Suppose $u_k$ converges to $u$ in $L^1(\Omega)$. Then 
          $$
        \AF (u,\Omega)\leq \lim_{k\rightarrow \infty}\inf\AF(u_k,\Omega)
          $$
         \item Suppose $e^{\alpha u_k}$ converges to $e^{\alpha u}$ in $L^1(\Omega)$. Then 
             $$
          \AF_\alpha( u,\Omega)\leq \lim_{k\rightarrow \infty}\inf\AF_\alpha( u_k,\Omega)
          $$
          	\item If $u\in BV_N(\Omega)$, then 
          $$ 
          \max\{||Du||_{N}(\Omega),vol(\Omega)\}\leq  \AF(u,\Omega)\leq vol(\Omega)+||Du||_{N}(\Omega)
          $$
          where $vol(\Omega)$ denotes the volume of $\Omega$ in $N$.
\end{enumerate}    
    \et 
    \bp The conclusion (1) and (2) follow directly from the definitions of $\AF(u,\Omega)$ and $\AF_\alpha (u,\Omega)$ respectively. As for the left inequality in (3), we just let $h\equiv 0$ or $X\equiv 0$ and take the supremum in \eqref{eq:AF1}. The right inequality in (3) just follows from the definition of BV functions. 
    \ep 
  
  \subsection{The $C^\infty$ approximation} In this subsection we show the $C^{\infty}$ approximation properties of $||Du||_{N}(\Omega)$, $\AF(u,\Omega)$ and $\AF_\alpha(u,\Omega)$ together. Notice that the method in (Theorem 1.17, \cite{Giu84}) can not be applied directly into the domains in general Riemannian manifolds because they may be not simply connected.  Neither does the method of the weighted BV function in Baldi \cite{Bal01} because its definition of weighted BV functions does not consider the Riemannian metric. 
  \bt\label{thm:app} Let $\Omega\subset N$ be a bounded open set and $u(x)\in BV_N(\Omega)$. 
  \begin{enumerate} 
  	  \item Then there is a sequence $\{u_k\}_{k=1}^\infty$ in $C^\infty(\Omega)$ s.t $u_k$ converges to $u$ in the $L^1(\Omega)$ sense and 
  	$$
  	\lim_{k\rightarrow \infty}||Du_k||_N(\Omega)=||Du||_N(\Omega)
  	$$
  	\item There is a sequence $\{u_k\}_{k=1}^\infty$ in $C^\infty(\Omega)$ s.t $u_k$ converges to $u$ in the $L^1(\Omega)$ sense and 
$$
  \lim_{k\rightarrow \infty}\AF(u_k,\Omega)=\AF(u,\Omega)
$$ 
 \item In addition suppose $\alpha>0$ and $\AF_\alpha(u,\Omega)$ is finite. Then there is a sequence $\{u_k\}_{k=1}^\infty$ in $C^\infty(\Omega)$ s.t $ e^{\alpha u_k}$ converges to $e^{\alpha u}$ in the $L^1(\Omega)$ sense and 
 $$
 \lim_{k\rightarrow \infty}\AF_\alpha (u_k,\Omega)=\AF_\alpha(u,\Omega)
$$
  \end{enumerate}
  \et 
  \bp \textit{\bf{The case of $\AF(u,\Omega)$}:} \\
 \indent Since $u\in BV_N(\Omega)$ and $\Omega$ is bounded, then $\AF(u,\Omega)$ is finite by (3) in Theorem \ref{thm:semi}. From Theorem \ref{Randon:measure} there is a Radon measure $\mu_0$ satisfying $\mu_0(\Omega')=\AF(u,\Omega')$ for any open set $\Omega'\subset \Omega$.\\
  \indent Fix $\Sc>0$. By Theorem \ref{thm:good:decomposition}, there is a collection of normal open balls $\{B_i\}_{i=1}^\infty$ such that (1) $\Omega\subset \cup_{i=1}^\infty B_i$; (2) there is an integer $\kappa(\Sc)>0$ such that $\{B_1,\cdots,B_n\}_{i=1}^{\kappa(\Sc)}$ is a pairwise disjoint collection with the estimate 
  \be\label{condition:open:ball}
  \mu_0(\Omega)-\Sc\leq \sum_{i=1}^{\kappa(\Sc)}\mu_0(B_i)\leq \mu_0(\Omega);\quad 
  \sum_{i=\kappa(\Sc)+1}^\infty \mu_0(B_i)\leq \Sc
  \ene Let $\{q_i(x)\}_{i=1}^\infty$ be a unit partition subordinate to the cover $\{B_i\}_{i=1}^\infty$, that is, $q_i\in C_0^\infty(B_i)$, $0\leq q_i\leq 1$ and $\sum_{i=1}^\infty q_i=1$. Assume $\Bh\in C_0(\Omega)$ and $X\in T_0\Omega$ satisfying 
  \be\label{eq:assumptiona}
  \Bh^2+\la X, X\ra\leq 1
  \ene 
  \indent Since each $B_i$ is an open normal ball in $N$, we can view $B_i$ as an open set in $\R^{n+1}$ with the metric $g=g_{kl}dx^kdx^l$. In what follows we use the same notation for the functions and vector fields on $B_i$ as those in the corresponding open set in Euclidean spaces.  Suppose $X$ is written as $X^k\F{\P}{\P x_k}$. Then on $B_i$ we have  
  \be\label{eq:assumptionb}
  \Bh^2+g_{kl}X^kX^l\leq 1
  \ene 
   For each $i>0$ we can choose $\sigma_i>0$ such that 
  \begin{gather}
    \int_{B_i}|\vp_{\sigma_i}*(uq_i)-uq_i|d\vf\leq \F{\Sc}{2^i}\label{condition:b}\\
    (q_i \Bh)^2+(q_i)^2 \la Y_{\sigma_i}, Y_{\sigma_i}\ra \leq 1+\Sc\label{condition:c} \quad \text{ by Theorem \ref{estimate:thm}}\\
    \begin{split}
    \int_{B_i}\vp_{\sigma_i}*(uq_i)&div(X))d\vf
    \leq \int_{B_i}(u div(q_i Y_{\sigma_i})d\vf\\
    &-\int_{B_i }u\la X, \nb q_i\ra)d\vf+\F{\Sc}{2^i}
    \end{split}\text{ by Lemma \ref{lm:apr}}\label{eq:st:term}
  \end{gather}
  where $Y_{\sigma _i}^k=\F{\vp_{\sigma_i}*(\sqrt{det(g)}X^k)}{\sqrt{det(g)}}$and $Y_{\sigma_i}=Y_{\sigma_i}^k\P_k$. \\
  \indent 
 Now we define $u_{\Sc}$ as 
  \be 
  u_{\Sc}=\sum_{i=1}^\infty \vp_{\sigma_i}*(uq_i)
  \ene 
  By \eqref{condition:b} one sees that 
  \be \label{integral:estimate}
  \int_\Omega |u_\Sc-u|d\vf\leq \sum_{i=1}^\infty \int_{B_i}|\vp_{\sigma_i}*(u q_i)-uq_i|d\vf\leq \Sc
  \ene 
  Consider
  \be \label{eq:key:term}
  \int_{\Omega}(\Bh+u_\Sc div(X))d\vf\leq \sum_{i=1}^\infty\int_{B_i} (\Bh q_i+\vp_{\sigma_i}*(uq_i)div(X))d\vf
  \ene 
From \eqref{eq:st:term} we have 
\begin{align*}
\int_{B_i} (\Bh q_i&+\vp_{\sigma_i}*(uq_i)div(X))d\vf \\
&\leq \int_{B_i}(\Bh q_i+u div(q_i Y_{\sigma_i}))d\vf
-\int_{B_i}u\la X,\nb q_i\ra d\vf+\F{\Sc}{2^i}\\
&\leq (1+\Sc)\mu_0(B_i)-\int_{B_i}u\la X,\nb q_i\ra d\vf+\F{\Sc}{2^i}
\end{align*}
Combining the condition \eqref{condition:open:ball} with \eqref{eq:key:term} we obtain 
\be\label{sider}
\begin{split}
\int_{\Omega}(\Bh+u_\Sc div(X))d\vf &\leq (1+\Sc)\sum_{i=1}^\infty \mu_0(B_i)+\Sc\\
&\leq (1+\Sc)(\mu_0(\Omega)+\Sc)+\Sc
\end{split}
\ene
In the first inequality we apply the fact that $\sum_{i=1}^\infty \int_{B_i} u\la X,\nb q_i\ra d\vf=0$ because of $\sum_{i=1}^\infty q_i(x)=1$. Taking supremum for all $h,X$ satisfying \eqref{eq:assumptiona} yields that 
\be \label{destimate}
\AF(u_\Sc,\Omega)\leq (1+\Sc)(\AF(u,\Omega)+\Sc)+\Sc
\ene 
Now take a sequence $\Sc_k\rightarrow 0$ as $k\rightarrow \infty$. By \eqref{integral:estimate} and \eqref{destimate} we obtain a smooth sequence $\{u_{\Sc_k}\}_{k=1}^\infty$ such that $u_{\Sc_k}$ converges to $u$ in the $L^1$ sense and 
\be\label{eq:threethree}
\lim_{k\rightarrow \infty}\sup\AF(u_{\Sc_k},\Omega)\leq \AF(u,\Omega)
\ene 
From Theorem \ref{thm:semi} we have $\lim_{k\rightarrow \infty}\inf\AF(u_{\Sc_k},\Omega)\geq \AF(u,\Omega)$. Thus $\{u_{\Sc_k}\}_{k=1}^\infty$ is the desirable sequence. We obtain the conclusion (1).\\
\indent \textit{\bf{The case of $||Du||_N(\Omega)$:}} \\
\indent Th proof of the conclusion (1) is similar to that of the conclusion (2). We just let $\Bh\equiv 0$ in the derivation from \eqref{eq:assumptiona} to \eqref{eq:threethree}. Recall that
\be
||Du||_N(\Omega)=\sup\{\int_{\Omega} u div(X) d\vf, X\in T_0(\Omega), \la X,X\ra\leq 1\}
\ene
 Now we define $u_{\Sc}$ as 
\be 
u_{\Sc}=\sum_{i=1}^\infty \vp_{\sigma_i}*(u q_i)
\ene 
where $\sigma_i$ satisfies \eqref{condition:b}--\eqref{eq:st:term}. 
 With the exact same derivation in the proof of the conclusion (2)  we can find a sequence 
 $\{\Sc_i\}_{i=1}^\infty$ converging to $0$ such that 
 $$
 \lim_{i\rightarrow+\infty}||Du_{\Sc_i}||_N(\Omega)=||Du||_N(\Omega)
 $$ 
 This is the conclusion (1). \\
  \indent \textit{\bf{The case of $\AF_\alpha(u,\Omega)$:}}\\
   The idea to derive the conclusion (3) is also similar to the proof of conclusion (2). But the approximating sequence shall have a little modification.\\
   \indent By Theorem \ref{Randon:measure}, there is a unique Radon measure $\mu_\alpha$ on $\Omega$ such that $\mu_\alpha(\Omega')=\AF_\alpha(u,\Omega')$ for $\forall$ open $\Omega'\subset \Omega$. Moreover $\mu_\alpha(\Omega)$ is finite. 
    Fix $\Sc>0$. By Theorem \ref{thm:good:decomposition}, there is a collection of open sets $\{B_i\}_{i=1}^\infty$ such that (1) each $B_i$ is an open normal ball in $\Omega$, $\Omega\subset \cup_{i=1}^\infty B_i$; (2) there is an integer $N(\Sc)$ such that $\{B_1,\cdots,B_n\}_{i=1}^{N(\Sc)}$ is a pairwise disjoint collection with the estimate 
   \be\label{condition:lab}
   \mt(\Omega)-\Sc\leq \sum_{i=1}^{N(\Sc)}\mt(B_i)\leq \mt(\Omega);\quad 
   \sum_{i=N(\Sc)+1}^\infty \mt(B_i)\leq \Sc
   \ene
   Now assume $h\in C_0(\Omega)$ and $X\in T_0\Omega$ satisfying 
   \be\label{eq:assumptionc}
   h^2+\la X, X\ra\leq 1
   \ene 
   Each $B_i$ can be viewed as an open set in $\R^n$ with the metric $g=g_{kl}dx^kdx^l$.\\
   \indent On each $B_i$ assume $X=X^k\P_k$. Then 
      \be 
        h^2+g_{kl}X^kX^l\leq 1
      \ene 
       Let $\{q_i(x)\}_{i=1}^\infty$ be a unit partition subordinate to the cover $\{B_i\}_{i=1}^\infty$. 
        For each $i>0$ we can choose $\sigma_i>0$ such that 
      \begin{gather}
     \int_{B_i}|\vp_{\sigma_i}*(e^{\alpha u}q_i)-e^{\alpha u}q_i|d\vf\leq \F{\Sc}{2^i}\label{condition:n:b}\\
     (q_i h')^2+(q_i)^2 g_{kl}Y_{\sigma_i}^kY_{\sigma_i}^l\leq 1+\Sc\label{condition:n:c}\\
        \begin{split}
     \int_{B_i}\vp_{\sigma_i}*(e^{\alpha u}q_i)&div(X))d\vf
     \leq \int_{B_i}(e^{\alpha u} div(q_i Y_{\sigma_i})d\vf\\
     &-\int_{B_i }e^{\alpha u}\la X, \nb q_i\ra d\vf+\F{\Sc}{2^i}
     \end{split}\label{condition:n:d}
     \end{gather}
     where $Y_{\sigma _i}^k=\F{\vp_{\sigma_i}*(\sqrt{det(g)}X^k)}{\sqrt{det(g)}}$, $h'=\F{\vp_{\sigma_i}*(\sqrt{det(g)}h)}{\sqrt{det(g)}}$ and  $Y_{\sigma_i}=Y_{\sigma_i}^k\P_k$.
    The proof of the above arguments is similar as that of the case of $\AF(u,\Omega)$ just replacing $u$ with $e^{\alpha u}$. In particular \eqref{condition:n:d} is from Lemma \ref{lm:apr}. \\
     \indent Now we define $u_{\Sc}$ as 
     \be \label{eq:con:aphla}
     e^{\alpha u_{\Sc}}=\sum_{i=1}^\infty \vp_{\sigma_i}*(e^{\alpha u}q_i)
     \ene  
     This definition is well-defined because the right hand is a finite positive summation at any point $x\in \Omega$. According to \eqref{condition:n:b} we have 
     \be 
     \int_\Omega |e^{\alpha u_\Sc}-e^{\alpha u}|d\vf\leq \sum_{i=1}^\infty \int_{B_i}|\vp_{\sigma_i}*(e^{\alpha u}q_i)-e^{\alpha u}q_i|d\vf\leq \Sc
     \ene 
     Replacing $u$ with $e^{\alpha v}$ in \eqref{eq:st:term} and taking the same derivation we obtain
     \be\label{eq:key:one}
     \begin{split}
         \int_{B_i}\vp_{\sigma_i}*(e^{\alpha u}q_i)div(X))d\vf&=\int_{\R^n}(e^{\alpha u} div(q_i Y_{\sigma_i})\\
         &-e^{\alpha u}\la Y_{\sigma_i}, \nb q_i\ra)d\vf
         \end{split}
      \ene
      On the other hand 
      \be\label{eq:key:two}
      \begin{split}
      \int_{B_i}\vp_{\sigma_i}*(e^{\alpha u}q_i) hd\vf&=\int_{B_i}\vp_{\sigma_i}*(e^{\alpha u}q_i) h\sqrt{det(g)}dx\\
      &=\int_{B_i}\vp_i e^{\alpha v} h'd\vf
      \end{split}
      \ene 
      where $h'=\F{\psi_{\sigma_i}*(\sqrt{det(g)}h)}{\sqrt{det(g)}}$. By \eqref{eq:key:one} and \eqref{eq:key:two} we compute 
      \begin{align*} 
      \int_{\Omega}&e^{\alpha u_\Sc} (h +\F{1}{\alpha}div(X))d\vf=\sum_{i=1}^\infty\int_{B_i}\vp_{\sigma_i}*(e^{\alpha u}q_i)(h+\F{1}{\alpha}div(X))d\vf\\
      &\leq \sum_{i=1}^\infty\{\int_{B_i}e^{\alpha u}(q_i h'+\F{1}{\alpha}div(q_iY_{\sigma_i}))d\vf-\F{1}{\alpha}\int_{B_i}e^{\alpha u}\la X, \nb q_i\ra d\vf+\F{1}{\alpha}\F{\Sc}{2^i}\}\\
      &\leq (1+\Sc)\sum_{i=1}^\infty\mt(B_i)+\F{\Sc}{\alpha}
      \end{align*}
    For the first term in the last inequality we combine  \eqref{condition:n:c} and the definition of $\AF_\alpha(u,B_i)$ together. For the second term in the last inequality we apply \eqref{condition:n:d} and the fact $\sum_{i=1}^\infty q_i\equiv 1$ on $\Omega$. Continuing applying the assumption \eqref{condition:lab} into the above estimate, we obtain 
    \be
    \begin{split}
    \int_{\Omega} e^{\alpha u_\Sc}(h+\F{1}{\alpha}div(X))d\vf&\leq (1+\Sc)(\mt(\Omega)+\Sc)+\F{\Sc}{\alpha}\\
    &=(1+\Sc)(\AF_\alpha( u,\Omega)+\Sc)+\F{\Sc}{\alpha}
    \end{split}
    \ene 
    Now we arrive a similar position as in \eqref{sider} when we show the conclusion (1). With a similar derivation we can achieve the conclusion (3). \\
    \indent The proof of Theorem \ref{thm:app} is complete.
  \ep 

  A corollary of the conclusion (2) of Theorem \ref{thm:app} is the compactness of BV function on Lipschitz domains in $N$. Its proof is similar to that in the Euclidean spaces. (See Theorem 2, Section 5.2.2 in \cite{EG15}). 
  \bt Suppose $\Omega$ is an open bounded set in $N$ with Lipschitz boundary. Suppose a sequence $\{u_i\}_{i=1}^\infty\in BV(\Omega)$ satisfies that 
  \be 
  \int_\Omega|u_i|d\vf+||Du_i||_N(\Omega)\leq c \quad\text{for all $i$}
  \ene 
  where $c$ is a fixed constant. Then there is a $u\in BV(\Omega)$ such that there is  a subsequence of $\{u_i\}_{i=1}^\infty$ converging to $u$ in $L^1(\Omega)$. 
  \et 
  \bp For any $\Sc>0$,  (2) in Theorem \ref{thm:app} implies that there is a set given by 
  \be \notag
  \begin{split}
   A_\Sc=\{u_j'\in C^\infty(\Omega):
   &\int_{\Omega}|u_j-u_j'|d\vf\leq \Sc,\\
   & ||Du'_j||_N(\Omega)=\int_{\Omega}|Du_j'|d\vf\leq ||Du_j||_N(\Omega)+\Sc
   \end{split}
  \ene 
  Thus $A_\Sc\in \mW^{1,1}(\Omega)$. Here $\mW$ denotes the Sobolev space. For more details see Appendix B. Arguing as Theorem 7.25 of  Chapter 7 in  \cite{GT01}, the fact that $\Omega$ is Lipschitz implies $\mW^{1,1}(\Omega)$ can be isometrically embedded into $\mW^{1,1}_0(\Omega')$. Here $\Omega'$ is a larger open set containing $\Omega$. By the Sobolev embedding theorem, $\mW^{1,p}_0(\Omega')$ is isometrically embedded into $L^1(\Omega')$. This implies that the set $\cup_{1>\Sc>0}A_\Sc$ is a precompact set in $L^1(\Omega)$. Letting $\Sc\rightarrow 0$ we can subtract the desirable sequence. The proof is complete.  
  \ep  
\section{The trace of BV functions} 
The traces of BV functions in {\RM}s play an important role in our computation on  generalized solutions to the Dirichlet problem \eqref{main:equation} (see Definition \ref{Def:gs}). We will collect related results on this topic. Most of their proofs will be skipped because of no essential modification comparation with those in corresponding Euclidean versions. For exact proofs, we refer to Chapter 2 in \cite{Giu84} and Chapter 5 in \cite{EG15}.  An application is given in Lemma \ref{lm:appr}.  
\subsection{The trace of BV functions}  The definition of the trace is from the following divergence theorem.
\bt[Theorem 5.9 \cite{EG15} or Lemma 2.4 \cite{Giu84}]\label{thm:trace}  Let $\Omega$ be a bounded open set in a Riemannian manifold $N$ with Lipschitz boundary $\P\Omega$. Then there is a linear functional  $\mT: BV_{loc,N}(\Omega) \rightarrow L^1(\P\Omega)$ such that for any $u\in BV_{loc,N}(\Omega)$ it holds that 
\be\label{eq:divergence:result}
\int_{\Omega} udiv(X)d\vf=-\int_{\Omega}\la X, \nu\ra d|Du|_N +\int_{\P\Omega} \mT u\la X,\vec{v}\ra d\vf_{\P\Omega}
\ene
for any smooth vector field $X$ with compact support in $N$. Here $|Du|_N$ is the Radon measure generated by $u$(see Theorem \ref{thm:st}), $\vec{v}$ is the inward normal vector of $\P\Omega$ and $d\vf_{\P\Omega}$ is the volume form of $\P\Omega$, $\nu$ is a vector field on the tangent bundle of $\Omega$ such that $\la \nu,\nu\ra=1$ a.e.-$|Du|_N$. 
\et 
Now we define the trace as follows.
\begin{Def}\label{def:trace} With the the above result, for $u\in  BV_{loc,N}(\Omega)$ the function $\mT u$ up to a measure zero set in $\P\Omega$  is called  the trace of $u$ on $\P\Omega$. $\mT$ is referred as the trace operator. 
\end{Def} 
The following two properties are not hard to check. 
\bl \label{lm:pro:trace}Let $\Omega$ be given in Theorem \ref{thm:trace}. Suppose $u\in BV_N(\Omega)$ and $\mT u$ is its trace on $\P\Omega$. It holds that
\begin{enumerate}
	\item  $\mT(u-j)=\mT u-j$ for a fixed constant $j$;
	\item (Proposition 2.6 in \cite{Giu84}) for a sequence $\{u_k\}_{k=1}^\infty\in BV_N(\Omega)$ such that $u_k$ converges to $u$ in $L^1(\Omega)$ and $||Du_j||_N(\Omega)$ converges to $||Du||_N(\Omega)$, $\mT u_k$ converges to $\mT u$ as $k\rightarrow +\infty$ in $L^1(\P\Omega)$. 
\end{enumerate}
\el 
\bp By Theorem \ref{thm:st} two Radon measures 
$|Du|_N$ and $|D(u-j)|_N$ are equal to each other. Fix $j\in \R$. By \eqref{eq:divergence:result}, we have 
\begin{align*}
\int_{\Omega} (u-j)div(X)d\vf&=-\int_{\Omega}\la X, \nu\ra d|Du|_N +\int_{\P\Omega} (\mT u-j)\la X,\vec{v}\ra d\vf_{\P\Omega}\\
&=-\int_{\Omega}\la X, \nu\ra d|D(u-j)|_N +\int_{\P\Omega} (\mT u-j)\la X,\vec{v}\ra d\vf_{\P\Omega}
\end{align*}
for any smooth vector field $X$. Due to the arbitrariness of $X$, \eqref{eq:divergence:result} implies that the trace of $u-j$ is $\mT u-j$. We obtain the conclusion (1).
As for the proof of the conclusion (2), we refer to Proposition 2.6 in \cite{Giu84} or the proof of step (5) of Theorem 5.9 in \cite{EG15}. 
\ep 
Now we illustrate the trace is indeed the ``boundary value'' of BV functions in the approximation sense.\\
\indent 
Since $\P\Omega$ is Lipschitz and bounded, there is a $r_0>0$ such that for every $r\in (0, r_0)$ if $y\in \Omega$ satisfies  $d(y,\P\Omega)\leq r$ there is only one $x\in \P\Omega$ with the property $d(x,y)=d(x,\P\Omega)$. Thus we define a set 
$$
\Gamma_{r}:=\P\Omega\PLH(0, \Sc)=\{(x,s): x\in \P\Omega, s\in (0,r)\}
$$
for $r<r_0$ with 
the corresponding metric $g=\sigma(x,r)+dr^2$. Here $\sigma(x,0)$ is just the induced metric of $\P\Omega$  and $\sigma(x,r)$ is smooth {\wrt} $x$ and {\Lp} {\wrt} $r$. With this local coordinate on  $\Gamma_{r_0}$ near $\P\Omega$, for every $f\in BV(\Omega)$ we define a  function on  $\P\Omega$ as 
\be \label{def:sc}
f_r(x): =f(p) \quad \text{ where } p =(x,r)\in \Gamma_{r_0}
\ene 
\bl \label{lm:cbv} Let $\Omega$ be an open bounded domain in a Riemannian manifold $N$ with {\Lp} boundary. Suppose $f\in BV _{loc,N}(\Omega)$. Let $f_r$ be the function given in \eqref{def:sc} defined on $\P\Omega$. Then there is a Lebesgue measure zero set $E$ in a small interval $(0,r_0)$
\be\label{local:existence}
\lim_{r\rightarrow 0, r\notin E}f_r(x)=\mT f(x)\quad \text{locally in\quad} L^1(\P\Omega)
\ene 
If $\P\Omega$ is unbounded, then $f_r$ is defined locally on $\P\Omega$ and \eqref{local:existence} is valid locally. That is, $r_0$ and $E$ are locally determined near a fixed point on $\P\Omega$. 
\el 

\bp We can assume that $\Omega$ is bounded. The general case could be easily derived by the unity partition technique.\\
\indent  Our proof just follows from the the step 4 and step 5 in that of Theorem 5.9 \cite{EG15}. Thus we only describe its sketch here. \\
\indent First assume $f\in \mC^\infty(\Omega)\cap BV (\Omega)$. Then $\{f(x,r)\}_{r_0>r>0}$ is a Cauchy sequence as $r\rightarrow 0^+$ because $f\in BV(\Omega)$ and 
\be \label{eq:st:ones}
 \int_{\P\Omega}|f_r-f_s|d\vf_{\P\Omega}\leq C ||Df||_N(\Gamma_{r})
\ene 
for any $r>s\in (0, r_0)$. 
Here $C$ is  a constant only depending on the metric on $\Omega$ and $\vf_{\P\Omega}$ is the induced volume form from $N$. 
Thus $\{f_r\}$ is the desirable Cauchy sequence in $L^1(\P\Omega)$ as $r\rightarrow 0^+$ converging to $\mT f$.  In \eqref{eq:st:ones} fixing $r$ and letting $s$ go to zero, one sees that 
\be \label{eq:condition:newtwo}
\int_{\P\Omega}|f_r-\mT f|d\vf_{\P\Omega}\leq C ||Df||_N(\Gamma_{r})
\ene 
for all $r\in (0,r_0)$. Let $E$ be the empty set. This gives  \eqref{local:existence}. \\
 \indent Next consider $f\in BV(\Omega)$. By the conclusion (1) in Theorem \ref{thm:app}, there is a sequence $\{f_k\}_{k=1}^\infty\in C^\infty \cap BV(\Omega)$ such that $f_k$ converges to $f$ in $L^1(\Omega)$ and $||Df_k||_{N}(\Omega)$ converges to $||Df||_N(\Omega)$ as $k\rightarrow \infty$. \\
 \indent By the Fubini theorem we conclude that except a Lebesgue measure zero set $E$ in $(0,r_0)$ such that for every $r\in (0, r_0)\backslash E$ it holds that 
 \begin{enumerate}
\item  $
 ||Df||_N(\P\Omega_r)=0
 $ 
 where  $\P\Omega_r=\{y\in \Omega: d(y,\P\Omega)=r\}$.
 \item $f_k(x,r)$ converges to $f(x,r)$ in $L^1(\P\Omega)$ because $f_k(p)$ converges to $f(p)$ in $L^1(\Gamma_{r_0})$ and $p=(x,r)\in \Gamma_{r_0}$; 
 \end{enumerate} 
Because of (1), Proposition 1.13 in \cite{Giu84} implies that 
$$\lim_{k\rightarrow +\infty}||Df_k||_N(\Gamma_r)=||Df_k||_N(\Gamma_r)$$
for every $r\in (0, r_0)\backslash E$. Now we apply \eqref{eq:condition:newtwo} for $\{f_k\}\in C^{\infty}(\Omega)\cap BV(\Omega)$. We have $\mT f_k$ converges to $\mT f$ in $L^1(\P\Omega)$. By the conclusion (2) of Lemma \ref{lm:pro:trace}  we take the limit as $k\rightarrow +\infty$ and obtain 
\be 
\int_{\P\Omega}|f_r-\mT f|d\vf_{\P\Omega}\leq C ||Df||_N(\Gamma_{r})
\ene 
Letting $r\rightarrow 0^+$ in $(0,r_0)\backslash E$ we obtain the conclusion. 
\ep 
First the trace only depends the conformal equivalent class of metrics. 
\bl \label{lm:conclass}The trace operator  $\mT$ in Definition \ref{def:trace} only depends on the conformal equivalent class of the metric $g$. 
\el
\bp  Suppose $\phi>0$ be a smooth function on $N$. Let $N$ and $N_\phi$ be two conformal Riemannian manifolds with the metric $g$ and $\phi^2g$ respectively. Let $\Omega$ be an open set in $N$ with {\Lp} boundary. \\
\indent  In the followings, the qualities with (without) the low index $\phi$ denote the ones related to the  manifold $(N, \phi^2 g)$ ($(N,g)$). \\
\indent  Let $\mT u$ be the trace of $u$ on $\P\Omega$ where $u\in BV _{loc,N}(\Omega)$. To show the Lemma, it is suffice to show that 
$\mT u=\mT_\phi u$.\\
\indent  By Remark \ref{Remark:locbv}, we have $u\in BV_{loc, N_\phi}(\Omega)$. Let $X$ be a smooth vector field. By \eqref{eq:divergence:result} 
\be\label{eq:nophi}
\begin{split}
\int_{\Omega} udiv(\phi^{n}X)d\vf&=-\int_{\Omega}\phi^{n}\la X, \nu\ra d|Du|_N \\
&+\int_{\P\Omega} \mT u\la \phi^{n}X,\vec{v}\ra d\vf_{\P\Omega}
\end{split}
\ene 
where $n=dim N$. By \eqref{con:vec}, the divergence theorem and Theorem  \ref{thm:BV} we have 
\begin{gather}
div(\phi^{n} X)dvol=div_\phi(X)d\vf_\phi\\
\la X,\nu_\phi\ra_\phi d|Du|_{N_\phi}=\phi^{n}\la X,\nu\ra d|Du|_N\\
\la X,\vec{v}_\phi\ra_\phi d\vf _{\P\Omega,\phi}=\la \phi^{n}X,\vec{v}\ra d\vf_{\P\Omega}
\end{gather}
\indent Based on the above identies we compare \eqref{eq:nophi} with \eqref{eq:divergence:result} in the case of $\mT_\phi u$. This gives that $\mT_\phi u=\mT u$ in the sense of $L^1(\P\Omega)$. The proof is complete. 
\ep 
The trace is helpful to compute the bounded variation on BV functions on {\Lp} boundaries.  
\bt[Proposition 2.8 in \cite{Giu84}]\label{thm:bd:trace} Suppose $\Omega^+$ and $\Omega^-$ are two open sets sharing a Lipschtiz boundary $\P\Omega$ in a Riemannian manifold $N$. Let $f_1\in BV(\Omega^+)$ and $f_2\in BV(\Omega^-)$. Let $\Omega=\Omega^+\cup \Omega^-$.  Define a function $f:\Omega\rightarrow \R$ by
\be
f=\left\{\begin{split}
	&f_1 &\quad \text{in}\quad \Omega^+\\
	&f_2 &\quad \text{in}\quad \Omega^{-}
\end{split}\right.
\ene 
Then $f\in BV(\Omega)$ and 
\be\label{eq:bd}
\int_{\P\Omega}|\mT f_1-\mT f_2|d\vf_{\P\Omega}=\int_{\P\Omega}d|Df|_N
\ene
Here $\mT$ denotes the trace and $d\vf_{\P\Omega} $ is the volume form of $\P\Omega$ with the induced metric. 
\et 
 \indent The following exension lemma shows that the trace could be any $L^1$ function on Lipschitz functions. 
\bl[Proposition 2.15 \cite{Giu84}]\label{lm:st:extend} Let $\Omega$ be a bounded open set with Lipschitz boundary. Let $\psi(x)\in L^1(\P\Omega)$. Then there is a function $f\in \mW^{1,1}(\Omega)$such that $\mT f  =\psi(x)$ on $\P\Omega$. In particular if $\psi(x)\in \mC(\P\Omega)$ then such $f$ is continuous and bounded.  Here $\mW$ denotes the sobolev space (see Definition \ref{def:sob:space}) and it is clear that $\mW^{1,1}(\Omega)\subset BV(\Omega)$. 
\el 
\bp The proof follows exactly from that of Proposition 2.15 in \cite{Giu84}. In particular from the construction of $f$ it is clear if $\psi(x)\in \mC(\P\Omega)$ then $f$ is bounded. Thus we skip its proof here. 
\ep 
\subsection{An application} Now we see a direct application of the trace of BV function. The following definition will be repeatedly used  in the remainder of this paper. 
\begin{Def}\label{Def:notation} Suppose $N$ is a Riemannian manifold. We write $Q$ for the manifold $N\PLH\R$ equipped with the product  $g+dr^2$ and write $Q_\alpha$ for the manifold $N\PLH\R$ equipped with the product $e^{2\F{\alpha r}{n}}(g+dr^2)$. 
	\end{Def}
\bl\label{lm:appr} Suppose $\Omega$ is a bounded open set in $N$. 
For $T>0$ define   \be\label{def:ut}
u_T(x):=\max\{\min(u(x),T), -T\}
\ene 
Let $U_T$ and $U$ be the subgraph of $u_T(x)$ and $u(x)$ respectively.  Then 
\begin{enumerate}
	\item $
	\lim_{T\rightarrow+\infty}||D\lambda_{U_T}||_Q(\Omega\PLH\R)= ||D\lambda_U||_Q(\Omega\PLH\R)
	$
	\item  $
	\lim_{T\rightarrow+\infty}||D\lambda_{U_T}||_{Q_\alpha}(\Omega\PLH\R)= ||D\lambda_U||_{Q_\alpha}(\Omega\PLH\R)
	$. 
\end{enumerate}
\el
\br\label{rk:bounded:approximation} Suppose $u\in BV_N(\Omega)$. Arguing as in Theorem 2.3 of \cite{Giu84} we conclude $$\lim_{T\rightarrow +\infty}||Du_T||_N(\Omega)=||Du||_N(\Omega)$$
\er
\bp We can assume  $||D\lambda_U||_Q(\Omega\PLH\R)$ and $||D\lambda||_{Q_\alpha}(\Omega\PLH\R)$ are finite. Otherwise the conclusion (1) and (2) are trivial by Theorem \ref{BV:semicontinuity}. This means that  $u\in L^1(\Omega)$ when we show the conclusion (1) and $e^{\alpha u}\in L^2(\Omega)$ when we show the conclusion (2). \\
\indent Now we begin to show the conclusion (1). Fix $T>0$.  Since $\lambda_{U_T}$ converges to $\lambda_U$ locally in $L^1(\Omega\PLH\R)$. From Theorem \ref{BV:semicontinuity} it is sufficient to show that 
\be \label{eq:uvst}
||D\lambda_U||_Q(\Omega\PLH\R)\geq 
\lim_{T\rightarrow+\infty}\sup||D\lambda_{U_T}||_Q(\Omega\PLH\R)
\ene  
 
Moreover beacause $\lambda_{U_T}$ is equal to $\lambda_U$ on the set $\Omega\PLH(-T,T)$, we have 
\begin{align}
||D\lambda_U||_Q(\Omega\PLH\R)
&\geq ||D\lambda_{U_T}||_Q(\Omega\PLH(-T,T))\label{computation:st}
\end{align} 
 Since $\lambda_{U_T}$ is a constant on the open set $\Omega\PLH(-\infty,-T)$ and $\Omega\PLH(T,+\infty)$, 
$$
||D\lambda_{U_T}||_Q(\Omega\PLH(-\infty,-T))=||D\lambda_{U_T}||_Q(\Omega\PLH(T,+\infty))=0
$$ 
We set 
\be\label{notation:ET}
E_{T+}=\{x\in \Omega:u(x)>T\}\quad    E_{T-}=\{x\in \Omega:u(x)<-T\}
\ene 
By Theorem \ref{thm:bd:trace}
\be\label{set:time}
||D\lambda_{U_T}||_Q(\Omega\PLH\{\pm T\})=\int_{E_{ T\pm }}d\vf 
\ene 
where $d\vf$ is the volume of $\Omega$. 
Since $u\in L^1(\Omega)$ and $\vf(\Omega)$ is bounded. Then the limit of \eqref{set:time} is $0$ as $T\rightarrow +\infty$. Continuing the computation in \eqref{computation:st} we obtain 
\begin{align*} 
||D\lambda_U||_Q&(\Omega\PLH\R)
\geq\lim_{T\rightarrow +\infty}\sup ||D\lambda_{U_T}||_Q(\Omega\PLH(-T,T))\\
&=\lim_{T\rightarrow +\infty}\sup ||D\lambda_{U_T}||_Q(\Omega\PLH(-R,R))+\lim_{T\rightarrow+\infty}||D\lambda_U||_Q(\Omega\PLH\{\pm T\})\\
&=\lim_{T\rightarrow +\infty}\sup ||D\lambda_{U_T}||_Q(\Omega\PLH(-R,R))
\end{align*}
This gives \eqref{eq:uvst}. Thus we arrive the conclusion (1). \\
\indent The proof of the conclusion (2) is similar to that of the conclusion (1). The differences come from the metric in $Q_\alpha$. By Theorem \ref{thm:bd:trace} 
\be\label{st:eq:rt}
||D\lambda_{U_T}||_{Q_\alpha}(\Omega\PLH\{\pm T\})=\int_{E_{ T\pm}}e^{\pm\alpha T}d\vf 
\ene 
Notice that 
\be 
\int_{E_{T+}}e^{\alpha T}d\vf \leq  \int_{E_{T+}}e^{\alpha u(x)}d\vf \quad \int_{E_{T-}}e^{-\alpha T}d\vf\leq e^{-\alpha T}vol(\Omega)
\ene 
Since $e^{\alpha u(x)}\in L^1(\Omega)$ and $\Omega$ is bounded, the limit of \eqref{st:eq:rt} is 0 as $T\rightarrow+\infty$. Now arguing as in the proof of the conclusion (1) we will show the conclusion (2). The proof is complete. 
\ep
  \section{Miranda's observation}
  In this section we study the relationship between two area functionals and the perimeter in corresponding manifolds. In Theorem \ref{key:thm} we show that the Miranda's observation mentioned in the introduction is  true for the conformal product functional $\AF_\alpha(u,\Omega)$ in $Q_\alpha$. 
  \subsection{The perimeter}
  We start with the definition of the perimeter. 
  \begin{Def}\label{Def:Cs} Let $\Omega$ be an open set in a Riemannian manifold $N$. Suppose $E$ is a Borel set in $N$ and $\lambda_E$ is its characteristic function. The perimeter of $E$ in $\Omega$ is given by 
  	\be 
  	\begin{split}
  		||D\lambda_E||_N(\Omega)
  		&=\sup\{\int_{\Omega}\lambda_E div(X)d\vf:X\in T_{0}\Omega, \la X,X\ra\leq 1\}
  	\end{split}
  	\ene 
  	If $E$ has locally finite perimeter, $||D\lambda_E||_N(\Omega')<\infty$ for each bounded open set $\Omega'$ in $\Omega$ (i.e. $\lambda_E\in BV_{loc, N}(\Omega)$), then $E$ is called a Caccioppoli set. 
  	\end{Def}
  \br\label{Remark:locbv:2} By Remark \ref{Remark:locbv} the fact that $E$ is a {\Ca} set only depends the conformal class of the metric in $N$. See also Lemma \ref{lm:conclass}. 
  \er 
  \indent Next we define the Hausdorff measure on Riemannian manifolds. 
  \begin{Def} 
   Let $N$ be a complete Riemannian manifold with dimension $n$. By the Nash Embedding Theorem, there is an isometric embedding $i:N\rightarrow \R^{n+s}$ for some positive integer $s$. Let $k>0$. For any set $E\subset N$, the $k$-dimensional Hausdorff measure of $E$, $H^k(E)$ is given by 
  	  $$
  	   H^k(E):= H^k(i(E))
  	  $$
  	  where in the right side $H^k$ denote the $k$-dimensional Hausdorff measure in $\R^{n+s}$. 
  \end{Def}
	A particular case is that for any open set $A$ in $N$, $H^n(A)=vol(A)$ according to the area formula where $vol(A)$ is the volume of $A$ with respect to the metric in $N$. For more details see the Example (2) of Subsection 8.6 in \cite{Simon83}.\\
	\indent  For a {\Ca} set its perimeter and other properties are unchanged if we make alternations of measure zero. In other words we are really concerned with the equivalence classes of a {\Ca} set. \\
  \indent The following result is a generalization of Proposition 3.1 in \cite{Giu84} with exactly the same proof by the Nash Embedding Theorem. 
   \begin{pro}\label{pro:evc} For any $x\in N$, $B_x(r)$ be the ball in $N$ centerred at $x$ with radius $r$. Let $inj(x)$ denote the injective radius of $x$,i.e the supremum of $r$ such that $B_r(x)$ is an embedded normal ball in $N$. \\
   	\indent If $E$ is a Borel set in  $N$, there exists a Borel set $\tilde{E}$ equivalent to $E$(that is, differs only by a set of $H^n$ measure zero) and such that 
   	\be \label{in:rep}
   	0<H^n(\tilde{E}\cap B_\rho(x))<vol(B_\rho(x)) 
    \ene
    for all $x\in \P\tilde{E}$ and $\rho\in (0,inj(x))$.
   	\end{pro}
  In the remainder of this paper we always assume \eqref{in:rep} holds for every {\Ca} set mentioned. 
  \subsection{The perimeter of subgraphs}
  First we show that the area functional of a BV function is the perimeter of its subgraph in corresponding product manifolds. \\
  \indent  The following result generalizes Theorem 14.6 in \cite{Giu84} in Euclidean space into general product manifold $Q$ and conformal product manifold $Q_\alpha$ (Definition \ref{Def:notation}).  
  \bt\label{perimeter} Use the notation in Definition \ref{Def:smooth}. Suppose $\Omega$ is an open bounded set in a Riemannian manifold $N$ with Lipschitz boundary. Let $u$ be a measurable function on $\Omega$ and $U$ be its subgraph. 
  \begin{enumerate} 
  	\item[(1)] If $u\in BV(\Omega)$, then 
$
\AF(u,\Omega)= ||D\lambda_U||_Q(\Omega\PLH\R)
$;
\item[(2)] If $\alpha>0$ and $e^{\alpha u}\in BV(\Omega)$, then
$ \label{conclusion:2A}
\AF_\alpha( u,\Omega)=||D\lambda_U||_{Q_\alpha}(\Omega\PLH\R)
$.
\end{enumerate}
  \et 
 First we obtain one side of the identites in Theorem \ref{perimeter} as follows.
  \bl\label{first:lm} Take the assumptions and notation in Theorem \ref{perimeter}. 
  \begin{enumerate}
  	\item Suppose $u\in BV(\Omega)$. Then 
  	$
  	||D\lambda_U||_{Q}(\Omega\PLH\R)\leq \AF(u,\Omega)
  	$. 
  	\item Suppose $\alpha >0$ and $e^{\alpha u}\in BV(\Omega)$. Then 
 $
||D\lambda_U||_{Q_\alpha}(\Omega\PLH\R)\leq \AF_\alpha (u,\Omega)
 $
  	\end{enumerate}
  \el  
  \bp  First assume $u\in C^1(\Omega)$. Then its subgraph $U$ has a $C^1$ boundary. According to the definition of the perimeter we have  
  \be 
  \begin{split}
  ||D\lambda_U||_{Q_{\alpha}}(\Omega\PLH\R)&=vol_{Q_\alpha}(\P U\cap (\Omega\PLH\R))\\
   ||D\lambda_U||_{Q}(\Omega\PLH\R)&=vol_Q(\P U\cap (\Omega\PLH\R))
   \end{split}
  \ene 
  where $vol$ denotes the volume and the lower index indicates the corresponding ambient manifold.   
  From the definition of $\AF(u,\Omega)$ and $\AF_\alpha(u,\Omega)$ we have 
  \begin{gather*}
vol_{Q_\alpha}(\P U\cap (\Omega\PLH\R))= \int_{\Omega}e^{\alpha u}\sqrt{1+|Du|^2}d\vf=\AF_\alpha(u,\Omega)\\
 vol_Q(\P U\cap (\Omega\PLH\R))= \int_{\Omega}\sqrt{1+|Du|^2}d\vf=\AF(u,\Omega)
  \end{gather*}
  Thus it holds that for $u\in C^1(\Omega)$ 
   \begin{gather}
||D\lambda_U||_{Q_{\alpha}}(\Omega\PLH\R)=\AF_\alpha(u,\Omega)\label{eq:c:1} \\ 
||D\lambda_U||_{Q}(\Omega\PLH\R)=\AF(u,\Omega)\label{eq:c:2}
  \end{gather}
  By (2) in Theorem \ref{thm:app}, for $u\in BV(\Omega)$ there exists a smooth sequence $\{u_i\}_{i=1}^\infty$ in $C^\infty(\Omega)$ such that $u_i$ converges to $u$ in $L^1(\Omega)$ and 
  $$
  \lim_{i\rightarrow \infty}\AF(u_i,\Omega)=\AF(u,\Omega)
  $$
    Let $U_i$ be the subgraph of $u_i$. It is easy to see that $\lambda_{U_i}$ converges to $\lambda_U$ in the $L^1_{loc}(Q)$. Thus with \eqref{eq:c:2} we obtain 
    \be 
  \begin{split}
  	||D\lambda_U||_{Q}(\Omega\PLH\R)&
  	\leq  \lim_{i\rightarrow\infty}\inf||D\lambda_{U_i}||_{Q}(\Omega\PLH\R)\\
  	&=\lim_{i\rightarrow +\infty }\inf\AF (u_i,\Omega)=\AF (u,\Omega)
  \end{split}
  \ene
  This gives the conclusion (1).\\
  \indent As for the conclusion (2) notice that if  $e^{\alpha u}\in BV(\Omega)$ then the definition of $\AF_\alpha(u,\Omega)$ implies that it is finite. By (3) in Theorem \ref{thm:app} there is a smooth sequence $\{u_i\}_{i=1}^\infty$ in $C^\infty(\Omega)$ such that $e^{\alpha u_i}$ converges to $e^{\alpha u}$ in $L^1(\Omega)$ and 
  \be \label{eq:alpha:convergence}
  \lim_{i\rightarrow \infty}\AF_\alpha( u_i,\Omega)=\AF_\alpha( u,\Omega)
  \ene 
 Thus with a similar argument as the proof of the conclusion (1) we conclude the conclusion (2).
   \ep 

Now we are ready to show Theorem \ref{perimeter}. Our proof is similar as that of Theorem 14.6 in \cite{Giu84}. 
  \bp Let $T$ be a fixed positive constant. First we assume $u\in [-T,T]$ and $u\in BV(\Omega)$. Suppose $h \in C_0(\Omega)$ and $X\in T_0(\Omega)$ satisfying 
  \be\label{condition:cw}
  e^{-2\alpha\F{(T+1)}{n}}(h^2+\la X,X\ra)\leq 1
  \ene 
  where $\la,\ra$ is the inner product on $N$.\\
  \indent Let $\eta(r)$ be a smooth function on $\R$ with its support in $[-(T+1),sup_{\Omega}u+1]$ such that $\eta\equiv 1$ in $[-T,\sup_{\Omega}u]$ and $|\eta(r)|\leq 1$. Let $\eta_1(r)$ be a smooth function with a compact support on $\R$ satisfying $\eta_1(r)=e^{-\alpha \F{r+T+1}{n}}$ on $[-(T+1),(T+1)]$.  Now we define a smooth vector field 
  \be
   X'=\eta_1(r)\eta(r)(h\P_r+ X)
  \ene 
  on $Q_\alpha$. 
  The inner product of $X'$ on $Q_\alpha$ is given by 
   \be 
   \la X',X'\ra_{g_\alpha}=e^{-2\alpha \F{(T+1)}{n}}\eta^2(r)(h^2+\la X,X\ra)\leq 1
   \ene 
   where $g_\alpha$ denotes  the metric of $Q_\alpha$. \\
 \indent Let $d\vf_\alpha$ and $d\vf$ be the volume form of $Q_{\alpha}$ and $\Omega$. We have 
   \be \label{eq:d:1}
   d\vf_{\alpha}= e^{\alpha\F{(n+1)r}{n}}d\vf dr 
   \ene 
   Let $div_\alpha$ and $div$ be the divergence of $Q_{\alpha}$ and $\Omega$ respectively. Notice that $X'$ has compact support in $\Omega\PLH\R$. By the definition of the perimeter we have 
   \be \label{middle:step:one}
   ||D\lambda_U||_{Q_\alpha}(\Omega\PLH\R)\geq \int_{\Omega\PLH\R}\lambda_U div_\alpha(X')d\vf_{\alpha}
   \ene 
   From the definition of $\eta_1(r)$ and $\eta(r)$ we have 
   \be \label{eq:d:3}
   X'\llcorner d\vf_{\alpha}=e^{-\alpha\F{(T+1)}{n}}\eta(r)e^{\alpha r}((-1)^nh(x)d\vf+X\llcorner d\vf dr)
   \ene 
   This yields that 
   \be\label{eq:d:4}
   \begin{split}
   div_{\alpha}(X')d\vf_{\alpha}&=d(X'\llcorner d\vf_{\alpha})\\
    &=e^{-\alpha\F{(T+1)}{n}}(\eta(r)e^{\alpha r})'h(x)d\vf dr\\
    &+e^{-\alpha\F{(T+1)}{n}}e^{\alpha r}\eta(r)div(X)d\vf dr
    \end{split}
   \ene
   Notice that 
   \be
   \int_{-\infty}^{u(x)}(\eta(r)e^{\alpha r})'dr=e^{\alpha u(x)}
   \ene 
   and 
   \be
   \int_{-\infty}^{u(x)}\eta(r)e^{\alpha r }dr=\left\{\begin{split}
   	&\F{ e^{\alpha u(x)}}{\alpha}+C, \quad \text{ $\alpha\neq 0$}\\
   	&u(x)+C,\quad \text{ $\alpha =0$}
\end{split}\right.
  \ene 
   where $C$ is a fixed constant. Now $F_{in}:\equiv \int_{\Omega\PLH\R}\lambda_U div_\alpha(X')d\vf_{\alpha}$. Therefore 
   \be
   F_{in}=\left\{\begin{split}
   	&\int_\Omega e^{-\alpha\F{(T+1)}{n}}(e^{\alpha u(x)}h(x)+\F{e^{\alpha u(x)}}{\alpha}div(X))d\vf , \quad \text{ $\alpha\neq 0$}\\
   	&\int_\Omega (h(x)+u(x)div(X))d\vf,\quad \text{ $\alpha = 0$}
   	\end{split}\right.
   \ene
   From \eqref{middle:step:one} and the condition \eqref{condition:cw}, we obtain 
   \be\label{eq:step:three}
   ||D\lambda_U||_{Q_\alpha}(\Omega\PLH\R) \geq \sup_{h,X\in\eqref{condition:cw}}F_{in}=\left\{\begin{split}
    	&\AF_\alpha(u,\Omega) , \quad \text{ $\alpha> 0$}\\
    	&\AF(u,\Omega),\quad \text{ $\alpha = 0$}
    \end{split}\right.
   \ene 
 for $|u|\leq T$. Combining Lemma \ref{first:lm} and \eqref{eq:step:three} together we obtain Theorem \ref{perimeter} in the case that $u$ is uniformly bounded.\\
\indent For $T>0$ we define    $$u_T(x):=\max\{\min(u(x),T), -T\}$$
Now suppose $u\in BV(\Omega)$, then $u\in L^1(\Omega)$. Thus $u_T$ converges to $u(x)$ in $L^1(\Omega)$ as $T\rightarrow +\infty$. By Theorem \ref{thm:semi} and  (1) in Lemma \ref{lm:appr}  we have   
   \be\label{conclusion:second}
   \begin{split}
   ||D\lambda_U||_Q(\Omega\PLH\R) &= \lim_{T\rightarrow+\infty}||D\lambda_{U_T}||_Q(\Omega\PLH\R)\\
   &=\lim_{T\rightarrow +\infty}\AF(u_T,\Omega)\geq \AF(u,\Omega) 
   \end{split} 
   \ene
   With the conclusion (1) of Lemma \ref{first:lm} we obtain 
     $$
       ||D\lambda_U||_Q(\Omega\PLH\R)=\AF(u,\Omega)
     $$
     for $u\in BV(\Omega)$. This is the conclusion (1).  \\
   \indent Now assume $e^{\alpha u(x)}\in BV(\Omega)$. Thus $\AF_\alpha(u,\Omega)$ is finite. Moreover $\Omega$ is bounded.  By Theorem \ref{thm:semi} and  (2) in Lemma \ref{lm:appr} we conclude that 
   \be\label{conclude}
     \begin{split}
   		 ||D\lambda_U||_{Q_\alpha}(\Omega\PLH\R) &= \lim_{T\rightarrow+\infty}||D\lambda_{U_T}||_{Q_\alpha}(\Omega\PLH\R)\\
   	&=\lim_{T\rightarrow +\infty}\AF_\alpha(u_T,\Omega)\geq \AF_\alpha(u,\Omega) 
   \end{split} 
   \ene 
     With the conclusion (2) of Lemma \ref{first:lm} we conclude 
     $$
      ||D\lambda_U||_{Q_\alpha}(\Omega\PLH\R)=\AF_\alpha(u,\Omega)
     $$  
     This is the conclusion (2). The proof is complete.
  \ep 
  Now we obtain a relationship between $\AF_\alpha(u,\Omega)$ and $||Du||_N(\Omega)$ as follows. 
  \begin{cor} \label{lm:est}Suppose $u$ is a measurable function on $\Omega$ with $|u|\leq T$ such that $\AF_\alpha(u,\Omega)$ is finite. Then 
  	\be 
  	\AF_\alpha(u,\Omega)\geq e^{-\alpha T}\max\{vol(\Omega), ||Du||_N(\Omega)\}  	\ene 
  	\end{cor}
  \bp  By Theorem \ref{thm:BV} and Theorem \ref{perimeter}, we have 
  \begin{align*}
   \AF_\alpha(u,\Omega)&=||D\lambda_U||_{Q_\alpha}(\Omega\PLH\R)\\
    &=\int_{\Omega\PLH\R}e^{\alpha r}d|D\lambda_U|_Q\\
    &\geq e^{-\alpha T}||D\lambda_U||_Q(\Omega\PLH\R)\\
    &=e^{-\alpha T}\AF(u,\Omega)     \end{align*}
    In the third line above, we use the fact that the support of the Radon measure $|D\lambda_U|_Q$ is on the set $(x,u(x))$. The conclusion follows from  (3) in Theorem \ref{thm:semi}. 
  \ep 
 
  \subsection{Miranda's observation} 
  In this subsection we show that a local minimizer of the conformal area functional $\AF_\alpha(u,\Omega)$ is a local minimizer of the perimeter of the perimeter in $Q_\alpha$. According to Giusti \cite{Giu80} it is Miranda \cite{Mir64} firstly to observe this phenomenon for product area functional in the product manifold.  \\
  \indent The following result is similar to Theorem 14.8 in \cite{Giu84}.
  \bl \label{key:lemma}Let $U$ be the subgraph of $u(x)$ in $Q_\alpha$. Let $\Omega$ be an open set in $N$.  Let $F\subset Q_\alpha$ be a {\Ca} set with finite perimeter  satisfying for a.e. $x\in \Omega$,
  $\lambda_F(x,t)=0$ for all $t>T_x$ where $T_x$ is a constant depending on $x$.\\
  \indent Then the function $\omega(x)$ 
   \be\label{eq:omega:def}
   e^{\alpha\omega(x)}=\alpha\lim_{k\rightarrow+\infty}(\int_{-k}^ke^{\alpha t}\lambda_{F}(x,t)dt)
   \ene 
   is a.e. well-defined and 
  \be\label{conclusion:es}
  \AF_\alpha(\omega,\Omega)\leq ||D\lambda_F||_{Q_\alpha}(\Omega\PLH\R)
  \ene 
  \el 
  \bp Due to the assumption in Lemma \ref{key:lemma}, it is obvious that $\omega(x)$ is well-defined. Suppose $h\in C_0(\Omega) $ and $X\in T_0(\Omega)$ satisfying 
  \be \label{condition:lessthanone}
  h^2+\la X,X\ra\leq 1
  \ene 
  where $\la,\ra$ is the inner product with respect to the metric $g$.  \\
  \indent Let $\eta(r)$ be a smooth function such that $0\leq \eta(r)\leq 1$ with compact support in $\R$. 
  Set $X'=e^{-\F{\alpha r}{n}}\eta(r)(X+h(r)\P_r)$. Then $\la X',X'\ra_\alpha\leq 1$ where  $\la, \ra_\alpha $ denotes the inner product of $Q_\alpha$. By Definition \ref{Def:BV} we have 
  \be \label{deq:st}
||D\lambda_F||_{Q_\alpha}(\Omega\PLH\R)\geq \int_{\Omega\PLH\R}\lambda_F(x,r)div_\alpha(X')d\vf_\alpha
   \ene 
   where $div_\alpha$ and $d\vf_\alpha$ are the divergence and the volume form of $Q_\alpha$.  Arguing as in \eqref{eq:d:1},\eqref{middle:step:one}, \eqref{eq:d:3} and \eqref{eq:d:4}  we obtain 
   \begin{align*}
   div_\alpha(X')d\vf_\alpha 
            &=\{(e^{\alpha r}\eta(r))'h(x)+e^{\alpha r}\eta(r)div(X)\}d\vf dr
   \end{align*}   where $div$ is the divergence of $\Omega$. Thus expanding \eqref{deq:st} gives that
   \be
   \begin{split}
   ||D\lambda_F||_{Q_\alpha}(\Omega\PLH\R)&\geq \int_{\Omega}h(x)\{\int_{-\infty}^\infty (e^{\alpha r}\eta(r))'\lambda_F(x,r) dr\}d\vf\\
   &+\int_\Omega div(X)\{\int_{-\infty}^\infty e^{\alpha r}\eta(r)\lambda_F(x,r)dr\} d\vf 
   \end{split}
   \ene 
 Replacing $\eta(t)$ with a sequence of $\{\eta_k(t): 0\leq \eta_k\leq 1\text{ with compact support}\}_{k=1}^\infty$ which converges to the constant function $1$ as $k\rightarrow +\infty$, we obtain 
   \be
 \begin{split}
 	||D\lambda_F||_{Q_\alpha}(\Omega\PLH\R)&\geq \int_{\Omega}h(x)\{\int_{-\infty}^\infty \alpha e^{\alpha r}\lambda_F(x,r) dr\}d\vf\\
 	&+\int_\Omega div(X)\{\int_{-\infty}^\infty e^{\alpha r}\lambda_F(x,r)dr\} d\vf \\
 	&=\int_{\Omega}e^{\alpha \omega(x)}(h(x)+\F{1}{\alpha}div(X))d\vf
 \end{split}
 \ene 
 Here we apply the fact for a.e. $x\in \Omega$,
 $\lambda_F(x,t)=0$ for all $t>T_x$ where $T_x$ is a constant depending on $x$. 
   Now taking the supremum of all $h$ and $X$ satisfying \eqref{condition:lessthanone} and applying the defintion of $\AF_\alpha(.,\Omega)$, one sees that
   \be
     	||D\lambda_F||_{Q_\alpha}(\Omega\PLH\R)\geq \AF_\alpha (\omega,\Omega )
     	\ene
    The proof is complete. 
  \ep Now we can conclude the Miranda's observation in the conformal product manifold $Q_\alpha$ given in Definition \ref{Def:notation} as follows. 
  \bt[\textbf{Miranda's observation}]\label{key:thm} Let $\Omega\subset N$ be an open bounded set and $\alpha>0$ be a fixed positive constant. Let $u$ be a measurable function. If $e^{\alpha u}\in BV(\Omega)$ and $u$ is a local minimum of the conformal area functional $\AF_\alpha(u,\Omega)$. Let  
  $ U$ be the subgraph of $u(x)$. Then $U$ locally minimizes the perimeter in the manifold $Q_\alpha$. 
  \et
  \bp  Suppose $u(x)$ is a measurable function on $\Omega$ with $e^{\alpha u(x)}\in BV(\Omega)$ and $u(x)$ is a local minimum of $\AF_\alpha(u,\Omega)$. Let $U$ denote the subgraph of $u(x)$ in $Q_\alpha$. \\
  \indent Let $F$ be a {\Ca} set with finite perimeter satisfying $F\Delta U\subset K$ where $K$ is a compact set in $\Omega\PLH\R\subset Q_\alpha$. Thus $F$ has the property that for a.e. $x\in \Omega$,
  $\lambda_F(x,t)=0$ for all $t>T_x$ where $T_x$ is a constant depending on $x$. Now from Lemma1 \ref{key:lemma} we obtain 
  \be \label{fact:first}
  \AF_\alpha (\omega, \Omega)\leq ||D\lambda_F||_{Q_\alpha}(\Omega\PLH\R)
  \ene
  where $\omega(x)$ is given by 
  $$e^{\alpha\omega(x)}=\alpha\lim_{k\rightarrow+\infty}(\int_{-k}^ke^{\alpha t}\lambda_{F}(x,t)dt)$$
  Since  $F\Delta U\subset K$, then $\omega(x)=u(x)$ outside a compact set in $\Omega$. By Theorem \ref{perimeter} one sees that 
  \be \label{fact:second}
 ||D\lambda_U||_{Q_\alpha}(\Omega\PLH\R)=\AF_\alpha(u,\Omega )\leq \AF_\alpha (\omega, \Omega)
  \ene
  Here we apply the fact that $u(x)$ is a local minimum of $\AF_\alpha(.,\Omega)$. Combining \eqref{fact:first} and \eqref{fact:second} together we obtain 
  \be
  ||D\lambda_U||_{Q_\alpha}(\Omega\PLH\R)\leq ||D\lambda_F||_{Q_\alpha}(\Omega\PLH\R)
  \ene
  The proof of Theorem \ref{key:thm} is complete. 
  \ep
  \section{Existence of generalized solutions}
  In this section we define  a generalized solution of the  Dirichlet problem \eqref{main:equation} following from the Miranda-Giusti's generalized solution theory in \cite{Giu84}.  
  \subsection{A generalized solution} We continue to use the notation in Definition \ref{Def:notation}. 
  Now fix two bounded open sets $\mathfrak{B}, \Omega$ in $N$ satisfying $\bar{\Omega}\subset\subset \mathfrak{B}$. Moreover $\Omega$ has Lipschitz boundary $\P\Omega$.
  \begin{Def}\label{Def:gs}
  	Let $\Omega\subset \subset\mathcal{B}$ be two open bounded set in $N$. Let $u(x), \psi(x)$ be two measurable functions taking values in $[-\infty, \infty]$ such that the subgraph $\Psi$ of $\psi(x)$ is a {\Ca} set in $Q_\alpha$. 
  	We say that $u(x)$ is a generalized solution to the {\DP} \eqref{main:equation} on $\Omega$ with boundary data $\psi(x)$ on $\P\Omega$ if
  	\begin{itemize} 
  		\item [(1)] the subgraph of $u(x)$, $U$,  coincides with $\Psi$ outside $\bar{\Omega}\PLH\R$ ;
  		\item [(2)]for any {\Ca} set $F$ satisfying $F\Delta U\subset K$ where $K$ is a compact set in $\bar{\Omega}\PLH \R$, it holds that 
  		\be \label{eq:minimize}
  		\int_{K}d|D\lambda_U|_{Q_\alpha} \leq  \int_{K}d|D\lambda_F|_{Q_\alpha}
  		\ene 
  	\end{itemize}
  	where $|D\cf_U|_{Q_\alpha}$ and $|D\cf_F|_{Q_\alpha}$ are Radon measures of generated by $\lambda_U$ and $\lambda_F$ in $Q_\alpha$ respectively (See Theorem \ref{thm:st}). 
  \end{Def} 
\br  In the followings such $u(x)$ is also referred as a generalized solution if without any confusion. 
\er 
In order to understand \eqref{eq:minimize}, we need to define the generalized trace of $u$ when its subgraph $U$ is a {\Ca} set in $Q_\alpha$. It is useful in the proof of Lemma \ref{lm:minimize}. 
\bl\label{lm:trace:ca} Suppose $\Omega$ is an open bounded set in $N$ with Lipschitz boundary. Let $u(x)$ be a measurable function on $N$ such that its subgraph $U$ is a {\Ca} set in $Q_\alpha$. Then 
\begin{enumerate}
	\item Let $u_k(x)$ be the function $\max\{\min(u(x),k), -k\}$.  Then $\{u_k(x) \}\in BV(\Omega)$ and their trace $\{ u_k(x)\}$  converges to a measurable function $\mg u(x)$ taking value in $[-\infty,\infty]$ a.e. on $\P\Omega$ as $k$ goes to $+\infty$. 
	\item Let $U_{\mg}$ be the subgraph of $\mg u(x)$ in $\P\Omega\PLH\R$. The trace of $\lambda_{U}$ on $\P\Omega\PLH\R$ satisfies that  
  $$
   \mT\lambda_U=\lambda_{U_{\mg}}  $$
   \end{enumerate}
\el 
\begin{Def}\label{def:g:trace}
When the subgraph of $u(x)$ is a {\Ca} set in $\Omega\PLH\R\subset Q_\alpha$, $\mg u(x)$ defined in the above lemma is called the generalized trace of $u(x) $ from $\Omega$. 
\end{Def}
\br If $u\in BV(\Omega)$ and $\Omega$ is Lipschitiz, the generalized trace of $u(x)$ is just the trace of $u$ by Remark \ref{rk:bounded:approximation} and Lemma \ref{lm:pro:trace}. In general the generalized trace of $u(x)$ can take the value in $[-\infty,\infty]$ and may not belong to $L^1(\P\Omega)$. For example see Remark \ref{rk:trace:infty}. 
\er 
\bp Let $U_k$ be the subgraph of $u_k(x)$ in $\mathcal{B}\PLH\R$. Thus $\lambda_{U_k}=\lambda_U$ on the open set $\Omega_k$ in $Q_\alpha$. Since $U$ is a {\Ca} set in $Q_\alpha$, then 
\begin{align}
\infty>||D\lambda_U||_{Q_\alpha}(\Omega\PLH(-k,k))&=||D\lambda_{U_k}||_{Q_\alpha}(\Omega\PLH(-k,k))\label{eq:mid:key}\\
 &=\AF_\alpha(u_k,\Omega)\quad\text{by Theorem \ref{perimeter}}\notag \\
 &\geq e^{-\alpha k}||Du_k||_N(\Omega) \quad\text{by Corollary \ref{lm:est}}\notag
\end{align}
Thus $u_k\in BV(\Omega)$. By Theorem \ref{thm:trace}, its trace $\mT u_k$ is well-defined. Let $U_k$ be the subgraph of $u_k$. Let $
\mT\lambda_{U_k}$ be the trace of $U_k$ on $\P\Omega\PLH\R\subset Q_\alpha$. By Lemma \ref{lm:cbv}  there is a measure zero set $E$ in a open interval $(0,r_0)$ such that 
\begin{gather}
\lim_{r\rightarrow 0^+, r\notin E}u_{k}(x,r)=\mT u_k(x) \quad x\in \P\Omega \text{ in  }  L^1(\P\Omega)\\
\lim_{r\rightarrow 0^+,r\notin E}\lambda_{U_k}(p,r)
=\mT\lambda_{U_k}(p)\quad p\in \P\Omega\PLH\R \text{ in  }  L^1(\P\Omega\PLH(-k,k))\
\end{gather}
Thus $\mT\lambda_{U_k}$ is the subgraph of $\mT u_k$ on $\P\Omega$ in $\P\Omega\PLH (-k,k)$. \\
 \indent By \eqref{eq:mid:key} and the definition of $U_k$ it is clear that $\mT\lambda_{U_k}$ converges to $\mT\lambda_U$ almost everywhere on $\P\Omega\PLH\R$ as $k\rightarrow +\infty$. Thus arguing as Lemma 16.3 in \cite{Giu84}, $
 \mT\lambda_{U} $ on $\P\Omega\PLH\R$ is a subgraph of some measurable function $\mg u(x)$ on $\P\Omega$. Thus  $\mT u_k(x)$ converges to $\mg u(x)$ a.e. on $\P\Omega$. 
\ep
 
\subsection{The existence of  generalized solutions} The main result of this section is stated as follows.
 \bt\label{mt:A}Let $\Omega\subset \subset\mathcal{B}$ be two open bounded set in a  Riemannian manifold $N$. Suppose  $\P\Omega$ is Lipschitz.  Let $\psi(x)$ be any measurable function such that its subgraph $\Psi$ is a {\Ca} set in $Q_\alpha$. Then there is a generalized solution to the Dirichlet problem \eqref{main:equation} with boundary data $\psi(x)$.  
\et 
\br \label{re:mk} The restriction on $\psi(x)$ is sufficient general. Suppose $\psi(x)\in L^1(\P\Omega)$ and $\P\Omega$ is Lipschitz. By Lemma \ref{lm:st:extend} there is a function, still denoted by $\psi(x)$, in $BV(\mathcal{B})$ such that its trace on $\P\Omega$ from $\Omega$ and $\mathcal{B}\backslash \bar{\Omega}$ are $\psi(x)$. Since $\mathcal{B}$ is bounded, the conclusion (3) in Theorem \ref{thm:semi} implies that $\AF(\psi,\mathcal{B})$ is finite. From the conclusion (1) in Theorem \ref{perimeter} $\Psi$ is a {\Ca} set in $Q$. Since $Q$ and $Q_\alpha$ are conformal to each other (see Definition \ref{Def:notation}), then $\Psi$ is also a {\Ca} set in $Q_\alpha$ by Remark \ref{Remark:locbv} and Remark \ref{Remark:locbv:2}.  In particular if $\psi(x)$ is continuous on $\P\Omega$ and bounded, then the extension function is still bounded by Lemma \ref{lm:st:extend}. 
\er
  \bp  For any $k>0$ we set 
  $$
  \psi_k(x)=\min\{k,\max\{\psi(x),-k\}\}
   $$
 Let  $\Psi_{k}$ be the subgraph of $\psi_k$. Because $\lambda_{\Psi_k}=\lambda_{\Psi}$ on $\Omega\PLH(-k,k)$. On the other hand since $\Psi$ is a {\Ca} set in $Q_\alpha$, then $$||D\lambda_{\Psi}||_{Q_\alpha}(\Omega\PLH(-k,k))=||D\lambda_{\Psi_k}||_{Q_\alpha}(\Omega\PLH(-k,k))<\infty$$ By (2)  in Theorem \ref{perimeter} this implies  $\AF_\alpha(\psi_k,\Omega)$ is finite. By Corollary  \ref{lm:est} $\psi_k\in BV(\mathfrak{B})$ because $|\psi_k(x)|\leq k$. Thus we can consider the following minimizing problem 
   \be
   \alpha_k=\min\{ \AF_\alpha (u,\mathfrak{B}): u\in BV(\Omega), |u|\leq k, u=\psi_k\text{ on } \mathfrak{B}\backslash\bar{ \Omega} \}
\ene 
Let $\{u_j\}_{j=1}^\infty\in BV(\Omega)$ be the sequence satisfying $|u_j|\leq k, u=\psi_k$ on $\mathfrak{B}\backslash\Omega$ such that 
   $$
  \lim_{j\rightarrow+\infty} \AF_\alpha( u_j,\mathfrak{B})=\alpha_k
   $$
  Again by Corollary \ref{lm:est} we have 
                  $$
               \max\{||Du_j||_N(\mathfrak{B}):j=1,\cdots,\infty\}\leq C(k,\alpha_k)
                  $$
Since the boundary of $\mathcal{B}$ is Lipshcitz, the compactness of BV functions implies that there is a subsequence of $\{u_j\}_{j=1}^\infty$ with $|u_j|\leq k$ and $u_j=\psi_k$ on $\mathfrak{B}\backslash\Omega$ , still denoted as $\{u_i\}_{j=1}^\infty$, such that $u_j\rightarrow u_k$ in $L^1(\mathfrak{B})$ as $j\rightarrow +\infty$. By the semicontinuity of $\AF_\alpha(u,\mathfrak{B})$, we have 
 \be 
 \alpha_k\leq \AF_\alpha(u_k,\mathfrak{B})\leq \lim_{j\rightarrow \infty}\inf\AF_\alpha ( u_j,\mathfrak{B})=\alpha_k
 \ene    
 with the property that 
 \be \label{eq:property}
  |u_k|\leq k\quad u_k=\psi_k \text{\quad on $\mathfrak{B}\backslash\Omega$ }
 \ene 
 Thus $\AF_\alpha(u_k,\mathfrak{B})=\alpha_k$. \\
 \indent  Now let $U_k$ be the subgraph of $u_k$ in $\mathfrak{B}\PLH\R\subset Q_\alpha$. Let  $W$ be the set $U_k\cup\Omega\PLH(-T,T)$. Suppose $k>2T$. By  Theorem \ref{key:thm} we obtain the estimate 
 \begin{align*}
 \int_{\mathcal{B}\PLH(-2T, 2T)}d|D\lambda_{U_k}|_{Q_\alpha}&\leq  \int_{\mathcal{B}\PLH(-2T,2T)}d|D\lambda_{W}|_{Q_\alpha}\\
 &\leq \int_{\mathcal{B}\PLH(-2T,2T)}d|D\lambda_{\Psi}|_{Q_\alpha}+2vol(\P (\Omega\PLH(-2T,2T)))\\
 &=c(T)
\end{align*}
Arguing as Lemma 16.3 in \cite{Giu84} we can extract a subsquence, still denote again by $u_k$, converging almost everywhere to a measurable function $u_\infty$.\\
 \indent 
  Now let $U_\infty$ be the subgraph of $u_\infty$. It is clear that $U_\infty$ coincides with $\Psi$ outside $\bar{\Omega}\PLH\R$.  Suppose $V$ is 
  a {\Ca} set in $\mathcal{B}\PLH\R$ coinciding with $U_\infty $ except some compact set $K\subset\bar{\Omega}\PLH(-T,T)$. Let $A$ be an open set satisfying $\Omega\subset \subset A\subset \subset \mathcal{B}$. Now set \be 
  V_k=\left\{\begin{split}
  	& V\text{ in } A\PLH(-T,T)\\
  	& U_k\text{ outside } A\PLH(-T,T)
  	\end{split}\right.
  \ene   
  By Theorem \ref{key:thm} $U_k$ locally minimizes the perimeter in $\mathcal{B}\PLH(-T,T)$. 
  \be\label{eq:uvst}
  \int_{\bar{A}\PLH [-T,T]}d|D\lambda_{U_k}|_{Q_\alpha}\leq \int_{\bar{A}\PLH[-T,T]}d|D\lambda_{V_k}|_{Q_\alpha}
  \ene 
 If $k>T$, it is clear that $V_k$ coincides with $U_k$ outside $\bar{\Omega}\PLH(-T,T)$.  By Theorem \ref{thm:bd:trace} and the definition of traces of BV functions we have  
   $$
   \int_{(\P A)\PLH (-T,T)}d|D\lambda_{U_k}|_{Q_\alpha}= \int_{(\P A)\PLH(-T,T)}d|D\lambda_{V_k}|_{Q_\alpha}
   $$
 Notice that  $\P(\bar{A}\PLH[-T,T])=\{(\P A)\PLH(-T,T)\}\cup (A\PLH\{\pm T\})$. By Theorem \ref{thm:bd:trace} the inequality \eqref{eq:uvst} is equivalent to  
  \be\label{eq:key:st}
  \begin{split}
  \int_{A\PLH(-T,T)}d|D\lambda_{U_k}|_{Q_\alpha}&\leq \int_{A\PLH(-T,T)}d|D\lambda_{V}|_{Q_\alpha}\\
  &+\int_{A\PLH\{\pm T\}}|\mT\lambda_{V}-\mT\lambda_{U_k}|d\vf_{Q_\alpha}
  \end{split}
  \ene
  where $\mT$ denotes the corresponding trace of $\lambda_V$ and $\lambda_{U_k}$  on from $A\PLH(-T,T)$ and $A\PLH(T,\infty)$ respectively.\\
  \indent Because $U_k$ and $U_\infty$ are {\Ca} set in $Q_\alpha$ for all $k=1,\cdots,\infty$. By the Fubiniz theorem for a.e. $T\in \R$  we have 
  \be\label{eq:one:step}
  \int_{A\PLH\{T\}}d|Df|_{Q_\alpha}=0\quad\text{ for\quad } f=\lambda_{U_k} \text{\quad or\quad} \lambda_{U_\infty}
  \ene 
  This implies that for a.e. $T$
  \be 
  \mT f(x, \pm T)=f(x,\pm T)\quad\text{ for\quad } f=\lambda_{U_k} \text{\quad or\quad} \lambda_{U_\infty}\quad a.e. x\in A
  \ene 
   and $\mT \lambda_V(x,\pm T)=\lambda_{U_\infty}(x,\pm T)$ for $a.e. x\in A$. Since $\lambda_{U_k}$ converges to $\lambda_{U_\infty}$ locally with the $L^1(Q)$ norm,  
   \be\label{eq:two:step}
   \lim_{k\rightarrow+\infty} \lambda_{U_k}(x,T)=\lambda_{U_\infty}(x,T)
   \ene  for a.e. $T\in \R$ and $x\in A$ as $k\rightarrow +\infty$. Now we choose a $T\in \R$ satisfying in \eqref{eq:one:step}, \eqref{eq:two:step} for all $k=1,\cdots, \infty$ and the compact set $K\subset \bar{\Omega}\PLH[-T,T]$. Let  $k\rightarrow \infty$ in \eqref{eq:key:st}. By Theorem \ref{thm:semi} it becomes 
  \be 
  \int_{A\PLH(-T,T)}d|D\lambda_{U_\infty}|_{Q_\alpha}\leq \int_{A\PLH(-T,T)}d|D\lambda_V|_{Q_\alpha}
  \ene 
  which gives 
\be 
\int_{K}d|D\lambda_{U_\infty}|_{Q_\alpha}\leq \int_{K}d|D\lambda_V|_{Q_\alpha}
\ene    
From the assumption of $V$ we obtain that $u$ is a generalized solution to the Dirichlet problem \eqref{main:equation} on $\Omega$ with boundary data $\psi(x)$. This completes the proof.   
  \ep 
  
   \section{The infinity value of generalized solutions} In this section we study generalized solutions that take the infinity values.  First we construct generalized solutions only taking inifinity values. Then we recall some facts on the almost minimal set. With its regularity theory we describe the boundary property of  the set in which a generalized solution takes the infinity value in Theorem \ref{thm:infty}.
  \subsection{Generalized solutions with infinity boundary data}
   We continue to use the notation in Definition \ref{Def:notation}. First in $Q_\alpha$ the translating motion of generalized solutions also gives generalized solutions. 
  \bl \label{lm:tran} Take the notation in Definition \ref{Def:gs}. Suppose $u(x)$ is a generalized solution with boundary data $\psi$.and $a\in \R$.  Then  $u(x)+a$ is a generalized solution  on $\Omega$ with boundary data $\psi(x)+a$ in $Q_\alpha$. 
  \el 
  \bp Fix $a\in \R$. We define a translating motion $T_a: N\PLH\R\rightarrow N\PLH\R$ as $T_a(x,t)=(x,t+a)$ where $x\in N$ and $t\in \R$. It is clear $T_a$ is an isometry of $Q$ with respect to the metric $g+dr^2$. \\
  \indent Suppose $F$ is  a {\Ca} set in $Q$. For any smooth vector field $X$ with compact support  and any Borel set $C$ in $Q$, we have 
  \be \label{eq:average}
  \int_{C}\lambda_F div_Q(X)d\vf_Q=\int_{T_{a}C}\lambda_{T_aF}div_Q(T^*_a X)d\vf_Q
  \ene 
 where $div_Q$ and $\vf_Q$ are the divergence and the volume form of $Q$ respectively.    By Theorem \ref{thm:st} for any Borel set $C$ in $Q$,  
  \eqref{eq:average} implies that 
  \be\label{eq:mist}
  \int_{C}d|D\lambda_{F}|_Q=\int_{T_{a }C}d|D\lambda_{T_aF}|_Q
  \ene  
  Recall that the metric of  $Q_\alpha$ is $e^{2\F{\alpha}{n}}(g+dr^2)$. Combining Theorem \ref{thm:BV} with \eqref{eq:mist} together we obtain 
  \be\label{eq:standard}
  e^{\alpha a}\int_{C}d|D\lambda_{ F}|_{Q_\alpha} =\int_{T_{a}C}d|D\lambda_{T_a F}|_{Q_\alpha}
  \ene 
  Here the function $\psi^{m-1}$ in Theorem \ref{thm:BV}  is just $e^{\alpha r}$.\\
  \indent Let $U$ be the subgraph of $u(x)$. Let $F$ be a
   {\Ca} set in $Q_\alpha$ satisfying $F\Delta T_aU\subset K$ where $K$ is a compact set in $\bar{\Omega}\PLH\R$. Thus $T_{-a}F\Delta U\subset T_{-a}K$. Because $u(x)$ is a generalized solution on $\Omega$, by Definition \ref{Def:gs} we have 
  \begin{align*}
  \int_{K}d|D\lambda_{T_aU}|_{Q_\alpha}&=e^{\alpha a}\int_{T_{-a}K}d|D\lambda_U|_{Q_\alpha}\\
  &\leq e^{\alpha a}\int_{T_{-a}K}d|D\lambda_{T_{-a}F}|_{Q_\alpha}\\
  &=  \int_{K}d|D\lambda_{F}|_{Q_\alpha}\quad \text{   by \eqref{eq:standard}}
  \end{align*}
  Obviously outside $\bar{\Omega}\PLH\R$ the subgraph of $u+a$ is just the subgraph of $\psi(x)+a$. By Definition \ref{Def:gs} $u(x)+a$ is the generalized solution on $\Omega$ with boundary data $\psi(x)+a$. The proof is complete. 
  \ep
   Suppose $u(x)$ is a generalized solution on $\Omega$  with boundary data $\psi(x)$ as in Definition \ref{Def:gs}.  We set 
  \begin{gather*}
  P_+=\{x\in\Omega:u(x)=+\infty\}\label{pinfty}\\
  P_-=\{x\in\Omega:u(x)=-\infty\}\label{ninfty}
  \end{gather*} 
 Now we give more assumptions on the boundary data $\psi(x)$. Assume that there are three open sets $ A_{\pm}, A_0$ in $\mathfrak{B}$ with the following properties:
  \begin{enumerate}
  	\item $\psi(x)=\pm \infty$ in $A_{\pm}$; 
  	\item $\psi(x)$ is continuous on $A_0$; 
  	\item $\P\Omega=\Gamma_0\cup\Gamma_{+}\cup
  	\Gamma_{-}\cup \mathfrak{N}$ where $\Gamma_0=\P\Omega \cap \P A_0,\Gamma_{\pm}=\P\Omega\cap \P A_{\pm}$ and $H^{n-1}(\mathfrak{N})=0$.
  	\item $\Gamma_{+}$ and $\Gamma_{-}$ are $C^2$-hypersurfaces. 
  	\end{enumerate}
  \bl \label{lm:minimize} Let $u(x)$ be a generalized solution to the Dirichlet problem \eqref{main:equation} on $\Omega$ with boundary data $\psi(x)$ on $\P\Omega$. Then it holds that 
  \begin{enumerate}
  	\item the function $u^+_{\infty}$ given by 
  $$
    u^+_\infty=+\infty \text{  on $P_+$ } \quad     u^+_\infty=-\infty \text{  on $\mathcal{B}\backslash P_+$ }  
  $$
  is a generalized solution to the Dirichlet problem \eqref{main:equation} on $\Omega$ with boundary data 
   $$
  \psi^+_\infty=+\infty \text{  on $A_+$ } \quad     \psi^+_\infty=-\infty \text{  on $\mathcal{B}\backslash A_+$ }  
  $$
  	\item the function $u^-_{\infty}$ given by 
  $$
  u^-_\infty=-\infty \text{  on $P_-$ } \quad     u^-_\infty=+\infty \text{  on $\mathcal{B}\backslash P_-$ }  
  $$
  is a generalized solution to the Dirichlet problem \eqref{main:equation} on $\Omega$ with boundary data 
  $$
  \psi^-_\infty=-\infty \text{  on $A_-$ } \quad     \psi^-_\infty=+\infty \text{  on $\mathcal{B}\backslash A_-$ }  
  $$
  \end{enumerate}
  \el 
  \br\label{rk:trace:infty} Let $\Psi^+_\infty$ be the subgraph of $\psi^{+}_{\infty}$ in $Q$. Following the proof of Lemma  \ref{lm:trace:ca}, the generalized trace of $\psi_\infty^+$ on $\P\Omega$ is  
   \be
  \mg \psi_\infty^{+}(x)=+ \infty \quad \text{ on } \Gamma_+ ;\quad     \mg \psi_\infty^{+}(x)=- \infty \text{  on   }  \P\Omega\backslash \Gamma_+
  \ene
  Thus the trace of $\lambda_{\Psi^+_\infty}$ on $\P\Omega\PLH\R$ in $Q$ is the {\fc} of the subgraph of  $\mg\psi_\infty^+$ over $\P\Omega$. \\
  \indent By Lemma \ref{lm:conclass} the manifold $Q$ can be replaced with $Q_\alpha$.  A similar conclusion also holds for $\psi^-_\infty$. 
  \er 
  \bp We only show the conclusion (1). The proof of the  conclusioin (2) is similar to that of the conclusion (1). So we skip it here.\\
  \indent  Suppose $u$ is a generalized solution on $\Omega$ with boundary data $\psi(x)$ on $\P\Omega$. By Lemma \ref{lm:tran}, $u(x)-j$ is a generalized solution of the Dirichlet problem with boundary $\psi(x)-j$.\\
  \indent  Let $U_j$, $U^+_\infty$, $\Psi_j$ and $\Psi^+_\infty$ be the subgraph of $u(x)-j$, $u^+_\infty(x)$, $\psi(x)-j$ and $\psi^+_\infty$ respectively. \\
  \indent Let $V$ be a {\Ca} set satisfying $V\Delta  U^+_\infty\subset K$ where $K$ is a compact set in $\bar{\Omega}\PLH[-T,T]$. Here $T$ is a positive constant chosen such that 
  \begin{gather}
  \int_{\Omega\PLH\{\pm T\}}d|D\lambda_{U_j}|_{Q_\alpha}= \int_{\Omega\PLH\{\pm T\}}d|D\lambda_{U^+_\infty}|_{Q_\alpha}=0\label{eq:assumption:one}\\
\lim_{j\rightarrow+\infty}  \lambda_{U_j}(x,T)=\lambda_{U^+_\infty}(x,T)\quad a.e. \quad x\in \Omega\label{eq:assumption:two}
 \end{gather}
 The proof of the above equalities are similar to those of \eqref{eq:one:step} and \eqref{eq:two:step}. 
  Now set \be 
 V_j=\left\{\begin{split}
 	& V\text{ in } \Omega\PLH(-T,T)\\
 	& U_j\text{ outside } \Omega\PLH(-T,T)
 \end{split}\right.
 \ene 
 Since $u(x)-j$ is a generalized solution, by Definition \ref{Def:gs} we have  
 \be\label{daming:eq}
 \int_{\bar{\Omega}\PLH[-T,T]}d|D\lambda_{U_j}|_{Q_\alpha}\leq \int_{\bar{\Omega}\PLH[-T,T]}d|D\lambda_{V_j}|_{Q_\alpha}
 \ene 
 Now with \eqref{eq:assumption:one} and Theorem \ref{thm:bd:trace},   \eqref{daming:eq} gives that 
 \be\label{eq:key:step}
 \begin{split}
 &\int_{\Omega\PLH(-T,T)}d|D\lambda_{U_j}|_{Q_\alpha}+\int_{\P\Omega\PLH(-T,T)}|\mT\lambda_{\Psi_j}-\mT\lambda_{U_j}|d\vf_{Q_\alpha}\\
 & \leq  \int_{\Omega\PLH(-T,T)}d|D\lambda_{V}|_{Q_\alpha}+\int_{\P\Omega\PLH(-T,T)}|\mT\lambda_{\Psi_j}-\mT\lambda_{V}|d\vf_{Q_\alpha} \\
 &+\int_{\Omega\PLH\{\pm T\}}|\mT \lambda_V-\lambda_{U_j}|(x,\pm T)d\vf_{Q_\alpha} 
 \end{split}
\ene
In the last line we use the fact $\mT \lambda_{U_j}=\lambda_{U_j}$ on $\Omega\PLH\{\pm T \}$ by \eqref{eq:assumption:one}.\\
\indent  Let $\mg \psi(x)$ be the generalized trace of $\psi(x)$ on $\P\Omega$ from $\mathcal{B}\backslash\bar{\Omega}$. By its definition Lemma \ref{lm:trace:ca} implies that 
\be \label{eq:lim:st}
\mg \psi (x)=+\infty,  \quad x\in \Gamma_+, \quad \mg \psi (x)<+\infty\quad x\in \P\Omega\backslash \Gamma_+
\ene
\indent Let $\mT \lambda_{\Psi}$ and $\mT \lambda_{\Psi_j}$ be the trace of $\lambda_{\Psi}$ and $\lambda_{\Psi_j}$ on $\P\Omega\PLH\R$. By the definition of $\Psi$ and $\Psi_j$, we have 
$
\lambda_{\Psi}(y,t)=\lambda_{\Psi_j}(y, t-j)
$
for every $y\in \mathcal{B}\backslash\bar{\Omega}$. By Lemma \ref{lm:cbv}, it is clear that $\mT \lambda_{\Psi_j}(x,t-j)=\mT \lambda_{\Psi}(x,t)$ for a.e. $x\in \P\Omega$. 
By Lemma \ref{lm:trace:ca}, $\mT \lambda_\Psi$ is the characteristic function of the subgraph of $\mg \psi$ in $\P\Omega\PLH\R$. Thus $\mT\lambda_{\Psi_j}$ is the {\fc} of the subgraph of $\mg \psi-j$ in $\P\Omega\PLH\R$. By \eqref{eq:lim:st} and Remark \ref{rk:trace:infty} letting $j$ go to $\infty$ we obtain that 
$$
\lim_{j\rightarrow+\infty} (\mg \psi-j)=\mg \psi^+_\infty
$$
This implies that 
\be\label{eq:reason:a}
\lim_{j\rightarrow+\infty}\mT \lambda_{\Psi_j}(p)=\mT\lambda_{\Psi^+_\infty}(p)\quad a.e. p\in \P\Omega\PLH\R
\ene 
where $\mT\lambda_{\Psi^+_\infty}$ is the trace of $\lambda_{\Psi^+_\infty}$ on $\P\Omega\PLH\R$. With a similar derivation we shall obtain on $\P\Omega\PLH\R$
  \be\label{eq:reason:b}
  \lim_{j\rightarrow+\infty}\mT \lambda_{U_j}=\mT\lambda_{U^+_\infty} \quad a.e. 
  \ene 
Applying \eqref{eq:reason:a} ,  \eqref{eq:reason:b}, \eqref{eq:assumption:one} and \eqref{eq:assumption:two} into \eqref{eq:key:step} and taking $j\rightarrow+\infty$, one sees that 
\be
\begin{split}
	&\int_{\Omega\PLH(-T,T)}d|D\lambda_{U^+_\infty}|_{Q_\alpha}+\int_{\P\Omega\PLH(-T,T)}|\mT\lambda_{\Psi^+_\infty}-\mT\lambda_{U^+_\infty}|d\vf_{Q_\alpha}\\
	& \leq  \int_{\Omega\PLH(-T,T)}d|D\lambda_{V}|_{Q_\alpha}+\int_{\P\Omega\PLH(-T,T)}|\mT\lambda_{\Psi^+_\infty}-\mT\lambda_{V}|d\vf_{Q_\alpha} \\
	&+\int_{\Omega\PLH\{\pm T\}}|\mT \lambda_V-\lambda_{U^+_\infty}|(x,\pm T)d\vf_{Q_\alpha} 
\end{split}
\ene
By \eqref{eq:assumption:one} and \eqref{eq:assumption:two} this implies that 
\be 
\int_{\bar{\Omega}\PLH[-T,T]}d|D\lambda_{U^+_\infty}|_{Q_\alpha}\leq \int_{\bar{\Omega}\PLH[-T,T]}d|D\lambda_{V}|_{Q_\alpha}
\ene 
Since $V\Delta U_\infty\subset K$ where $K$ is an arbitrary compact set in $\bar{\Omega}\PLH[-T,T]$, $u^+_\infty$ is a generalized solution to the Dirichlet problem \eqref{main:equation} with boundary data $\psi^+_\infty$. The proof is complete. 
 \ep
\subsection{Almost minimal set}
We shall recall some basic facts on almost minimal sets for later use. The papers of Duzzar-Steffen \cite{DS93}  and Tamanini \cite{Tam82} are our main references.  Although their results are discussed in Euclidean space,  their versions in Riemannian manifolds can be easily obtained following from the technique in Section 36 in \cite{Simon83}. 
\begin{Def} \label{lm:partial:minimal}Let $W$ be an open set in a Riemannian manifold $G$ with $dim G=n+1$. Let $inj_W$ denote the injective radius of $W$ in $G$. Let $E$ be a {\Ca}  set in $W$.  We say $E$ is an almost minimal set in $W$ if  it holds that 
\be\label{eq:amb}
\int_{B_\rho (x)}d|D\lambda_E|_{G}\leq \int_{B_\rho(x)}d|D\lambda_F|_{G}+C\rho^{n+\beta}
\ene 
for every point $x$ in any compact set $A\subset W$ and any {\Ca} set  $F\Delta E\subset B_\rho(x)$. Here $\beta\in [0,1)$ is a given  constant,  $\rho<\min\{inj_W, dist(x,G\backslash W)\}$, $C$ is  a positive constant depending on $W$. \\
\indent The boundary $\P E$ (See Proposition \ref{pro:evc}) is called as an almost minimal boundary.  If $C=0$, $\P E$ is called the minimal boundary and $E$ is a minimal set. 
\end{Def}
\br By Proposition \ref{pro:evc} we always require $E$ satisfying $
   0<H^n(E\cap B_\rho(x))<vol(B_\rho(x)) 
    $       for all $x\in \P E$ and sufficiently small $\rho$. 
\er
An example of almost minimal boundaries is the boundary of smooth domains. Our proof imitates  that of Example A.1 in Appendix A of \cite{Eic09} by Eichmair and applies the fact that a $\mC^2$ boundary has locally bounded mean curvature. 
\bl\label{lm:amb} Let $\Omega$ be an open set in a Riemannian manifold $G$. Suppose $\Gamma\subset \P\Omega$ is a $\mC^2$ connected hypersurface in $G$. For each point $x\in \Gamma$, there exists an open set $W$ near $x$ such that $\Omega\cap W$ is an almost minimal set in $W$. 
\el 
\bp We just give the sketch of the proof here. And the notation in Appendix A in \cite{Eic09} will be used without explanation. \\
\indent Notice that the above result is a local result. Without loss of generality, we can assume $\Omega$ is bounded and $\Gamma=\P\Omega$.  Set  $$s(x)=\left\{\begin{split}
&-dist (x,\P\Omega) \,x\in \bar{\Omega}\\
&+dist(x,\P\Omega)\, x\in N\backslash \bar{\Omega }
\end{split}\right.\quad \Gamma_\mu=\{x\in G, s(x)\in (-\mu, \mu)\}$$
Arguing as Lemma 14.16 in \cite{GT01}, $d(x)$ is smooth on  $\Gamma_{\mu}$ for sufficiently small $\mu$. Taking $\mu$ small enough such that for each $y\in \Gamma_\mu$ there is only one $x\in \P\Omega$ with $d(y)=dist(x,y)$. Thus $Dd$ is a smooth unit vector field on $\Gamma_\mu$. Now we view the boundary of a {\Ca} set as a varifold. For any {\Ca} set $F$
\be\label{eq:inst}
M_{B_r(x)}(\P F)=\int_{B_r(x)}d|D\lambda_F|_{G}
\ene 
Now let $x\in  \P\Omega$ and take $r>0$ such that $B_r(x)\subset \Gamma_\mu$. Here $M$ denotes the mass of the varifold (see Simon \cite{Simon83}).\\
\indent Suppose $F$ is a {\Ca} set such that $F\Delta \Omega$ containing a compact set in $B_r(x)$.  Notice that $\P F=\P \Omega+\P(F\Delta \Omega)$. Let $\omega^*= Ds\llcorner d\vf$ where $dvol$ is the volume form of $G$ and $Ds$ is the gradient of $s(x)$.  Suppose $\phi$ is any smooth function with compact support on $B_r(x)$ with $\phi\leq 1$. 
\begin{align*}
M_{B_r(x)}(\P F)&=\sup_{\omega'\in D^n(\Gamma_\mu),|\omega'|\leq 1}(\P \Omega+\P(F\Delta \Omega))(\phi\omega')\\
&\geq \sup_{\phi(x)}(\P\Omega)(\phi\omega')-\sup_{\phi}\P(F\Delta \Omega)(\phi\omega')\\
&\geq M_{B_r(x)}(\P\Omega)-(F\Delta \Omega)(d\omega^*)\\
&\geq M_{B_r(x)}(\P\Omega)-Cvol(B_r(x))
\end{align*}
because $d\omega^*=div(Ds)dvol$ and $|div(Ds)|\leq C$ on $\Gamma_\mu$. With \eqref{eq:inst}, the above estimate is just \eqref{eq:amb}.  Let $W$ be $B_r(x)$. The proof is complete. 
\ep 
Next we define the regular set of a {\Ca} set.
\begin{Def}  Suppose $F$ is a {\Ca} set in a Riemannian manifold $G$. Define the 
	regular set $reg(\P F):=\{x\in \P F: \exists \rho>0 \,\text{with}\, \P F\llcorner B_\rho(x) \text{is  a}\\ C^{1,\beta} \text{ graph}\}$ where $\beta\in (0,1)$. The singular set in $\P F$ is its complement, written as $sing(\P F)$. 
\end{Def}
The following two facts about almost minimal boundaries are standard. 
\bt \label{thm:amb:regularity}  Let $\P F$ be an almost minimal boundary in a Riemannian manifold $G$ with $dim G=n+1$.  
\begin{enumerate} 
	\item (Theorem 1 in \cite{Tam82} and Theorem 5.6 in \cite{DS93}) Then $sing(\P F)=\emptyset$ for $n\leq 6$. If $n = 7$ $sing(\P F)$ consists of isolated points. If $n>7$, $H^{n-7+\beta}(sing \P F)=0$ for $\forall \beta>0$. 
	\item For any compact set $K\subset F$, there exists a $r_0:=r_0(K)>0$ s.t. for all $r\in (0,r_0)$
	$$
	vol(\P F\cap B_r(x))\geq Cr^{n}\quad \forall x\in K \cap \P F
	$$
	where $C$ is a positive constant only depending on $r_0$ and the metric $g$ on $K$.
\end{enumerate}
\et 
The proof of the conclusion (2) above is exactly to that of Proposition 5.14 in \cite{Giu84} if we take $r_0$ as sufficiently small as possible. Thus we skip its proof here. 
\subsection{The property of $P_{\pm }$}First we define the mean curvature of a smooth boundary as follows. 
\begin{Def} Let $W$ be an open domain in a Riemannian manifold $G$ with a $\mC^2$ boundary $\P W$.  Let $\vec{v}$ be the outward normal vector of $\P W$. The mean curvature of $\P W$ is $div(\vec{v})$ written as $H_{\P W}$.  Here $div$ is the divergence of $G$.
\end{Def}
\br  In Appendix \ref{appendix:mc} we collect some facts of mean curvature. 
\er 
The main result of this section is given as follows. 
\bt\label{thm:infty}
Let  $\Omega$ and $\mathfrak{B}$ be two bounded open domains in a Riemanian manifold $N$ satisfying $\Omega\subset \subset \mathcal{B}$ and $\P\Omega$ is $\mC^2$. Let $n$ be the dimension of  $N$. Suppose $u(x)$ is a generalized solution to the Dirichlet problem in \eqref{main:equation} on $\Omega$ with bounded boundary data $\psi(x)$. Define 
\begin{gather*}
P_+=\{x\in\Omega:u(x)=+\infty\}\label{pinfty}\\
P_-=\{x\in\Omega:u(x)=-\infty\}\label{ninfty}
\end{gather*} The following statements hold: \begin{enumerate}
	\item  $P_+\PLH\R$ is an almost minimal set  in $\mathcal{B}\PLH \R\subset Q_\alpha$. Moreover the part of  $P_+\PLH\R$ in  $\Omega\PLH\R
	\subset Q_\alpha$ is an minimal set. The same conclusion holds for $P_-\PLH\R$. 
	\item if in addition $\P\Omega$ is mean convex, that is,
	$H_{\P\Omega}\geq 0$,  \begin{enumerate}
		\item then $\P P_+$ and $\P P_-$ are closed smooth minimal surfaces in $\bar{\Omega}$ when $n\leq 7$; 
		\item then $\P P_+$ and $\P P_-$ are closed smooth minimal surfaces in $\bar{\Omega}$ except a singular set $S$ with Hausdorff measure $H^{k}(S)=0$ for $k>n-7$ when $n>7$. 
	\end{enumerate}
\end{enumerate}
\et
\br  \label{Rk:P:Ex}In Theorem 16.6 \cite{Giu84} $P_+$ is a local minimal set  in $\mathcal{B}$. But in our case we can only obtain that $P_+\PLH\R$ is an local almost minimal set in $Q_\alpha$. This difference is the essential reason that we need Condition (2) to guarantee generalized solutions are classical. The example in Appendix \ref{Example:Section} says such condition is necessary. 
\er 
\bp We just prove the case of $P_+\PLH\R$ and $\P P_+$. The case of $P_-\PLH\R$ and $\P P_-$ is similar so we skip their proofs here.\\
\textbf{The proof of the conclusion (1):} \\
\indent By Lemma \ref{lm:minimize} $u_\infty^+$ is a generalized solution to the Dirichlet problem of \eqref{main:equation}  with boundary data $\psi^+_\infty$ on $\P\Omega$. Because $\psi(x)$ is finite, we have $\psi^+_\infty(x)=-\infty$ outside $\bar{\Omega}$. Thus $P_+\PLH\R$ is the subgraph of $u_\infty^+$ with boundary data $\psi_\infty^+$ in $Q_\alpha$. Thus we also write $U_{\infty}$ for $P_+\PLH\R$ in the following proof.\\
\indent Now the discussion is divided into two cases.  The first case is  $p\in \Omega\PLH\R$. There is a $r>0$ with $B_r(p)\subset \Omega\PLH\R$ where $B_r(p)$ is the ball centerred at $p$ with radius $r$ in $Q_\alpha$.  By  Definition \ref{Def:gs}, we have 
\be\label{con:one:st}
\int_{B_r(p)}d|D\lambda_{U_\infty}|_{Q_\alpha}\leq \int_{B_r(p)}d|D\lambda_{F}|_{Q_\alpha}
\ene 
for any {\Ca} set $F$ satisfying $F\Delta U_\infty\subset B_r(p)$. Thus $P_+\PLH\R$, namely $U_\infty$, is a minimal set in $\Omega\PLH\R$. \\
\indent The second case is $p\in \P\Omega\PLH\R$. Let $E$ be any {\Ca} set coinciding with $U_\infty$ outside some compact set in $B_r(p)$. Since $U_\infty$ locally minimizes the perimeter in $\bar{\Omega}\PLH\R$, then 
\be\label{eq:mi:a}
\int_{B_r(p)}d|D\lambda_{U_\infty}|_{Q_\alpha}\leq \int_{B_r(p)}d|D\lambda_{(\Omega\PLH\R)\cap E}|_{Q_\alpha}
\ene 
 Because $\P\Omega$ is smooth, by Lemma \ref{lm:amb} there is a $r_0>0$  such that $\Omega\PLH\R$ is an almost minimal set in the ball $B_r(p)$ for all $r<r_0$. Thus there are two constants $C=C(r_0)>0$ and $\beta\geq 0$ such that  
 \be \label{eq:mi:b}
 \int_{B_r(p)}d|D\lambda_{(\Omega\PLH\R)}|_{Q_\alpha}\leq \int_{B_r(p)}d|D\lambda_{(\Omega\PLH\R)\cup E}|_{Q_\alpha}+C(r_0)r^{n+\beta}
\ene 
By Lemma 15.1 in \cite{Giu84} one sees that 
 \be\label{eq:mi:c}
 \begin{split}
 \int_{B_r(p)}d|D\lambda_{(\Omega\PLH\R)\cup E}|_{Q_\alpha}&+\int_{B_r(p)}d|D\lambda_{(\Omega\PLH\R)\cap E}|_{Q_\alpha}\\
 &\leq  \int_{B_r(p)}d|D\lambda_{(\Omega\PLH\R)}|_{Q_\alpha}+ \int_{B_r(p)}d|D\lambda_{E}|_{Q_\alpha}
 \end{split}
 \ene 
 Combining \eqref{eq:mi:a}, \eqref{eq:mi:b} and \eqref{eq:mi:c} together we obtain 
 \be \label{eq:mi:d}
 \int_{B_r(p)}d|D\lambda_{U_\infty}|_{Q_\alpha}\leq \int_{B_r(p)}d|D\lambda_{ E}|_{Q_\alpha}+C(r_0)r^{n+\beta}
 \ene 
 Thus 
$P_+\PLH\R$, i.e. $U_\infty$, is an almost minimal set in $\mathcal{B}\PLH\R \subset Q_\alpha$. Together with \eqref{con:one:st} we obtain the conclusion (1). \\
 \textbf{The proof of the conclusion (2):} \\
\indent Since $ P_+\PLH\R$  is an almost minimal set in $Q_\alpha$,  Theorem \ref{thm:amb:regularity} implies that 
\begin{enumerate}
	\item [(a)]If $n\leq 6$, $\P P_+\PLH\R$ is regular. It is locally a $\mC^{1,\beta}$ embedded hypersurface in $Q_\alpha$ for some $\beta>0$. 
	\item [(b)] If $n=7$, the singular set of $\P P_+\PLH\R$ is  a collection of isolated points. 
	\item  [(c)]If $n>7$, let $S$ be the singular set of $\P P_+\PLH\R$. Then $H^{k}(S)=0$ for $k>n-7$. 
	\end{enumerate}
\indent The case of $n=7$ should be examined more carfully.  Let $(x_0, r_0)$ be a isolated singular point in $\P P_+\PLH\R$.  Then there is an open set $V$ in $Q_\alpha$ containing $(x_0, r_0)$ such that $V\cap (\P P_+\PLH\R)$ is $\mC^{1,\beta}$ except the point $(x_0,r_0)$. Thus there is a neighborhood $V_0$ of $x_0$ in $N$ such that $\P P_+\cap V_0$ is $\mC^{1,\beta}$ embedded. But this implies that near $(x_0,r_0)$ $\P P_+\PLH\R$ is $\mC^{1,\beta}$   embedded.  Thus all points in $\P P_+\PLH\R$ are regular.  \\
\indent It is clear that $p=(x,r)$ is a regular point in $\P P_+\PLH\R$ if and only if $x$ is a regular point in $\P P_+$. The above derivation yields that
\begin{enumerate}
	\item [(e)] If $n\leq 7$, $\P P_+$ is a $\mC^{1,\beta}$ embedded hypersurface in $N$ for some $\beta>0$. 
	\item [(f)] If $n>7$, let $S$ be the closed singular set of $\P P_+$. Then $H^{k}(S)=0$ for $k>n-7$. 
\end{enumerate}
\indent The last step is to show that the regular part of $\P P_+$ is smooth and minimal. This is equivalent to show $H_{\P P_+}=0$ a.e. on the regular part of $\P P_+$. \\
\indent First we collect two facts on the regular part of $\P P_+\PLH\R$. Just noticing that the normal vector of $\P P_+\PLH\R$ is perpendicular to $\P_r$ in $Q$ and $Q_\alpha$ we obtain  
 \begin{gather} 
 H^\alpha_{\P P_+\PLH\R}(p)= e^{-\F{r}{n}}H_{\P P_+\PLH\R}(p)\quad \text{by Lemma \ref{lm:con:ref}}\label{eq:fact:one}\\
 H_{\P P_+\PLH\R}(p)=H_{\P P_+}(x)\label{eq:fact:two}
 \end{gather}
 where $p=(x,r)$ belongs to the regular part of $\P P_+\PLH\R$ and one side of these identies exists. Here $H_{\P P_+\PLH\R}$ $(H^\alpha_{\P P_+\PLH\R})$ is the mean curvature of $\P P_+\PLH\R$ in $Q$ $(Q_\alpha)$. 
 Our discussion is divided into two cases. \\
 \indent The first case is when  $x\in \P P_+\cap \Omega$. In the conclusion (1) $P_+\PLH\R$ is a minimal set in $\Omega\PLH\R$. Thus in the regular part of $\P P_+\PLH\R$ we have 
 $
 H^\alpha_{\P P_+\PLH\R}(p)=0
 $ a.e. Thus by \eqref{eq:fact:two} $H_{\P P_+}(x)=0$ a.e. $x\in \P P_+\cap \Omega$.  \\
 \indent The second case is when $x\in \P P_+\cap \P\Omega$.  Let $x$ be a regular point in $\P P_+\cap \P \Omega$. By our assumption, $H_{\P \Omega}(x)\geq 0$ with respect to the outward normal vector. On the other hand, $P_+\PLH\R$ locally minimizes the perimeter in $\bar{\Omega}\PLH\R$. Thus $H^\alpha_{\P P_+\PLH\R}\leq 0$ a.e. on the regular part of $(\P P_+\cap \P\Omega)\PLH\R$ with respect to the outward normal vector. By \eqref{eq:fact:one} and \eqref{eq:fact:two}, we have $H_{\P P_+}(x)\leq 0$ a.e. $x\in \P P_+\cap \P\Omega$ with respect to the outward normal vector. Notice that $\P P_+$ is tangent to $\P\Omega$ at every point in $ \P P_+\cap \P \Omega$. Fix one such $x\in \P P_+\cap \P \Omega$. In a neighborhood of $x$,   $\P P_+$ and $\P\Omega$ are two $\mC^{1,\alpha}$ graphs satisfying $H_{\P \Omega} \geq 0$ and $H_{\P P_+}\leq 0$. Applying the weak version of the strong maximum principle (Theorem 8.19 in \cite{GT01}), we obtain that $\P\Omega$ should coincide with $\P P_+$ near $x$ and $H_{\P  P_+}=0$ a.e.. \\
\indent In summary in the regular part of $\P P_+$ we have $H_{\P P_+}=0$ a.e.  Together with the statements in (e) and (f) we obtain the conclusion (2). The proof is complete.  
 \ep  

  \section{Regularity of generalized solutions}
 In this section we study the regularity property of a generalized solution to the Dirichlet problem of \eqref{main:equation}. Then we study the condition in which generalized solutions of the Dirichlet problem of \eqref{main:equation}  is classical. \\
 \indent All of them are summarized as follows. 
 \bt\label{mt:B}  Let $N$ be a Riemannian manifold with $dim N=n$. Suppose  $\Omega$ is a bounded open domain in $N$ with $\mC^2$ boundary satisfying 
 \begin{enumerate} 
 	\item $ H_{\P\Omega}=div(\vec{v})\geq 0$ on $\P\Omega$ where $\vec{v}$ is the outward normal vector of $\P\Omega$ ; 
 	\item \begin{enumerate}
 		\item if $n\leq 7$, no closed embedded minimal hypersurface exists in $\bar{\Omega}$;
 		\item if $n>7$, no closed embedded minimal hypersurface with a closed singular set $S$ with $H^{k}(S)=0$ for $k>n-7$ exists in $\bar{\Omega}$ where $H^k$ is the $k$-dimensional hausdorff measure on $N$;
 	\end{enumerate}
 \end{enumerate}
 Then the Dirichlet problem of \eqref{main:equation} admits a unique solution $u\in \mC^2(\Omega)\cap \mC(\bar{\Omega})$ for any continous function $\psi(x)$ on $\P\Omega$.
 \et 
 \br The condition (2) can be viewed as an obstacle to the sovability of the Dirichlet problem of \eqref{main:equation}. For example see Appendix \ref{Example:Section}. 
 \er
 \subsection{The interior regularity}
  Now we show the interior regularity of locally bounded generalized solutions to the Dirichlet problem in \eqref{main:equation}. Let $\mC^k(A)$ denote the set of $k$-th differential functions on the open set $A$. Let $\mW^{1,1}(A)$ be the Sobolev space on $A$. For more details see Appendix \ref{appendix:sob}.
 \bt \label{thm:ctwo}Fix $\alpha>0$.  Suppose $v(x)$ is a locally bounded generalized solution to the Dirichlet problem in \eqref{main:equation} on $\Omega$ with boundary data $\psi(x)$ in $\P\Omega$. Then $v(x)\in \mC^2(\Omega)$. 
 \et 
 \br  Our proof follows the spirit of Theorem 14.13 in \cite{Giu84}. 
 \er 
 \bp  Let $\Sigma$ be the graph of $u(x)$ in $Q_\alpha$. Because $u(x)$ is locally bounded,  $\Sigma$ locally minimizes perimeter in $\Omega\PLH\R \subset Q_\alpha$. Thus (1) in Theorem \ref{thm:amb:regularity} implies that $\Sigma$ is smooth except a singular closed set $\Sigma'$ satisfying $H^{n-7}(\Sigma')$=0. Let $S$ denote the projection of $\Sigma'$ from $Q_\alpha$ into $\Omega$. Then $S\subset \Omega$ is a closed set with $H^{n-7}(S)=0$. \\
 \indent Set 
 $
\ms=\Omega\backslash S
 $.  Let $\Sigma_{\ms}$ be the part of $\Sigma$ restricted on $L\PLH\R$. Thus $\Sigma_\ms$ is smooth.  We claim that 
 \bl\label{lm:key:mid} $\Sigma_\ms$ is a smooth graph in $Q_\alpha$ with  $v\in \mC^2(\ms)$. 
 \el 
\bp Notice that  $Q_\alpha$ is conformal to $Q$. For the convenience of computations,  we work in $Q$ instead of $Q_\alpha$. Let $\vec{v}$ be the downward normal vector of $\Sigma_\ms$ in $Q$. Since $v(x)$ locally minizes the functional $\AF_\alpha(v,\Omega)$.  Aruging as in Lemma 2.2 of    \cite{Zhou16c}  $\Ta=\la \vec{v},\P_r\ra$ satisfies 
 \be 
 \Delta \Ta+(|A|^2+\bar{R}ic(\vec{v},\vec{v}))\Ta+\la \nb \Ta, \P_r\ra=0 
 \ene 
 where $\bar{R}ic$ is the Ricci curvature of $Q$, $\Delta$ is the Laplacian operator on $\Sigma$, $\nabla$ is the covariant derivative of $\Sigma_\ms$. Thus on each connected component of $\Sigma_\ms$, $\Ta\equiv 0$ or $\Ta>0$.\\
 \indent  Suppose $\Ta\equiv 0$ on a connected component of $\Sigma_\ms$. 
 Fix $p=(x,u(x))\in \Sigma_\ms$ satisfying $\Ta(p)=0$. Then there is a connected open set $V\subset \Sigma_\ms$ containing $p$ satisfying $\Ta\equiv 0$ on $V$. Let $\Gamma$ be the projection of $V$ on $\Omega$. Therefore $H^{n-1}(\Gamma)>0$. If $z\in \Gamma$, the vertical line through $z$ contains a point $q\in \Sigma\backslash \Sigma'$ with $\Ta(q)=0$. If this line does not meet $\Sigma'$, it lies entirely on $\Sigma$. This is impossible because $u$ is locally bounded. Thus $\Gamma \subset S$.  Thus $H^{n-1}(S)>0$. But this gives a contradiction.  Thus $\Ta>0$ on $\Sigma_\ms$.  Since $\Sigma_\ms$ is smooth, then $v(x)\in C^2(\ms)$ where $\ms=\Omega\backslash S$. \ep
 \indent Fix any $x_0\in S$. By Theorem \ref{DP:TMCE} there is a constant $r^*>0$ such that  for every $r\in (0, r^*)$ the Dirichlet problem 
 \be\label{eq:sDp}
 \left\{\begin{split}
 	div(\F{Du}{\Wg}) &=\F{\alpha}{\Wg},\quad &x\in B_r(x_0)\subset \Omega \\
 	u(x)&=\psi(x)\quad &x\in \P B_r(x_0)
 \end{split}\right. 
 \ene
 is uniquely solvable in $\mC^2(B_r(x_0))\cap \mC(\bar{B}_r(x_0))$ for any continuous function $\psi(x)$ on $\P B_r(x_0)$. \\
 \indent Next we show $v(x)\in \mW^{1,1}(\mathcal{B}_{r^*}(x_0))$. Let $V$ be the subgraph of $v(x)$ in $Q_\alpha$.  Lemma \ref{lm:key:mid} says that 
 \begin{align*}
 \int_{\Omega\PLH\R} d|D\lambda_V|_{Q_\alpha}
 &\geq\int_{(\Omega\backslash S)\PLH\R}d|D\lambda_{V|_{Q_\alpha}} \\
 &=\int_{\Omega\backslash S}e^{\alpha v}\sqrt{1+|Dv|^2}d\vf
 \end{align*}
 where $d\vf$ is the volume on $\Omega$.  Since $v(x)$ is locally finite, thus $$\int_{B_{r^*}(x_0)\backslash S} |Dv|d\vf<\infty $$
 Because $H^{n-1}(S)=0$, we have $v(x)\in \mW^{1,1}(\mathcal{B}_{r^*}(x_0))$. \\
 \indent Choose one $r_0\in (0, r^*)$ such that 
  \be \label{trace:assumption}
  \mT v(x)=v(x)\quad x\in \P \mathcal{B}_0
  \ene where 
  $\mathcal{B}_0$ is  the ball $B_{r_0}(x_0)\subset \Omega$ and $\mT v$ is the traces of $v(x)$ from $\mathcal{B}_0$.\\
  \indent  Since the singular closed set $S$  satisfies that $H^{n-7}(S)=0$, from the definition of the Radon measure we can find a sequence  of open sets $\{S_i\}$ in $\Omega$ such that 
 $$
 S_i\subset\subset S_{i+1},\quad i=1,2,\dots, \cap_{i=1}^\infty S_{i}=S
 $$
 with $H^{n-1}(S_i\cap \P \mathcal{B}_0)\rightarrow 0$ as $i\rightarrow +\infty$. Suppose a sequence of smooth functions $\{\psi_i\}_{i=1}^\infty$ on $\P\mathcal{B}_0$ satisfies 
 \be\label{trace:s}
 \psi_i=v(x) \text{\quad in \quad} \P\mathcal{B}_0\backslash S_i, \quad \sup_{\P \mathcal{B}_0}|\psi_i|\leq 2\sup_{\P \mathcal{B}_0}|v(x)| 
 \ene 
 Again we use the fact that $v(x)$ is locally bounded in $\Omega$.\\
 \indent  Let $\{u_i\}_{i=1}^\infty$ be the classical solution of the Dirichlet problem \eqref{eq:sDp} with smooth boundary data $\{\psi_i\}_{i=1}^\infty$ on $\P\mathcal{B}_0$ . By \eqref{estimate:C0} and \eqref{trace:s}, we have the estimate
 \be 
 \max_{\mathcal{B}_0}|u_i(x)|\leq C(r^*, n,\max_{\P \mathcal{B}_0}|\psi_i|)=C(r^*, n, 2\sup_{\P \mathcal{B}_0}|v(x)| )
 \ene 
 Following the derivation in \eqref{final:estimate:a}, we conclude 
 \be \label{cq:est}
 \max_{\mathcal{B}_0}\{|u_i|, |Du_i|, |D^2u_i|\}\leq C(\sup_{\P\mathcal{B}_0}|\psi(x)|, r^*)
 \ene 
 By the Ascoli-Arzela theorem there is a subsequece, still denoted by $u_j$,  will converge uniformly on compact subsets of $\mathcal{B}_0$ to a locally $\mC^2$ function $u(x)$ satisfying 
 \be 
 div(\F{Du}{\omega})=\F{\alpha}{\omega}\quad \text{on}\quad \mathcal{B}_0
 \ene 
 By \eqref{cq:est}, $u(x)\in \mW^{1,1}(\mathcal{B}_0)$.  Let $\mT:\mW^{1,1}(\mathcal{B}_0)\rightarrow L^1(\P\mathcal{B}_0)$ denote the trace operator. Notice that $u_i\in \mW^{1,1}(\mathcal{B}_0)$ and $u_i$ also converges to $u(x)$ in $\mW^{1,1}(\mathcal{B}_0)$ as $i\rightarrow +\infty$ by the Dominimated Convergence Theorem. Thus 
 $$
 \mT u_i=\psi_i\rightarrow \mT u \quad \text{as}\quad i \rightarrow +\infty \text{  in  } L^1(\P \mathcal{B}_0)
 $$
 On the other hand, from the definition in \eqref{trace:s} and $H^{n-1}(S_i\cap \P \mathcal{B}_0)\rightarrow 0$, $\mT u_i$ converges to $v(x)$ in $L^1(\P \mathcal{B}_0)$. Hence $\mT u(x)=\mT v(x)$ on $\P\mathcal{B}_0$  in the $L^1(\P\mathcal{B}_0)$ sense. \\
 \indent Now we claim $u(x)=v(x)$ on $\mathcal{B}_0$ under the condition $\mT u=\mT v$. Since $u\in \mC^2(\mathcal{B}_0)$ satisfying $div(\F{Du}{\omega})=\F{\alpha}{\omega}$ where $\omega=\sqrt{1+|Du|^2}$. Now we define a vector field in $\Omega\PLH\R$ in the product manifold $Q$ (not the conformal product manifold $Q_\alpha$ ) as 
 $$
 X=e^{\alpha r}\F{\P_r-Du}{\omega}
 $$
 Let $div_p$ be the divergence of $Q$. Suppose $\{e_1, \cdots, e_n\}$ is a local orthormal frame on $\Omega$. Let $\bnb$ deote the covariant derivative of $Q$. Then 
 \begin{align*} 
 div_p(X)&=-e^{\alpha r}\la \bnb_{e_i}\F{Du}{\omega},e_i\ra+\alpha\la\P_r,X\ra \\
 &=e^{\alpha r}(-div(\F{Du}{\omega})+\F{\alpha}{\omega})=0
 \end{align*}
 Let $f\in \mC(\bar{\mathcal{B}}_0)\cap \mC^2(\mathcal{B}_0)$. Let $\Sigma_f$ be the graph of $f$ over $\mathcal{B}_0$. Let $U$ be the domain in $\Omega\PLH\R$ enclosed by $\Sigma_f$,  $\Sigma_u$ (the graph of $u(x)$) and $\P \mathcal{B}_0\PLH\R$. Applying the divergence theorem for $X$ in $U$, we have 
 \begin{align}
 0&=\int_{U} div_p(X)d\vf dr\notag\\
 &=\int_{\mathcal{B}_0} e^{\alpha u(x)}\omega d\vf-\int_{\mathcal{B}_0}e^{\alpha f}\la\F{ \P_r-Du}{\omega}, (\P_r-Df)\ra d\vf\notag\\
 &+\int_{\P\mathcal{B}_0}(f-\mT u)\la X,\gamma \ra d\vf_{\P\mathcal{B}_0}\label{identity}
 \end{align}
 where $\gamma$ is the inward normal vector of $\P\mathcal{B}_0$ and $Df$ is the gradient of $f$ in $\Omega$. \\
 \indent Since $\mC(\bar{\mathcal{B}}_0)\cap \mC^2(\mathcal{B}_0)$ is dense in $\mW^{1,1}(\mathcal{B}_0)$, there is a sequence $\{v_j\}_{j=1}^\infty \in \mC(\bar{\mathcal{B}}_0)\cap \mC^2(\mathcal{B}_0)$ such that $v_j$ converges to $v$ in $\mW^{1,1}(
 \mathcal{B})$ such that $\mT v_j\rightarrow \mT v$ as $j\rightarrow\infty$. Recall that  $\mT u=\mT v$. Replacing $f$ with $v_j$ in \eqref{identity} and letting $j\rightarrow \infty$ give that 
 \be \label{equ:st}
 0=\int_{\mathcal{B}_0} e^{\alpha u(x)}\sqrt{1+|Du|^2} d\mu-\int_{\mathcal{B}_0}e^{\alpha v(x)}\la\F{ \P_r-Du}{\omega}, (\P_r-Dv)\ra d\mu
 \ene 
 Define the function 
 \be 
 \tilde{v}=\left\{\begin{split}
 	&u(x)\quad x\in \mathcal{B}_0\\
 	&v(x)\quad \text{  otherwise  }
 \end{split}
 \right.
 \ene 
 By \eqref{trace:assumption}, we have 
 \be 
   \mT^+\tilde{v}=\mT^-\tilde{v}
 \ene 
 where $\mT^+$ ($\mT^-$) is the trace of $\tilde{v}$ from $\mathcal{B}_0$ (outside $\mathcal{B}_0$). Thus combining Lemma \ref{lm:trace:ca} and Theorem \ref{thm:bd:trace} together, we obtain 
 \be \label{boundary:vanish}
 \int_{\P \mathcal{B}_0\PLH\R}d|D\lambda_{\tilde{V}}|_{Q_\alpha}=0
 \ene 
 where $\tilde{V}$ is the subgraph of  $\tilde{v}$.
 On the other hand because $v(x)$ is a generalized solution in $\Omega$, we obtain 
 \be 
 \int_{\bar{\mathcal{B}}_0\PLH\R}d|D\lambda_V|_{Q_\alpha}\leq   \int_{\bar{\mathcal{B}}_0\PLH\R}d|D\lambda_{\tilde{V}}|_{Q_\alpha}
 \ene 
 By \eqref{boundary:vanish}, this implies that 
 \begin{align*}
 \int_{\mathcal{B}_0\PLH \R}d|D\lambda_V|_{Q_\alpha} &\leq \int_{\mathcal{B}_0\PLH\R}d|D\lambda_{\tilde{V}}|Q_{\alpha}
 \end{align*}
Since $v\in W^{1,1}(\mathcal{B}_0)$ and $u\in C^2(\mathcal{B})$, the above inequality is equivalent that 
 \be 
 \int_{\mathcal{B}_0}e^{\alpha v(x)}\sqrt{1+|Dv|^2}d\vf \leq \int_{\mathcal{B}_0}e^{\alpha u(x)}\omega d\vf
 \ene 
 Due to \eqref{equ:st}, the above equality should hold and implies that 
 \be 
 \F{\P_r-Du}{\omega}=\F{\P_r-Dv}{\sqrt{1+|Dv|^2}}\quad a.e.\quad  x\in \mathcal{B}_0
 \ene 
 Recall that $\omega=\sqrt{1+|Du|^2}$. 
 Thus $Du=Dv$ a.e. $x\in \mathcal{B}_0$. Thus $u=v+C$ for a constant $C$ on $\mathcal{B}_0$. The fact $\mT u=\mT v$ on $\P\mathcal{B}_0$ implies that $C=0$.\\
 \indent We obtain that $u=v$ on $\mathcal{B}_0$. Thus $v(x)\in \mC^2(\mathcal{B}_0)$. Recall that $\mathcal{B}_0=B_{r_0}(x_0)$ where $x_0\in S$. Due to the arbitrariness of $x_0$ we conclude that $v(x)\in \mC^2(\Omega)$. 
 \ep 
 \subsection{The proof of Theorem \ref{mt:B}} The first step is to show that the corresponding generalized solution is locally bounded. 
 \bl\label{lm:mid}  Under the assumption of Theorem \ref{mt:B}, there is a locally bounded generalized solution $u(x)$ on $\bar{\Omega}$ with bounded boundary data $\psi(x)$ . 
 \el 
  \br  In this lemma we do not use the fact that $\P \Omega$ is mean convex.
 \er 
 \bp Notice that $\P\Omega$ is $\mC^2$. Then by Lemma \ref{lm:st:extend} and Remark \ref{re:mk} we can extend $\psi(x)$ is a bounded BV function (still written as $\psi(x)$) on a larger bounded open set containing $\Omega$ such that its subgraph is a {\Ca} set in $Q_\alpha$ and its trace on $\P\Omega$ is $\psi(x)$.  By Theorem \ref{mt:A} there is a generalized solution $u(x)$ with the continuous boundary data $\psi(x)$. \\
 \indent Recall that $P_{\pm}$ are the sets 
 $
 \{x\in \bar{\Omega}: u(x)=\pm\infty \}
 $. 
 By Theorem \ref{thm:infty}, $\P P_{\pm}$ are closed embedded minimal surfaces with a closed singular set $S$ satisfying $H^{k}(S)=0$ for $k>n-7$ in $\bar{\Omega}$ if $n>7$ and $\P P_{\pm}$ are embedded minimal surfaces for $n\leq 7$. By the assumption (2) in Theorem \ref{mt:B} $P_{\pm}$ are empty in the sense of perimeter. Namely $P_{\pm}$ is a $H^{n-1}$ measure zero set  where $n$ is the dimension of $N$. \\
 \indent Now assume $P_{+}$ is not empty. Then $u(x)$ is not locally bounded near some point in $\bar{\Omega}
 $. Without loss of generality we can assume that there is a sequence $\{x_j\}$ in $\Omega$ converging to $x_0$ in $\bar{\Omega}$ such that $u(x_j)>j$ as $j\rightarrow +\infty$. Let $z_j=(x_j, u(x_j))$. \\
 \indent  Because $\psi(x)$ is bounded  there is a $R>0$ such
 that the  normal ball $B_{R}(z_j)$ in $Q_\alpha$ does not intersect the graph of $\psi(x)$. The following proposition is useful. 
 \bpo\label{key:pro} It holds that $U\cap B_R(z_j)$ is an almost minimal set in $B_R(z_j)$.
 \epo 
 \bp 
 The definition in Defintion \ref{Def:gs} implies that $U$ has the least perimeter in $\bar{\Omega}\PLH\R \cap B_{R}(z_j)$. Recall that $\P\Omega \PLH\R$ is $\mC^2$ in $Q_\alpha$.  Arguing as in \eqref{con:one:st}, \eqref{eq:mi:a}, \eqref{eq:mi:b} and \eqref{eq:mi:c}, we will obtain \eqref{eq:mi:d} for $U$. Thus we obtain the conclusion. 
 \ep 
 \indent Let $U_j$ be the subgraph of $u_j(x)=u(x)-j$. By Proposition \ref{key:pro}  each $U_j$ is an almost minimal set in $B_{R_1}(z_0)$.\\
 \indent By Conclusion (2)  in Theorem \ref{thm:amb:regularity} we have 
 \be 
 \vf(\P U_j\cap B_r(z_0))> cr^{n}
 \ene 
 for some $c>0$ and all $r< C(R)$ where $C(R)$ is a positive constant depending on $R$. Thus 
 \be 
 \vf(\P U_j\cap B_{2R}(z_0))\geq cr^{n}
 \ene 
 where $z_0=(x_0,0)$ and any $r\in (0, C(R))$. 
 As $j\rightarrow \infty$ $\lambda_{U_j}$ converges weakly to $\lambda_{P_+\PLH\R}$ in the BV function sense of $Q_\alpha$. Thus $\vf((\P P_+\PLH\R)\cap B_{2R}(z_0))\geq cr^n$. This gives a contradiction since $P_+$ is a $H^{n-1}$ measure zero set. Thus $P_+$ is  empty in $\Omega$.  \\
 \indent A similar derivation yields that $P_-$ is also empty in $\Omega$. Thus $u(x)$ is locally bounded. 
 \ep 
 Now we conclude the boundary continuity when $\psi(x)$ is continuous on $\P\Omega$. 
 \bl\label{last:step} Let $u(x)$ be the generalized solution on $\Omega$ with boundary data $\psi(x)$ under the assumption in Theorem \ref{mt:B}. Then $u(x)$ is continuous on $\bar{\Omega}$ and is equal to $\psi(x)$ on $\P\Omega$. 
 \el 

\bp  By the assumption of Theorem \ref{mt:B} $\psi(x)$ is continuous on $\P\Omega$. By Lemma \ref{lm:mid} $u(x)$ is locally bounded on $\bar{\Omega}$. \\
\indent  Suppose  $x_0\in \P\Omega$ and $\lambda=\lim\sup_{x\in \Omega, x\rightarrow x_0} u(x)> \psi(x_0)$. Then there is a sequence $\{x_j\}$ and $\lambda>0$ in $\Omega$ converging to $x_0$ and 
\be\label{eq:mid:qu}
\lim_{j\rightarrow  +\infty}u(x_j)=\lambda> \psi(x_0)
\ene 
 Let $z_0$ be the point $(x_0,\lambda)$ in $Q_\alpha$.   By Lemma \ref{lm:st:extend} we can view $\psi(x)$ is a bounded continuous BV function on a larger open set containing $\Omega$.  There is a $R>0$ such that there is a $R>0$ such
 that the normal ball $B_{R}(z_0)$ in $Q_\alpha$ does not intersect the graph of $\psi(x)$. \\
\indent By the Nashing isometric  embedding theorem,  we can view $Q_\alpha$ as a smooth submanifold in certain $\R^{n+k}$ with the induced metric for some positive integer $k$. Now we blow up $U\cap B_R(z_0)$ in $\R^{n+k}$ as follows: 
\be 
\begin{split}
   U_j &=\{z\in \R^{n+k}: j^{-1}z+z_0\in U\cap B_R(z_0)\}\\
   S_j &=\{z\in \R^{n+k}: j^{-1}z+z_0\in \bar{\Omega}\PLH\R\}
   \end{split}
\ene 
Similar as the derivation in Theorem 37.4 in \cite{Simon83} $U_j$ will converge weakly to a minimal cone $C$ in $T_{z_0}Q_\alpha$ in $\R^{n+k}$ and $S_j$ converges to a closed plane $S$ in $T_{z}Q_\alpha$. Then $C$ is contained in the half-space determined by $S$. By Theorem 15.5 in \cite{Giu84} $C$ is just a closed half-space in $T_zQ_\alpha$. Thus $\P U$  is $\mC^{1,\beta}$ near $z_0$ and can be written as a graph of a $\mC^{1,\beta}$ function $w(x)$ over $\P \Omega\PLH\R$ with $\omega(z_0)=0$. \\
\indent Since $H_{\P\Omega}\geq 0$,  \eqref{eq:fact:one} implies that  
  $$H^\alpha_{\P\Omega\PLH\R}=e^{-\F{r}{n}}H_{\P\Omega}\geq 0$$ with respect to the outward normal vector of $\P\Omega\PLH\R$ in $Q_\alpha$. Because $\P U$ is smooth, we can assume that $\vec{v}'$ be the normal vector of $\P U$ near $z_0$ which points outward to $(\mathcal{B}\backslash \bar{\Omega})\PLH\R$ at $z_0$. The fact that  $U$ locally minimizes the perimeter in $\bar{\Omega}\PLH\R\subset Q_\alpha$ yields that
   \be 
H^\alpha_{\P U}=div_{Q_\alpha}(\vec{v}')\leq 0 
\ene 
near $z_0$. Notice that $\P U$ is tangent to $\P\Omega\PLH\R$ at $z_0$. By the weak version of the strong maximum principle (see Theorem 8.19 in \cite{GT01})  $\P U$ coincides $\P\Omega\PLH\R$ near $z_0$.  This contradicts to the fact $\lambda=\lim\sup_{x\in \Omega, x\rightarrow x_0} u(x)$. Thus we conlcude
  $$
  \lim_{j\rightarrow  \infty}u(x_j)\leq \psi(x_0)
  $$
With a similar argument $ \lim_{j\rightarrow \infty }u(x_j)\geq \psi(x_0)$. Set $u(x)=\lim_{x_j\rightarrow x_0}u(x_j)$ for $x\in \P\Omega$. Then $u(x)$ is continuous until the boundary and $u(x)=\psi(x)$ for each $x\in \P\Omega$.\\
\ep 
 \indent By Theorem \ref{mt:A} there is a generalized solution $u(x)$ with continuous boundary data $\psi(x)$. Theorem \ref{thm:ctwo}  implies $u\in \mC^2(\Omega)$. By Lemma \ref{last:step} $u(x)$ is continuous to the boundary. Namely $u(x)\in \mC(\bar{\Omega})$.  Since $u(x)$ locally minimizes the conformal area functional $\AF_\alpha(u,\Omega)$,  $u(x)$ satisfies $div(\F{Du}{\omega})=\F{\alpha}{\omega}$ on $\Omega$.  The uniqueness of the solution to the Dirichlet problem \ref{main:equation} is obvious. 
 We complete the proof of Theorem \ref{mt:B}. 
 \appendix
 \section{A decomposition result of Radon measures}\label{app:BC}
 Now we consider a decomposition of Radon measures on {\RM}s. The reason we derive it here is that a domain in Riemannian manifold may not be simply connected any more. Many techiniques in Euclidean spaces can not be applied directly. The main references of this section are Chapter 1 in \cite{Simon83} and Section 2.8 in \cite{Fed69}.  \\
 \indent Throughout this section let $N$ be a complete Riemannian manifold with $dim N=n$.  For every point $x\in N$ we denote  the open (closed) embedded normal ball (see Definition \ref{defopen}) centered at $x$ with radius $r$ by $B_r(x)(\bar{B}_r(x))$. 
 \begin{Def} Let $\mathcal{F}$ be a collection of closed normal balls such that the radius of these balls is a bounded set. Let $A$ denote the set of all centers of those balls. We say that $\mathcal{F}$ covers $A$ finely if the infimum of the radius of balls containing every point in $A$ is $0$. 
 	\end{Def}
 The following theorem is a statement of Theorem 2.8.14, Federer \cite{Fed69} in the case of Riemannian manifolds.  
 \bt[\textbf{\textit{{\Bc} Covering Theorem}}]\label{thm:BVT} Let $\Omega\subset N$ be a bounded open set. There is a positive constant $\kappa=\kappa(n,\Omega)$ such that the following property holds. Let $\mathcal{F}$ be a collection of closed embedded normal balls in $\Omega$ with uniformly bounded radius. Let $A$ be the set of all centers of these balls in $\mathcal{F}$. If $\mathcal{F}$ covers $A$ finely, then there are $\kappa$ subcollections $\{\mathcal{F}_i\}_{i=1}^{\kappa}$ of $\mathcal{F}$ such that the balls in each $\mathcal{F}_i$ are pairwise disjoint and $A\subset \cup_{i=1}^{\kappa}\cup_{\bar{B}\in\mathcal{F}_i}\bar{B}$. 
 \et 
A straightforward verfication shows that 
 \begin{cor}\label{cor:BVT} Suppose $\Omega$ is a bounded open set in $N$. Let $\mu$ be a Radon measure on $\Omega$ with $\mu(\Omega)<\infty$. Let $\mathcal{F}$ be a collection of closed normal balls covering $\Omega$ finely. Then there is a countable pairwise disjoint collection of closed normal balls $\{\bar{B}_{r_j}(x_j)\in \mathcal{F}: j=1,\cdots,\infty \}$ with $\mu(\Omega\backslash\cup_{j=1}^\infty \bar{B}_{r_j}(x_j))=0$. 
 \end{cor}
Now we obtain a useful decomposition of Randon measures in Riemannian manifolds as follows. 
 \bt \label{thm:good:decomposition}
 Let $\Omega$ be an open bound set in a Riemannian manifold. Fix any $\Sc>0$ and $r_0>0$. Suppose $\mu$ is a Radon measure satisfying $\mu(\Omega)<\infty$. Then there is a collection of countable open normal  balls in $\Omega$ defined by 
 \be 
 \mathcal{B}=\{B_k=B_{r_k}(x_k):k=1,\cdots,\infty, x_k\in \Omega, r_k\leq r_0, \mu(\P B_k)=0\}
 \ene 
 and an positive integer $\kappa_0=\kappa_0(\Sc, n,\Omega)$ such that 
 \begin{enumerate}
 	\item $\{B_1,\cdots, B_{\kappa_0}\}$ is a pairwise disjoint subcollection of $\mathcal{B}$ with 
 	$$
 	\mu(\Omega)-\Sc\leq \sum_{k=1}^{\kappa_0}\mu(B_k)=\mu(\cup_{k=1}^{\kappa_0}B_k)\leq \mu(\Omega)
 	$$
 	\item the subcollection $\{B_k:k=\kappa_0+1,\cdots,\infty\}$ of $\mathcal{B}$ satisfies that 
 	$$
 	\sum_{k=\kappa_0+1}^\infty \mu(B_k)\leq \kappa \Sc
 	$$
 	where $\kappa=\kappa(n,\Omega)$ is the positive integer given in Theorem \ref{thm:BVT}.  
 \end{enumerate}
 \et
 \bp Let $d$ be the distance given by the metric on $\Omega$. We define a collection of closed normal balls as follows.
 \be \label{cover:F}
 \mathcal{F}=\{\bar{B}_r(x): x\in \Omega, r<\min\{r_0,d(x,\P\Omega)\}, \mu(\P\bar{B}_r(x))=0\}
 \ene 
 Since $\mu(\Omega)<\infty$,  the Fubini's theorem implies that $ \mu(\P\bar{B}_r(x))=0$ for any $x\in \Omega$ and a.e. $r\in (0, \min\{r_0,d(x,\P\Omega)\})$ . Thus $\mathcal{F}$ covers $\Omega$ finely. \\
 \indent Fix $\Sc>0$. By Corollary \ref{cor:BVT} there is  $\kappa_0=\kappa_0(\Sc,n,\Omega)$ and a pairwise disjoint subcollection of closed balls $\{\bar{B}_{r_i}(x_i)\}$ in $\mathcal{F}$ such that 
 \be \label{est:finite:open}
 \mu(\Omega)-\F{\Sc}{4}\leq \sum_{k=1}^{\kappa_0}\mu(B_{r_k}(x_k))=\mu(\cup_{k=1}^{\kappa_0}B_{r_k}(x_k))\leq \mu(\Omega)
 \ene   
 because $\mu(\P\bar{B}_{r}(x))=0$ for each $\bar{B}_{r}(x)\in \mathcal{F}$. \\
 \indent Namely there is a pairwise disjoint collection of finite open balls 
 \be \label{f:open:ball}
 \{B_{r_1}(x_1),\cdots, B_{r_{\kappa_0}}(x_{\kappa_0})\}
 \ene
 satisfying 
 \be\label{second:property}
 \mu(\Omega\backslash \cup_{k=1}^{\kappa_0}\bar{B}_{r_k}(x_k))\leq \F{\Sc}{4}
 \ene   Now define an open set $\Omega_\eta$ as 
 \be
 \Omega_\eta:\equiv \{x\in \Omega:d(x,\Omega\backslash\cup_{k=1}^{\kappa_0}\bar{B}_{r_k}(x_k))<\eta\}
 \ene 
 where $\eta$ is a sufficiently small positive constant such that $\mu(\Omega_\eta)\leq\F{\Sc}{2}$. Similar as in \eqref{cover:F} we define a collection of closed normal balls in $\Omega_\eta$ as 
 \be 
 \mathcal{F}_\eta=\{\bar{B}_r(x): x\in \Omega_{\eta}, r<\min\{r_0,d(x,\P\Omega_\eta)\}, \mu(\P\bar{B}_r(x))=0\}
 \ene   
 By Theorem \ref{thm:BVT} there are $\kappa=\kappa(\Omega, n)$ subcollections  $\{\mathcal{F}_{\eta,k}\}_{k=1}^{\kappa}$ such that the closed balls in each $\mathcal{F}_{\eta,k}$ are pairwise disjoint and 
 \be 
 \Omega_{\eta}\subset\cup_{k=1}^{\kappa}\cup_{\bar{B}_r(x)\in \mathcal{F}_{\eta,k}}\bar{B}_r(x)
 \ene 
 Moreover for each $k=1,\cdots,\kappa$ \be\sum_{\bar{B}_r(x)\in \mathcal{F}_{\eta,k}}\mu(\bar{B}_r(x))\leq\mu(\Omega_\eta) \leq \F{\Sc}{2}\ene 
 Notice that there are only countable closed normal balls in each subcollection $\mathcal{F}_{\eta,k}$. For each ball $\bar{B}_r(x)$ in each collection $\mathcal{F}_{\eta,k}$ we can replace it with a large open ball $B_{r_x}(x)$ with $r_x<\min\{r_0,d(x,\P\Omega_\eta), 1.5r\}$ such that 
 \be\label{condtion:collection}
 \sum_{\bar{B}_r(x)\in \mathcal{F}_{\eta,k}}\mu(\bar{B}_r(x))\leq\sum_{\bar{B}_r(x)\in \mathcal{F}_{\eta,k}}\mu(B_{r_x}(x))\leq  \Sc
 \ene 
 This gives $\kappa$ collections of open normal balls as follows 
 \be 
 \mathcal{F}'_{\eta,k}:=\{B_{r_x}(x): \bar{B}_{r}(x)\in \mathcal{F}_{\eta,k}, \bar{B}_{r}(x)\subset B_{r_x}(x)\subset\Omega_\eta\}
 \ene
 with the condition \eqref{condtion:collection} for  $k=1,\cdots, \kappa$. Now we relabel all open balls in $\{\mathcal{F}_{\eta,k}\}_{k=1}^{\kappa}$ and list them as follows
 \be \label{second:ball}
 \begin{split}
 	\{B_{r_k}(x_k):k&=\kappa_0+1,\cdots,\infty\}\\
 	&=\{B_{r_x}(x):B_{r_x}(x)\in \mathcal{F}'_{\eta,k},k=1,\cdots,\kappa\}
 \end{split}
 \ene 
 Obviously $\Omega_\eta\subset\cup_{k=\kappa_0+1}^{\infty}B_{r_k}(x_k)$ according to our defintion. The condition \eqref{condtion:collection} yields that 
 \be\label{property:two}
 \sum^\infty_{k={\kappa_0+1}}\mu(B_{r_k}(x_k))\leq \kappa\Sc
 \ene 
 Combining the open balls in \eqref{f:open:ball} with the property \eqref{second:property} and the open balls in \eqref{second:ball} with the property \eqref{property:two} together we obtain the desirable collection of open normal balls. The proof is complete.  
 \ep 
    \section{Local exsitence of the Dirichlet problem of mean curvature equations}\label{appendix:sob}
    In this section we show the existence of the Dirichlet problem \eqref{main:equation} on a sufficiently small normal ball in Riemannian manifolds. \\
    \indent  Throughout this section $N$ is a Riemannian manifold with dimension $n$. For a  set $A$ in $N$, $\mC(A),\mC^k(A), \mC^\infty(A)$ denote the sets of continuous functions, $k$-differential functions and smooth functions on $A$ respectively. For a point $x\in N$ $B_r(x)$ denotes the normal ball centerred at $x$ with radius $r$.\\
    \indent The main result of this section is described as follows. 
    \bt\label{DP:MCE} Fix a point $x_0$ in
    $N$ and  a  function $f(x)\in \mC(N)$. Then there is a positive constant $r^*=r^*(f, Ric,n)>0$ depending on the Ricci curvature near $x_0$ such that for every $r\in (0,r*)$ the Dirichlet problem 
    \be \label{name:eq} 
    	\left\{\begin{split}
    div(\F{Du}{\omega}) &= f(x),\quad &x\in B_{r}(x_0)\\
    u(x)&=\psi(x)\quad &x\in  \P B_{r}(x_0)
    \end{split}
\right.
    \ene 
    is uniquely solvable in $\mC^2(B_r(x_0))\cap \mC(\bar{B}_r(x_0))$ for any $\psi(x)$ in $\mC (\P B_{r}(x_0))$. 
     \et 
    In the following we give some basic definitions, which is from Section 2.1 of Chapter 2 in Hebey \cite{Heb96}.\\
    \indent  Let $g$ be the Riemannian metric of $N$. Let $k$ be an integer and $u\in \mC^\infty(N)$ and $\nb^k$ denote the kth covariant derivative of $u$. By definition one has that 
     $$
     |\nb^k u|^2=g^{i_1j_1}\cdots g^{i_kj_k}(\nb^k u)_{i_1\cdots i_k}(\nb^ku)_{j_1\cdots j_k}
     $$
     where $g=g_{ij}dx^i dx^j$ and $(g^{ij})=(g_{ij})^{-1}$. 
     The $L^p$ norm of $u(x)$ is 
              $$
                ||u||_p=\int_{N}|u|^pd\vf
              $$
 In local coordiantes $d\vf=\sqrt{det(g_{ij})}dx$ is the volume form of $N$ where $dx$ stands for the Lebesgue's volume element of Euclidean space $\R^n$.  The $L^p$ space is the completion of $\mC^\infty(N)$ with respect to the norm $||.||_p$.\\
\indent For $k$ an integer and $p\geq 1$ real, we denote by $\mathcal{C}_k^p(N)$ the sets given by
    $$
     \mathcal{C}_k^p(N)=\{u\in \mC^k(N): s.t. \quad \forall j=0,\cdots,k,\int_N|\nb^j u|^p d\vf<\infty\}
    $$
   \begin{Def}\label{def:sob:space} The Sobolev space $\mW^{k,p}(N)$ is the completion of $\mC_k^p(N)$ with respect to the norm 
      $$
      ||u||_{\mW^{k,p }}=\sum_{j=0}^k(\int_N|\nb^j u|^p d\vf)^{\F{1}{p}}
      $$
    The Sobolev space $\mW_0^{k,p}(N)$ is the completion of the functions in $\mC^p_k(N)$ having compact supports with respect to the norm $||.||_{\mW^{k,p}}$. 
   	\end{Def}
   \indent  First we see the Sobolev inequalities from Euclidean spaces are valid on small normal balls. We denote the constant $C$ only depending on $i,j,k$ by $C(i,j,k)$.
     \bt\label{thm:si:3} Fix $x_0\in N$. The dimension of $N$ is $n$. There is a $r_1>0$ such that the following Sobolev inequalities hold for any $u\in  \mW^{1,p}_0 (B_r(x_0))$ with every $r<r_1$, 
     \begin{gather}
     ||u||_{\F{np}{n-p}}\leq C(n,p,r_1)||Du||_p \text{\quad  for }\quad 1\leq p<n, \label{si:1}\\
     ||u||_p\leq C(r_1) r||Du||_p  \text{\quad  for }\quad 1\leq p<\infty,\label{si:2}
     \end{gather}
     Here $r_1$ depends on the metric near $x_0$. 
         \et 
       
         \br Let $\Omega$ be an open set in $\R^n$. Suppose $u\in W_0^{1,p}(\Omega), 1\leq p <\infty$. According to (7.44) and Theorem 7.10 of Chapter 7 in \cite{GT01}, we have 
          \begin{gather}
          ||u||_{\F{np}{n-p}}\leq C(n,p)||Du||_p \text{\quad  for }\quad 1\leq p<n, \label{si:3}\\
          ||u||_p\leq (\F{1}{\omega_n}|\Omega|)^{\F{1}{n}}||Du||_p  \text{\quad  for }\quad 1\leq p<\infty,\label{si:4}
          \end{gather}
          where $|\Omega|$ is the volume of $\Omega$ and $\omega_n$ is the area of the unit sphere in $\R^n$. 
         \er 
    \bp Fix $x_0\in N$.  First we take $r$ sufficiently. The metric on $B_r(x_0)$ can be written as 
    \be 
    g=g_{ij}(y)dy^idy^j
    \ene  
    Here $y\in B_r(0)$ , in the Euclidean ball in $\R^n$ with radius $r$. By Lemma 1.4 of Chapter 1 in \cite{Heb96} by Hebey there is a constant $r_1$ such that for all $r\in (0,r_1)$ it holds that 
      \begin{gather}
      \F{1}{4}\delta_{ij}\leq g_{ij}(y)\leq 4\delta_{ij}\label{s:7}
      \end{gather}
      \indent Since $\mC_k^p(B_{r_1}(x_0))$ is complete in $W_0^{1,p}(B_{r_1}(x_0))$, we can $u\in \mC_k^p(B_{r_1}(x_0))$. Without confusion we also write $u$ for the corresponding function on the Euclidean ball $B_r(0)$. Let $D_Eu$ denote the Euclidean gradient of $u$, that is, $|D_Eu|^2=\sum_{i=1}^n(\P_i u)^2$. According to \eqref{s:7}  we have the estimate
       \begin{gather}
     \F{1}{4} |D_E u|^2(y)\leq |D u|^2=g^{ij}\P_i u\P_j u\leq 4 |D_E u|^2(y)\label{s:8}\\
     \F{1}{2^n}\leq  \sqrt{det(g_{ij})}(y) \leq 2^n, \label{s:9}
     \end{gather}
       for $y\in B_{r_0}(0)$. \\
       \indent Let $dx$ be the Lebesgue volume form in $\R^n$.  Suppose $1\leq p<n$. With \eqref{si:3}, \eqref{s:8} and \eqref{s:9} we derive 
       \begin{align*}
       ||u||_{\F{np}{n-p}}&=(\int_{B_{r_0}(0)}|u|^{\F{np}{n-p}}\sqrt{\det(g_{ij})}|dx)^\F{n-p}{np}\\
               &\leq C(n,p)||Du||_{\F{np}{n-p}}
       \end{align*}
   We obtain \eqref{si:1}.  Now assume $1\leq p<\infty$. Applying \eqref{si:4} a similar derivation as above should show \eqref{si:2} because $\vf(B_r(x_0))=Cr^n$.
    \ep 
    Next we obtain a $\mC^0$  estimate for the solution to the Dirichlet problem  \eqref{name:eq} in sufficiently small normal balls. 
    \begin{pro}\label{b:five} Fix $x_0\in N $ and $f(x)\in \mC(N)$. Let $r_1$ be the constant given in Theorem \ref{thm:si:3}.  Let $\mu_1$ denote the constant $\max_{B_{r_1}(x_0)}|f|$. Then there is a $r_2=r_2(\mu_1,r_1)<r_1$ such that the following properties holds.  For every $r\in (0, r_2)$ if $u(x)\in \mC(\bar{B}_r(x_0))\cap \mC^2(B_r(x_0))$ is the solution to the Dirichlet problem to \eqref{name:eq} on $B_r(x_0)$, then 
      \be\label{uf:estimate}
      \sup_{B_r(x_0)}|u|\leq C(n,r_1,\mu_1)(\sup_{\P B_r(x_0)}|\psi(x)|+1)
      \ene
     \end{pro}
 \bp Our proof follows from Theorem 10.10 of Chapter 10 in \cite{GT01} .  \\
 \indent For any function $h(x)$, let $h^+(x)=\max\{h(x),0\}$ and $h^-(x)=\max\{-h,0\}$. First we claim that  there is a constant $r^*(\mu_1,r_1)<r_1$ such that for every $r\in (0, r^*)$
 \begin{gather}
    ||u^+||_1\leq C(r_1,\mu_1)(\sup_{\P B_r(x_0)}\psi^+(x)+1)\label{eq:f}
 \end{gather}
 where $u$ is the solution of \eqref{name:eq}. Here and in the followings all $||.||_p$ are computed for measurable functions on $B_r(x_0)$.\\
 \indent It is sufficient to show  \eqref{eq:f}. Without loss of generality, we assume that $r<r_1$ and $\psi(x)\leq 0$ on $\P B_r(x_0)$.  Since $u^+\in \mW_0^1(B_r(x_0))$. Recall that $\omega=\sqrt{1+|Du|^2}$. Let $d\vf$ be the volume form of $(N,g)$. Then by \eqref{name:eq}
 \be\label{mid}
  \int_{B_r(x_0)}\la \F{Du}{\omega}, Du^+\ra d\vf  \leq \int_{B_r(x_0)}|f(x)u^+|d\vf\leq \mu_1||u^+||_1
  \ene 
  With the sobolev inequality \eqref{si:2} this implies that 
$$  
 ||u^+||_1\leq C(r_1)r ||Du^+||_1\leq C(r_1)r\mu_1||u^+||_1+C(r_1)r_1 vol(B_{r_1}(x_0))
$$
Now we take $r_2=r_2(r_1,\mu_1)<r_1$ such that $C(r_1)r\mu_1\leq \F{1}{2}$. Then for all $r<r_2$ we conclude \eqref{eq:f} for $\psi(x)\leq 0$. \\
\indent Now for every $r\in (0,r_2) $ we claim that 
 \be\label{u:estimate}
 \sup_{B_r(x_0)}u^+\leq C(n,r_1,\mu_1)(\sup_{\P B_r(x_0)}\psi^+(x)+1)
\ene 
 First we assume $\psi^+(x)=0$. Then $u^+(x)\in \mW_0^{1,p}(B_r(x_0))$ for any $p$. \\
 \indent We set $v(x)=u^+(x)+1$. Let $\beta>1$ be any fixed constant. Let $v'(x)=v^\beta(x)-1$.   Replacing $u^+$ with $v'(x)$ in \eqref{mid} and applying \eqref{name:eq} we obtain
 $$
  \beta \int_{B_r(x_0)}v^{\beta-1}|Du^+|d\vf \leq \mu_1\int_{B_r(x_0)}(v^{\beta}-1)d\vf+\beta\int_{B_r(x_0)}v^\beta(x)d\vf
  $$
   On the other hand applying \eqref{si:1} on $v^\beta(x)-1$ we obtain 
   \be 
   ||v^\beta(x)-1||_{\F{n}{n-1}}\leq \beta C(n,1) \int_{B_r(x_0)}v^{\beta-1}|Du^+|d\mu
   \ene 
Combining them together  we obtain 
   \begin{align*}
   ||v||_{\beta\F{n}{n-1}}\leq\beta^{\F{1}{\beta}}C(n,\mu_1,r_1)^\F{1}{\beta}||v||_{\beta}
   \end{align*}
  Repeating the above process with $v^{\beta^k}$ where $k\in Z^+$ and letting $k\rightarrow \infty$ gives the estimate in \eqref{u:estimate} with $\sup_{\P B_r(x_0)}u^+=0$. For exact details we refer to Theorem 8.15 in \cite{GT01}. For the general case of $\psi^+$, we replace $u$ with $u-L$ with $L=\sup_{\P B_r(x)}\psi^+$.  \\
  \indent With a similar derivation as above we  can obtain that 
   \be\label{u:estimate:two}
  \sup_{B_r(x_0)}u^-\leq C(n,r_1,\mu_1)(\sup_{\P B_r(x_0)}\psi^-(x)+1)
  \ene 
  Notice that $r_2$ is the same contant as above. Combining \eqref{u:estimate} with \eqref{u:estimate:two} we conclude \eqref{uf:estimate}. The proof is complete. 
  \ep 
  The following proposition is a corollary of Theorem 15.53 in \cite{Aub82} by Aubin. 
  \begin{pro} \label{B:six}Fix $x_0\in N$ and $f\in \mC (N)$. There is a constant $r_3=r_3(n, Ric)$ such that for  $r\in (0, r_3)$ the ball $B_r(x_0)\subset N$ is embedded and  the mean curvature of $\P B_r(x_0)$ with respect to the outward normal vector satisfies 
  	\be \label{eq:st:4}
  	  H_{\P B_r(x_0)}> \max_{B_{r_3}(x_0)}|f(x)|
  	\ene 
  	\end{pro}
  \bp Theorem 15.53 in the Aubin's book \cite{Aub82} implies that $H_{\P B_r(x_0)}$ satisfies 
  \be 
  H_{\P B_r(x_0)}=\F{n}{r}+O(r), \text{\quad  as \quad}  r\rightarrow 0
  \ene 
  when $r$ is sufficiently small. 
 Notice $f(x)$ is a continous function on $N$.  Taking $r_3$ be sufficientlly small shall imply the conclusion. 
  \ep 
  Recall that $\omega=\sqrt{1+|Du|^2}$. Let $\{\P_i\}^n_{i=1}$ be a local coordinate vector field on $B_r(x_0)$.  Denote $\P_i u$ by $u_i$ and $g^{ij}u_i$ by $u^i$.  We define \be\label{def:g}
  \tilde{g}^{ij}=\sigma^{ij}-\F{u^iu^j}{\omega^2}\ene  For a $\mC^2$ function $u$ in $N$, $u_{ij}$ is the covariant derivative $\nb_i\nb_ju$. In fact 
  $$
  div(\F{Du}{\omega})=\F{1}{\omega}\tg^{ij}u_{ij}
  $$
 The following lemma is very useful for the gradient estimate of mean curvature equation.  
  \bl   [Lemma 2.8 in \cite{zhou18}] It holds that 
  \be 
     \tg^{ij}\omega_{ij}=(|A|^2+Ric(\F{Du}{\omega}, \F{Du}{\omega}))\omega+\la \F{Du}{\omega}, \nb (H\omega)\ra+\F{2}{\omega}\tg^{ij}\omega_i\omega_j
  \ene 
  where $|A|^2=\F{1}{\omega^2}\tg^{ik}\tg^{jl}u_{kj}u_{lj}$, $Ric$ is the Ricci curvature of $N$ and $H=div(\F{Du}{\omega})$. 
  \el 
  Let us see a result on the interior estimate of the gradient in terms of boundary estimates of the gradient and $\mC^0$ estimate on mean curvature equations. 
  \bt \label{estimate:gradient} Fix $x_0\in N$. Suppose $u(x)\in \mC^2(B_r(x_0))\cap \mC^1(\bar{B}_r(x_0))$ satisfying 
  $$
  \tg^{ij}\omega_{ij}=\psi(x,u,Du)
  $$ 
  where $\psi(x,u,Du)=1$ or $f(x)\omega$. Then there is a positive constant $K>0$ such that 
  \be \label{conclusion:B}
\sup_{\bar{B}_r(x_0)}|Du| \leq e^{2K \max_{\bar{B}_r(x_0)}|u|}\sup_{\P B_r(x_0)}(1+|Du|)
  \ene  
  Here $K$ is a constant depending on the Ricci curvature on $B_{r}(x_0)$ (and $\max_{\bar{B}_r(x_0)}|f(x)|$ in the case of $\psi(x,u,Du)=f(x)\omega)$)
      \et 
    \bp Let $K$ be a fixed positive constant determined later. Let $\eta$ be the function $e^{K u}$.  Suppose the function $\eta\omega$ achieves the maximum in $\Omega$ at some point $y_0\in B_r(x_0)$. Let $L$ denote the operator $\tg^{ij}\nb_i\nb_j$.  First observe that $\omega_i \eta+\omega\eta_i=0$ at $y_0$. Moreover 
    \begin{align}
    0\geq L(\omega\eta)(y_0)&= (L\omega-\F{2}{\omega}\tg^{ij}\omega_i\omega_j)\eta+\omega L\eta\notag\\
    &\geq \eta(\F{L\eta}{\eta}+Ric(\F{Du}{\omega},\F{Du}{\omega})+|A|^2+\la \F{Du}{\omega^2}, \nb(\psi(x,u,Du))\ra \label{st:ued}
    \end{align}
 We can assume $\omega(y_0)\geq \sqrt{2}$. Otherwise nothing needs to be proven. Thus
 $$
 \tg^{ij}u_iu_j (y_0)=(\sigma^{ij}-\F{u^iu^j}{\omega^2})u_iu_j=\F{|Du|^2}{\omega^2}\geq \F{1}{2}
 $$ Notice that at $y_0$ we have 
   \begin{align}
   \F{L\eta}{\eta}&=K^2\tg^{ij}u^iu^j+ K \psi(x,u,Du)\notag\\
                 &\geq \F{1}{2}K^2+ K \psi(x,u,Du)
   \end{align}
  Since $D\omega= -D\eta=-K Du$ at $y_0$.  No matter in the case of $\psi(x,u,Du)=1$ or $f(x)\omega$ it is clear that at $y_0$. 
  \be 
 K \psi(x,u,Du)+\la \F{Du}{\omega^2},  \nb \psi(x,u,Du)\ra \geq  -C(\mu_1)K
  \ene 
  where $C(\mu_1)$ is a positive constant depending on $\mu_1=\max_{\bar{B}_r(x_0)}|f(x)|$. Now \eqref{st:ued} can be simplified as 
$$  
0\geq L(\omega\eta)(y_0)\geq \eta (\F{1}{2}K^2-\max_{\bar{B}_r(x_0),\la X, X\ra\leq 1}|Ric(X, X)|-C(\mu_1)K)
$$
Now we choose choose $K$ sufficiently large such that   \eqref{st:ued} becomes 
\be 
0\geq L(\omega\eta)(y_0)\geq 0
\ene 
Here $K$ depends on the Ricci curvature (and $\max_{\bar{B}_r(x)}|f(x)|$ in the case of $\psi(x,u,Du)=f(x)\omega$).\\
\indent  Since $L$ is an elliptic operator, this leads to a contradiction to the strong maximum principle of elliptic operators. For such fixed  $K$,  then 
 \be 
 e^{Ku}\omega(x)\leq \sup_{\P B_r(x_0)}e^{Ku}\omega\text{\quad for \quad} x \in \bar{B}_r(x_0) 
 \ene 
 It gives the conclusion \eqref{conclusion:B}. The proof is complete. 
\ep		                                 			
  Now we are ready to show Theorem \ref{DP:MCE}. 
  \bp Our proof is exactly the same in the Dirichlet problem of minimal surface equation in Theorem 16.9 in \cite{GT01}.  First let $r^*$ be $\min\{r_3, r_2\}$ where $r_3$ and $r_2$ be the positive constant in Proposition \ref{B:six} and Proposition \ref{b:five} respectively.  It is obvious that $r^*<r_1$ and $r^*=r^*(n, Ric, \mu_1)$ where $r_1$ is the constant in Theorem \ref{thm:si:3}. \\
  \indent Now assume $r\in (0, r^*)$ and $\psi\in \mC^{2}(\P\Omega)$.  If for $u_\sigma\in \mC(\bar{B}_{r}(x))\cap \mC^2(B_r(x))$  is the solution    \be 
  \left\{\begin{split}
  	div(\F{Du}{\omega}) &= \sigma f(x),\quad &x\in B_{r}(x)\\
  	u(x)&=\sigma \psi(x)\quad &x\in \P B_{r}(x)
  \end{split}
  \right.
  \ene 
  for every $\sigma\in [0, 1]$ we claim that 
  \be \label{final:estimate:a}
  \max_{\Omega}\{|u_\sigma|, |Du_\sigma|, |D^2u_\sigma|\}\leq \mu_0=\mu_0(\sup_{\P\Omega}|\psi(x)|, \sup_{B_r(x)}|f(x)|)
  \ene 
Since $r<r^*<r_2$, Proposition \ref{b:five} implies that 
 \be 
     \max_{\Omega}|u_\sigma|\leq C(n,r_1,\mu_1)(\sup_{\P B_r(x)}|\psi(x)|+1)\quad \text{for all $\sigma\in [0, 1]$}
 \ene 
 Notice that $u_\sigma$ satisfies that 
  $$
   \tg^{ij}u_{ij}=f(x)\omega
  $$
 Therefore by Theorem \ref{estimate:gradient}, we have 
 \be 
 \sup_{\bar{B}_r(x)}|Du_\sigma| \leq e^{2K \max_{\bar{B}_r(x)}|u_\sigma|}\sup_{\P B_r(x)}(1+|Du_\sigma|)
 \ene 
 On the other hand arguing as Theorem 14.6 in \cite{GT01}  we conclude 
 \be 
 \sup_{\P B_r(x)} |Du_\sigma|\leq C(n,\mu_1, \sup_{\P B_1(x)}|\psi(x)|)
 \ene 
 because of $r<r^*<r_3$ and \eqref{eq:st:4}. This gives a uniform estimate of the gradient $|Du_\sigma|$ for  $\sigma\in [0, 1]$.  The classical Schauder estimate in Riemannian manifolds gives the estimate of $|D^2u_\sigma|$ in \eqref{final:estimate:a}.\\
 \indent  Applying Theorem 11.3 in \cite{GT01} we obtain the existence of the solution to the Dirichlet problem \eqref{name:eq} assuming $\psi(x)\in \mC^2(\P B_r(x))$ for all $r\in (0, r^*)$. \\
 \indent For general $\psi(x)\in \mC(\P B_r(x))$ since $\P B_{r}(x)$ is smooth,  there is a sequence $\{\psi^{k}(x)\}_{k=1}^\infty$ in $\mC^2(\P B_r(x))$ such that $\psi^k(x)$ converges to $\psi(x)$ uniformly as $k\rightarrow \infty$ on $\P B_r(x)$. With the above argument, there is a sequence $\{u^{k}(x)\}$ as  the existence of the solutions to the Dirichlet problem \eqref{name:eq} with $u^{k}(x)=\psi^k(x)$ on $\P B_r(x)$. Since $\{\psi^k\}_{k=1}^\infty$ and $\psi$ are uniformly bounded, $\{u^k\}$ satisfies \eqref{final:estimate:a} uniformly. By the Arzela-Ascoli theorem  $\{u^k\}$ converges to $u(x)$ in the $\mC^2$ norm. And $u(x)$ is a solution of the Dirichlet problem \eqref{name:eq} with $u(x)=\psi(x)$ for $\psi(x)\in \mC(\P B_r(x))$. \\
 \indent The proof is complete. 
  \ep 
  A direct consequence of Theorem \ref{DP:MCE} is given as follows. It  says that the Dirichelt problem \eqref{main:equation} is solvable if the domain is a sufficiently small normal ball in Riemannian manifolds. This plays a key role in the proof of  Theorem \ref{thm:ctwo}. 
   \bt\label{DP:TMCE} Fix $x_0\in N$ and $\alpha>0$. Then there is a positive constant $r^*>0$ depending on the sectional curvature near $x_0$, $\alpha$ and $n$ such that for every $r\in (0,r*)$ the Dirichlet problem 
  \be \label{name:eq:a} 
  \left\{\begin{split}
  	div(\F{Du}{\omega}) &= \F{\alpha}{\omega},\quad &x\in B_{r}(x_0)\\
  	u(x)&=\psi(x)\quad &x\in \P B_{r}(x_0)
  \end{split}
  \right.
  \ene 
  is uniquely solvable in $\mC^2(B_r(x_0))\cap \mC( \bar{B}_r(x_0))$ for any $\psi(x)$ in $\mC (\P B_{r}(x_0))$. Here $B_r(x_0)$ is the ball centered at $x_0$ with radius $r$ in $N$. 
  \et 
  \bp Fix $x_0\in N$. By Theorem \ref{DP:MCE} there is a constant $r^*=r^* (Ric, n,\alpha)>0$ such that for every $r$ in $(0, r^*)$ the Dirichlet problem
  \be 
  div(\F{Du}{\omega})=\pm \alpha \text{\quad on\quad} B_r(x_0),  \quad u(x)=\psi(x)  \text{\quad on\quad} \P B_r(x_0)
  \ene 
  is solvable in  $\mC^2(B_r(x_0))\cap \mC( \bar{B}_r(x_0))$ for any $\psi(x)$ in $\mC(\P B_r(x_0))$. Now we write $u^+(x)$ ($u^{-}(x)$) for its solution as in the case that the right hand side is $+\alpha$ ($-\alpha$) respectively.\\
  \indent 
  For any $\sigma\in [0,1]$, suppose $u_\sigma\in \mC^2(B_r (x_0))\cap \mC(\bar{B}_r(x_0))$ is the solution of the following Dirichlet problem: 
   \be \label{A:name:eq} 
  \left\{\begin{split}
  	div(\F{Du}{\omega}) &= \sigma\F{\alpha}{\omega}\quad &x\in B_{r}(x_0)\\
  	u(x)&=\psi(x)\quad &x\in \P B_{r}(x_0)
  \end{split}
  \right.
  \ene 
  By the comparision principle, we have  $u^+(x)\leq  u_\sigma\leq u^-(x)$ on $B_r(x)$. Namely we have the estimate 
  \be\label{estimate:C0}
  \max_{\bar{B}_r(x_0)}|u_\sigma(x)|\leq C(r^*, n,\max_{\P B_r(x_0)}|\psi|,\alpha)
   \ene 
  according to Proposition \ref{b:five}. Proceeding with the same derivation in the Theorem \ref{DP:MCE}, we obtain
   \be
  \max_{\Omega}\{|u_\sigma|, |Du_\sigma|,|D^2u_\sigma|\}\leq C(\sup_{\P\Omega}|\psi(x)|, r^*,\alpha)
  \ene 
\indent  Thus 
   the existence of the solution to the Dirichlet problem \eqref{name:eq:a} follows from the continuation method in (Theorem 11.3, \cite{GT01})  if assuming $\psi(x)\in \mC^2(\Omega)$. The general case can be derived as similar as that in Theorem \ref{DP:MCE}.  The proof is complete. 
    \ep 
    \section{An example}\label{Example:Section} 
    We give an example to illustrate the condition (2) in Theorem \ref{mt:B} can not be removed in general. 
    \bt The sphere $S_n$ in $\R^{n+1}$ is given by 
    $$
    S_n=\{(x_1,\cdots,x_{n+1})\in \R^{n+1}: x_1^2+\cdots+x_{n+1}^2=1\}
    $$with the induced metric $\sigma_n$. Let $S_n^+$ be the open upper half sphere given by
    $$
    S_n^+= \{(x_1,\cdots,x_{n+1})\in \R^{n+1}: x_1^2+\cdots+x_{n+1}^2=1, x_{n+1}>0\}
    $$
    Then the Dirichlet problem 
    	\be\label{dt:mce}
  div(\F{Du}{\omega}) =\F{n}{\omega} \text{ on } S_n^+;\quad  u(x)=\psi(x)\text{ on }\P S_n^+
    \ene
    for any continuous function $\psi(x)$ on $\P S_n^+$ admits no solution in $ \mC(\bar{S}_n^+)\cap
    \mC^2(S_n^+)$.  Here $\omega=\sqrt{1+|Du|^2}$. 
    \et 
    \br  Notice that the set 
    $$
    \P S_n^+:=\{(x_1,\cdots, x_{n+1})\in \R^{n+1}: x^2_1+\cdots+x^2_{n+1}=1, x_{n+1}=0
    \}$$
    is a minimal smooth hypersurface in $(S_n,\sigma_n)$. The above theorem satisfies the assumptions of Theorem \ref{mt:B}  except its condition (2). 
    \er 
    \bp Suppose there is a continuous function $\psi(x)$ on $\P S_n^+$ such that the Dirichlet problem \eqref{dt:mce} has a solution $u(x)\in
     \mC(\bar{S}_n^+)\cap \mC^2(S_n^+)$. It is clear that  $\Sigma=(x,u(x))$ is a minimal surface in the conformal product manifold 
     $$
       M:= (S_n\PLH\R, e^{2r} (\sigma_n+dr^2))
     $$
    Let $s=e^r$. Then $M$ is isometric to the Euclidean space $\R^{n+1}\backslash\{0\}$ with its polar coordinate representation
     $$
        (S_n\PLH (0,+\infty), s^2\sigma_{n}+ds^2)
     $$
     via the  map $f(x,r)=(x, e^r)$ from $M$ to $\R^n\backslash\{0\}$. 
     Thus $f(\Sigma)$ is an embedded minimal surface in $\R^{n+1}$ not touching the origin with its boundary on the plane $x_{n+1}=0$ in $\R^{n+1}$. Then $x_{n+1}=0$ on the boundary of $f(\Sigma)$ and it is not a constant function. Therefore there is a point $p=(p_1, \cdots, p_{n+1})$ in the interior of $f(\Sigma)$ such that  $p_{n+1}\neq 0$  attains the minimum or maximum of the $x_{n+1}$-coordinate on $f(\Sigma)$. From the maximum principle then $f(\Sigma)$ should be contained in the totally geodesic plane $x_{n+1}=p_{n+1}$. This gives a contradiction since the boundary of $f(\Sigma)$ belongs to the plane $x_{n+1}=0$. \\
     \indent Thus no classical solution to the Dirichlet problem \eqref{dt:mce} exists.  
    \ep 
    \br\label{Rk:not:gs}
    By Theorem \ref{mt:A} the generalized solution to \eqref{dt:mce} exists. Define $u_0(x)=-\infty$ on $S_n^+$ and $u_0(x)$ is a continuous function on $S_n\backslash S_n^+$ such that $u_0(x)=\psi(x)$ on $\P S_n^+$. From the above proof $u_0(x)$ is a generalized solution to \eqref{dt:mce} with boundary data $\psi(x)$. \\
    \indent Thus $P_+$ is equal to $S_n^+$. Moreover the above proof shows that $\P S_n^+\PLH\R$ has locally least perimeter in $M=\R^{n+1}\backslash \{0\}$. But $S_n^+$ is not a minimal set in $S_n$ because its boundary $\P S^+_{n}$ is unstable in $S^+_n$. 
   \er 
    
    \section{Some facts on mean curvatures}\label{appendix:mc}
    In this section we collect some facts on mean curvature. Our reference is \cite{Zhou16c}.\\
    \indent Let $M$ be a {\RM} with a metric $\sigma$. Let $\Sigma$ be a smooth hypersurface in $M$ with a normal vector $\vec{v}$. Suppose a local orthonormal frame $\{e_1,e_2,\cdots,e_n\}$. Then the mean curvature vector of $\Sigma$, $\vec{H}$, is given by 
    $$
    \vec{H}=\la \bnb_{e_i}e_i,\vec{v}\ra\vec{v}=-div_M(\vec{v})\vec{v}
    $$
    where $\bnb$ and $div_M$ are the covariant derivative and the divergence of $M$ respectively.
    \begin{Def}  $H=div_M(\vec{v})$ is called the mean curvature of $\Sigma$ with respect to the normal vector $\vec{v}$. 
    	\end{Def}
    \indent For example the mean curvature of the standard sphere $S_n$ in $\R^{n+1}$ with the outward normal vector is $n$.
    \bl[Lemma 3.1 in \cite{Zhou16c}]\label{gs:lm} Let $f$ be a smooth function on $M$. Let $\tilde{M}$ denote the same smooth manifold as $M$ equipped with the metric $e^{2f}g$. Let $\Sigma$ be a smooth hypersurface in $M$ (so is in $\tilde{M}$). Then the mean curvatures of $\Sigma$ in $\tilde{M}$ and in $M$ have the relation that
    \be 
    \tilde{H}=e^{-f}(H+(m-1) df(\vec{v}))
    \ene 
    Here $\vec{v}$ is the normal vector of $\Sigma$ in $M$, $m-1=dim \Sigma$, $H$ ($\tilde{H}$) denotes the mean curvature of $\Sigma$ in $M$ $(\tilde{M})$ with respect to the normal vector  $\vec{v}$ ($e^{-f}\vec{v}$).
    \el 
    Let $N$ be a complete manifold with a metric $g$. Let $Q$ be a product manfiold $N\PLH\R$ with the metric $g+dr^2$ and $Q_\alpha$ be the conformal product manifold $N\PLH\R$ equipped with the metric $e^{2\F{\alpha}{n}r}(g+dr^2)$. From Lemma \ref{gs:lm} we conclude that  
     \bl \label{lm:con:ref}Let $\Omega$ be a smooth hypersurface in $Q$ with the normal vector $\vec{v}$. Let $H$ be the mean curvature of $\Sigma$ in $Q$ with respect to $\vec{v}$. Then the mean curvature of $\Sigma$ in $Q_\alpha$, $H^\alpha$, is 
    \be 
    H^\alpha =e^{-\alpha\F{r}{n}}(H+\alpha\la \vec{v},\P_r\ra)
    \ene 
    where $\la,\ra$ is the inner product in $Q$. 
    \el  Suppose  $u\in \mC^2(\Omega)$ and $\Sigma$ is  the graph of $u(x)$. Let $\vec{v}$ be the downward normal vector of $\Sigma$. Then the mean curvature of $\Sigma$ with repsect to $\vec{v}$ in $Q$ is 
    $
    H=div(\F{Du}{\omega})
    $
    where $div$ is the divergence of $N$ and $\omega=\sqrt{1+|Du|^2}$. Lemma 2.2 in \cite{Zhou16c} indicates that 
    \bl If $u\in \mC^2(\Omega)$ satisfies $div(\F{Du}{\omega})=\F{\alpha}{\omega}$, then its subgraph is minimal in $Q_\alpha$. 
    \el 
   
    \bibliographystyle{abbrv}	
    \bibliography{Paper}
	   	\end{document}